%% file: mdof_ir.tex
\documentclass{article}

\usepackage[margin=1in]{geometry}
\usepackage{amsmath}
\usepackage{amssymb}
\usepackage{graphicx}
\usepackage{subcaption}
\usepackage{authblk}
\usepackage[firstinits=true, sorting=none, style=numeric-comp, backend=bibtex, maxnames=999]{biblatex}
\usepackage{hyperref}
\usepackage{cleveref}
\usepackage{xcolor}
\usepackage{multirow}
\usepackage{soul} \usepackage{comment} \usepackage{placeins} 

\newcommand{\removeBib}[1]{\AtEveryBibitem{\clearfield{#1}}\AtEveryBibitem{\clearlist{#1}} \AtEveryBibitem{\clearname{#1}}} 

\removeBib{language}
\removeBib{month}
\removeBib{day}

\DeclareFieldFormat{titlecase}{\MakeTitleCase{#1}}

\newrobustcmd{\MakeTitleCase}[1]{\ifthenelse{\ifcurrentfield{booktitle}\OR\ifcurrentfield{booksubtitle}\OR\ifcurrentfield{maintitle}\OR\ifcurrentfield{mainsubtitle}\OR\ifcurrentfield{journaltitle}\OR\ifcurrentfield{journalsubtitle}\OR\ifcurrentfield{issuetitle}\OR\ifcurrentfield{issuesubtitle}\OR\ifentrytype{book}\OR\ifentrytype{mvbook}\OR\ifentrytype{bookinbook}\OR\ifentrytype{booklet}\OR\ifentrytype{suppbook}\OR\ifentrytype{collection}\OR\ifentrytype{mvcollection}\OR\ifentrytype{suppcollection}\OR\ifentrytype{manual}\OR\ifentrytype{periodical}\OR\ifentrytype{suppperiodical}\OR\ifentrytype{proceedings}\OR\ifentrytype{mvproceedings}\OR\ifentrytype{reference}\OR\ifentrytype{mvreference}\OR\ifentrytype{report}} {#1}
{\MakeSentenceCase{#1}}}

\providecommand{\keywords}[1]{\textbf{\textit{Keywords---}} #1}

\title{Efficient Model Reduction and Prediction of Superharmonic Resonances in Frictional and Hysteretic Systems}
\author[1]{Justin H. Porter~\thanks{jp88@rice.edu}}
\author[1]{Matthew R. W. Brake~\thanks{brake@rice.edu}}

\affil[1]{Department of Mechanical Engineering, Rice University, Houston,
  TX 77005}

\begin{document}

\include{preamble.tex}

\pagebreak

\pagenumbering{arabic} 

\maketitle{}
\begin{abstract}

Modern engineering structures exhibit nonlinear vibration behavior as designs are pushed to reduce weight and energy consumption. Of specific interest here, joints in assembled structures introduce friction, hysteresis, and unilateral contact resulting in nonlinear vibration effects. In many cases, it is impractical to remove jointed connections necessitating, the understanding of these behaviors. This work focuses on superharmonic and internal resonances in hysteretic and jointed systems. Superharmonic resonances occur when a nonlinear system is forced at an integer fraction of a natural frequency resulting in a large (locally maximal) response at an integer multiple of the forcing frequency. When a second vibration mode simultaneously responds in resonance at the forcing frequency, the combined phenomena is termed an internal resonance. First, variable phase resonance nonlinear modes (VPRNM) is extended to track superharmonic resonances in multiple degree of freedom systems exhibiting hysteresis. Then a novel reduced order model based on VPRNM (VPRNM ROM) is proposed to reconstruct frequency response curves faster than utilizing the harmonic balance method (HBM). The VPRNM ROM is demonstrated for a 3 degree of freedom system with a 3:1 internal resonance and for the jointed Half Brake-Reu{\ss} Beam (HBRB), which exhibits a 7:1 internal resonance. For the HBRB, new experimental results are used to validate the modeling approaches, and a previously developed physics-based friction model is further validated, achieving frequency predictions within 3\%. For the considered cases, VPRNM ROM construction is up to 4 times faster than HBM, and the evaluation of the VPRNM ROM is up to 780,000 times faster than HBM. Furthermore, the modeling framework provides insights into the mechanisms of superharmonic resonances in jointed structures, showing that both tangential slipping and normal direction clapping of the joint play important roles in exciting the superharmonic resonances in the HBRB.

\end{abstract}

\section{Introduction}

Modern engineering practice seeks to optimize designs to reduce weight while improving efficiency and reliability. 
Optimized designs can demonstrate nonlinear vibration behavior due to friction in bolted connections \cite{brakeMechanicsJointedStructures2017} and thin panels \cite{bhattuExperimentalAnalysis2023}. 
At the same time, vibration testing to fully understand this nonlinear behavior is costly \cite{brakeMechanicsJointedStructures2017}, and a poor understanding of nonlinear vibration behavior hinders further optimization of designs. 
For jointed structures, nonlinear vibration behavior is generally characterized by a decrease in resonant or modal frequency and an increase in modal damping with an increase in amplitude. 
While not perfect, recent models have shown clear improvements in understanding this behavior \cite{porterPredictivePhysicsbasedFriction2023}, and many computationally efficient modeling approaches are available \cite{krackNonlinearModalAnalysis2015, lacayoUpdatingStructuralModels2019, balajiQuasistaticNonlinearModal2020}.
However, frictional and hysteretic systems exhibit the more complicated nonlinear vibration behaviors of superharmonic and internal resonances as observed both experimentally \cite{chenMeasurementIdentificationNonlinear2022, scheelChallengingExperimentalNonlinear2020, IMAC_LTW, smithInfluenceWearNonlinear2024, claeysModalInteractionsDue2016, bhattuExperimentalAnalysis2023} and numerically \cite{krackMultiharmonicAnalysisDesign2012, krackReducedOrderModeling2013, claeysModalInteractionsDue2016, krackEfficacyFrictionDamping2016}.
Understanding superharmonic and internal resonances is critical to prevent turbomachinery failures because these phenomena can alter maximum vibration amplitudes \cite{krackMultiharmonicAnalysisDesign2012} and be used for identifying crack failures \cite{sinouDetectionCracksRotor2008}.
Internal and superharmonic resonances are also important for MEMS \cite{asadiStrongInternalResonance2021} and wind turbines \cite{liFlapwiseDynamicResponse2012, inoueNonlinearVibrationAnalysis2012, ramakrishnanResonancesForcedMathieu2012}.
In the case of superharmonic and internal resonances, computationally efficient approaches for modeling frictional and hysteretic systems break down, requiring more computationally intensive simulation techniques.

When a linear structure is harmonically excited at a constant frequency, the structure responds in steady-state only at the forcing frequency. 
For nonlinear structures, responses generally include additional harmonics at integer multiples of the forcing frequency \cite{krackHarmonicBalanceNonlinear2019}.
Superharmonic resonances are locally maximal responses of these higher harmonics at integer multiples of the forcing frequency. Generally, superharmonic resonances occur when the forcing frequency is approximately the natural frequency of a mode divided by an integer and are characterized by a phase shift of the superharmonic \cite{volvertPhaseResonanceNonlinear2021}.
Internal resonances occur when two vibration modes respond near resonance with an integer ratio of frequencies. 
For instance, if one mode has a primary resonance near the same forcing frequency as a second mode has a superharmonic resonance then the combined response is called an internal resonance. 
Internal resonances can also occur with subharmonic resonances (resonant responses at the forcing frequency divided by an integer) and ultra-subharmonic resonances (resonant responses with an integer ratio to the forcing frequency) \cite{nayfehNonlinearInteractionsAnalytical2000, volvertPhaseResonanceNonlinear2021}. However, subharmonic and ultra-subharmonic resonances are rarely reported for frictional systems, so are not discussed further here.

For simple systems, superharmonic and internal resonances can be modeled analytically with perturbation methods \cite{nayfehNonlinearInteractionsAnalytical2000, nayfehNonlinearOscillations1995}.
For frictional systems, analytical methods are less frequently applied, but some recent analytical work has considered modal interactions \cite{quinnIdentificationModeCoupling2023}.
Other efficient analysis techniques for superharmonic and internal resonances include invariant manifolds \cite{boivinNonlinearModalAnalysis1995, pesheckNonlinearModalAnalysis2001, jiangConstructionNonlinearNormal2005, liNonlinearAnalysisForced2022, bettiniModelReductionSpectral2024} and Melnikov-type analyses \cite{cenedeseHowConservativeBackbone2020, guckenheimerAveragingPerturbationGeometric1983}.
However, numerically described friction nonlinearities prevent easy applications of these efficient and analytical methods. 
As an alternative, the harmonic balance method (HBM) can be utilized to calculate steady-state harmonic responses for models with many degrees of freedom (DOFs) \cite{krackHarmonicBalanceNonlinear2019}, but HBM is computationally intensive motivating more efficient numerical modeling approaches.

Nonlinear normal modes (NNMs), as defined by \cite{kerschenNonlinearNormalModes2009, peetersNonlinearNormalModes2009}, provide a computationally efficient approach for analyzing isolated nonlinear modes. The approaches allow for numerical calculation of backbones (amplitude dependent modal frequency and vibration shape) for conservative systems \cite{kerschenNonlinearNormalModes2009, peetersNonlinearNormalModes2009}. 
NNMs have been extended to damped systems by complex nonlinear modes (CNMs) \cite{laxaldeComplexNonlinearModal2009} and the extended periodic motion concept (EPMC) \cite{krackNonlinearModalAnalysis2015} that both calculate amplitude dependent damping in addition to frequency and vibration shape.
EPMC has been used in numerous studies for analyzing the nonlinear trends of frictional systems \cite{jahnComparisonDifferentHarmonic2019, sunNonlinearModalAnalysis2021, sunComparisonDifferentMethodologies2021, mullerNonlinearDampingQuantification2022, yuanComputationDampedNonlinear2022, wuDesignWavelikeDry2021, porterPredictivePhysicsbasedFriction2023, balajiQuasistaticNonlinearModal2020, schwarzValidationTurbineBlade2020}.
Given a nonlinear modal backbone, efficient reduced order models (ROMs) based on single nonlinear mode theory have been developed for constructing transient responses \cite{krackNonlinearModalAnalysis2015, krackComputationSlowDynamics2014} and harmonically forced frequency response curves (FRCs) \cite{krackReducedOrderModeling2013, krackMethodNonlinearModal2013, schwarzValidationTurbineBlade2020} when a single nonlinear mode dominates the response. 
However, when two modes exhibit an internal resonance, EPMC breaks down. Additionally, ROMs based on single nonlinear mode theory break down in the presence of internal resonances \cite{krackReducedOrderModeling2013, krackMethodNonlinearModal2013, schwarzValidationTurbineBlade2020, quaegebeurNonlinearCyclicReduction2021}.
This breakdown motivates recent work that has investigated computationally efficient approaches for understanding superharmonic resonances.

For forced steady-state responses, phase resonance nonlinear modes (PRNM) have been proposed to efficiently track superharmonic, subharmonic, and ultra-subharmonic resonance behavior across forcing levels \cite{volvertPhaseResonanceNonlinear2021}.
PRNM exploits the phase shift of the resonant superharmonic, subharmonic, or ultra-subharmonic to establish a resonant phase criteria analytically \cite{volvertPhaseResonanceNonlinear2021, volvertResonantPhaseLags2022}. This phase lag can then be tracked with a slightly modified HBM approach calculating only a single HBM solution per amplitude level. 
Additional analyses utilizing perturbation methods and experimental tests based on PRNM have also been conducted \cite{abeloosControlbasedMethodsIdentification2022a, zhouIdentificationSecondaryResonances2024}.
As an alternative, extrema tracking has been developed to directly capture the maximum or minimum amplitude point without relying on a phase criteria while using a more computationally efficient approach than HBM
\cite{razeTrackingAmplitudeExtrema2024}.
Though extrema tracking is more accurate than phase resonance approaches, implementations generally require second derivative information about the nonlinearity. 
This results in increased computational cost compared to phase resonance based approaches \cite{razeTrackingAmplitudeExtrema2024}, and it likely limits the current applicability for popular friction models that are piecewise linear (or approximated as such) and cannot easily provide useful second derivatives.
Therefore, a phase resonance approach has more immediate opportunities for modeling jointed structures.
Applications of PRNM have focused on analytical and conservative nonlinearities \cite{volvertPhaseResonanceNonlinear2021, volvertResonantPhaseLags2022, abeloosControlbasedMethodsIdentification2022a}, but jointed structures exhibit hysteretic behavior \cite{mathisReviewDampingModels2020}.
To allow for analysis of hysteretic systems, variable phase resonance nonlinear modes (VPRNM) has been proposed to track superharmonic resonances \cite{porterTrackingSuperharmonic2024}. Hysteretic systems exhibit a non-constant phase of superharmonic resonances and thus cannot be analyzed with a constant phase criteria. 
VPRNM has only been investigated for single DOF (SDOF) systems, motivating further analysis for multiple DOF (MDOF) and real frictional systems. 
For a more thorough discussion of PRNM and VPRNM, see \cite{porterTrackingSuperharmonic2024}.

This paper extends VPRNM \cite{porterTrackingSuperharmonic2024} to MDOF systems, develops a VPRNM ROM for efficiently capturing superharmonic resonances, and compares results to new experiments with a jointed structure.
First, the superharmonic resonance behavior is demonstrated in \Cref{sec:motiv_example} to motivate amplitude control as a basis for understanding superharmonic resonances. 
The equations for VPRNM are introduced in \Cref{sec:theory} and modifications for high dimensional systems are discussed.
The computationally efficient VPRNM ROM is derived and discussed in \Cref{sec:full_rom_theory} and applied to a 3 DOF system in \Cref{sec:3dof_sys}.
To validate the approach, the Half Brake-Reu{\ss} Beam (HBRB) \cite{chenMeasurementIdentificationNonlinear2022} is introduced in \Cref{sec:hbrb_system} and HBM, VPRNM ROM, and new experimental results are compared in \Cref{sec:hbrb_results} utilizing the best available modeling techniques from \cite{porterPredictivePhysicsbasedFriction2023}.
Lastly, conclusions are presented in \Cref{sec:mdof_conclusions}.
The results in \Cref{sec:3dof_sys} and \Cref{sec:hbrb_results} highlight the potential of the VPRNM ROM for capturing superharmonic resonances and demonstrate the computational time savings of the approach.

\section{Motivating Examples} \label{sec:motiv_example}

\begin{figure}[!h]
	\centering
	\begin{subfigure}[c]{0.45\linewidth } 
		\centering
		\includegraphics[width=\linewidth]{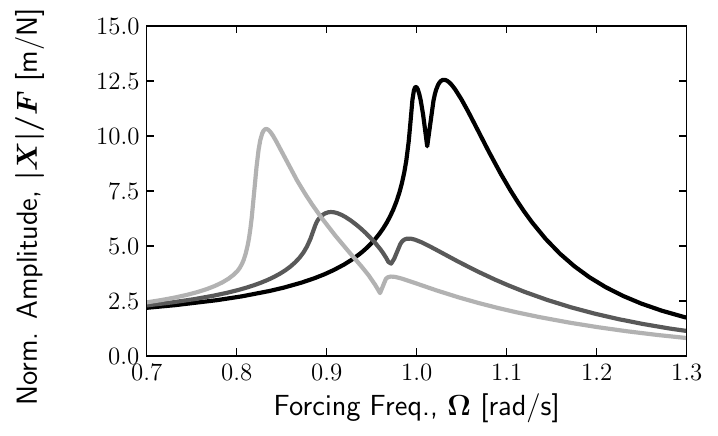}
		\caption{}
		\label{fig:motiv_constantF_tot}
	\end{subfigure}	\quad
	\begin{subfigure}[c]{0.45\linewidth } 
		\centering
		\includegraphics[width=\linewidth]{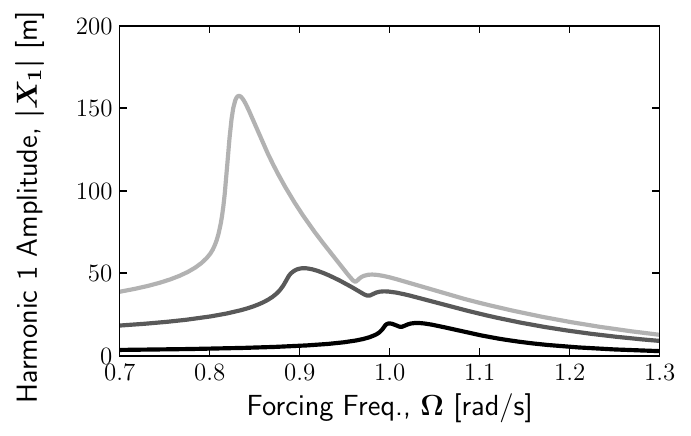}
		\caption{}
		\label{fig:motiv_constantF_h1}
	\end{subfigure}	\\
	\begin{subfigure}[c]{0.45\linewidth } 
		\centering
		\includegraphics[width=\linewidth]{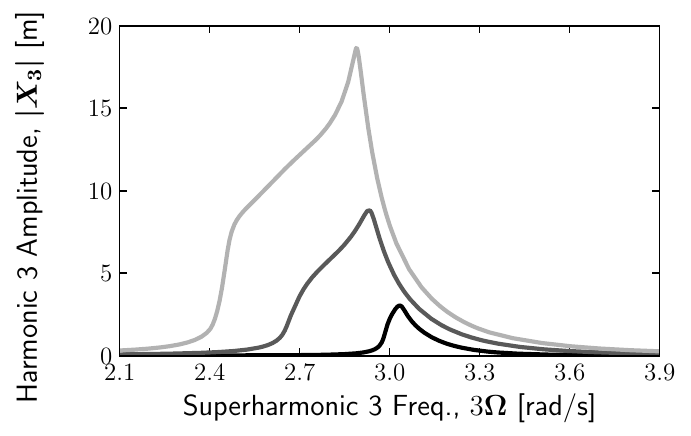}
		\caption{}
		\label{fig:motiv_constantF_h3}
	\end{subfigure}	\quad
	\begin{subfigure}[c]{0.45\linewidth } 
		\centering
		\includegraphics[width=\linewidth]{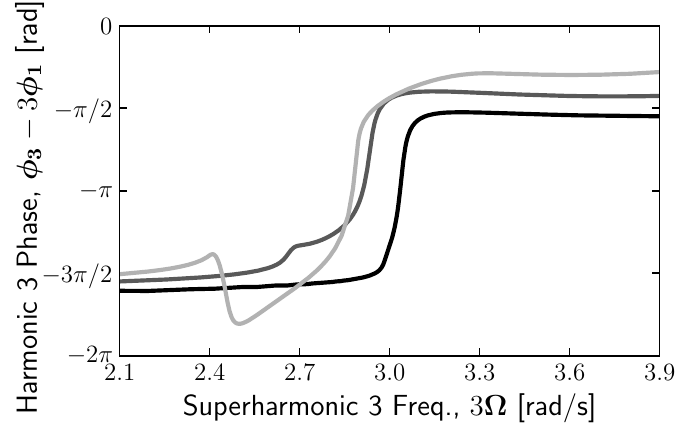}
		\caption{}
		\label{fig:motiv_constantF_phase}
	\end{subfigure}		\caption{Frequency responses for DOF 1 of the 3 DOF system described in \Cref{sec:3dof_sys} for force levels of 1.6 N, 8.0 N, and 16 N from darkest to lightest. The subplots are (a) peak response amplitude for a cycle normalized by force magnitude, (b) magnitude of the first harmonic response, (c) magnitude of the third harmonic response, and (d) phase difference between the third and first harmonics. Note that the forcing frequency ranges are identical for all plots, but plots (c) and (d) plot the superharmonic response frequency, which is three times the forcing frequency. Only plot (a) is normalized by force.} \label{fig:motiv_constantF}
\end{figure}

The present work seeks to understand superharmonic resonances in the neighborhood of primary resonances that result in internal resonances. 
To illustrate this phenomena, FRCs at three force levels for a 3 DOF system with a frictional element are illustrated in \Cref{fig:motiv_constantF}. These responses are calculated for the system described in \Cref{sec:3dof_sys} with HBM as described in \Cref{sec:theory}. This paper seeks to understand the multiple peaks and local minima caused by the superharmonic resonance and develop new computationally efficient approaches for this regime where existing methods break down \cite{krackReducedOrderModeling2013, krackMethodNonlinearModal2013, schwarzValidationTurbineBlade2020}. 
The phase shift in \Cref{fig:motiv_constantF_phase} can be used to isolate the superharmonic resonance using a form of phase resonance that was initially proposed for conservative systems in \cite{volvertPhaseResonanceNonlinear2021} and extended to hysteretic systems with VPRNM in \cite{porterTrackingSuperharmonic2024}.

\begin{figure}[!h]
	\centering
	\begin{subfigure}[c]{0.45\linewidth } 
		\centering
		\includegraphics[width=\linewidth]{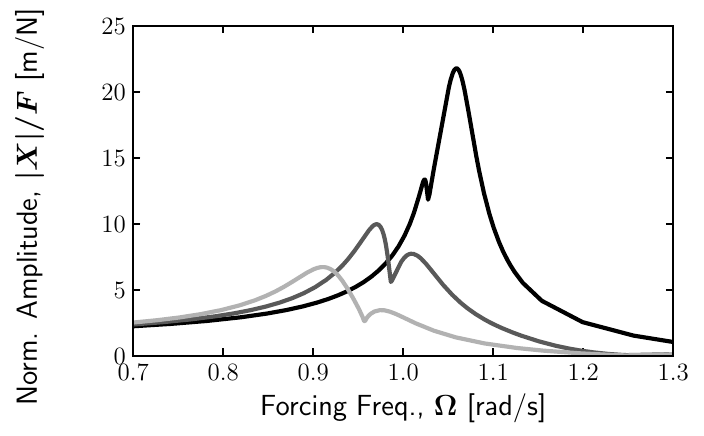}
		\caption{}
		\label{fig:motiv_constantA_tot}
	\end{subfigure}	\quad
	\begin{subfigure}[c]{0.45\linewidth } 
		\centering
		\includegraphics[width=\linewidth]{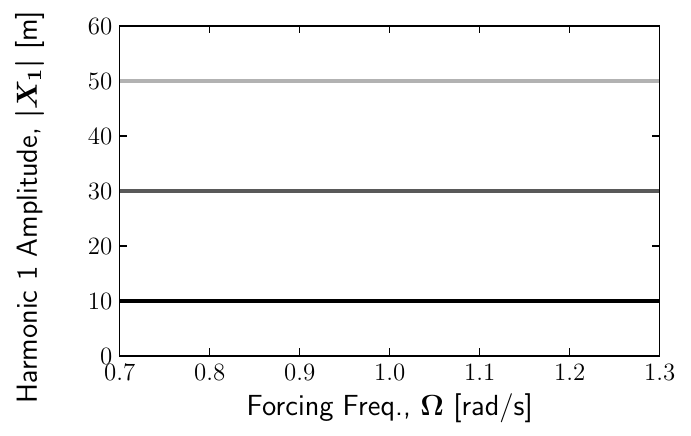}
		\caption{}
		\label{fig:motiv_constantA_h1}
	\end{subfigure}	\\
	\begin{subfigure}[c]{0.45\linewidth } 
		\centering
		\includegraphics[width=\linewidth]{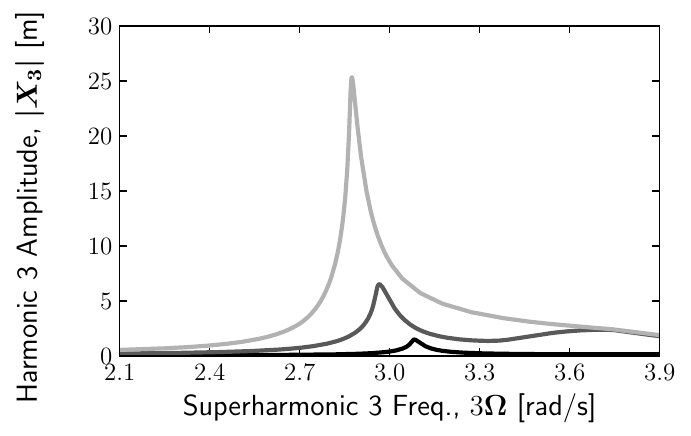}
		\caption{}
		\label{fig:motiv_constantA_h3}
	\end{subfigure}	\quad
	\begin{subfigure}[c]{0.45\linewidth } 
		\centering
		\includegraphics[width=\linewidth]{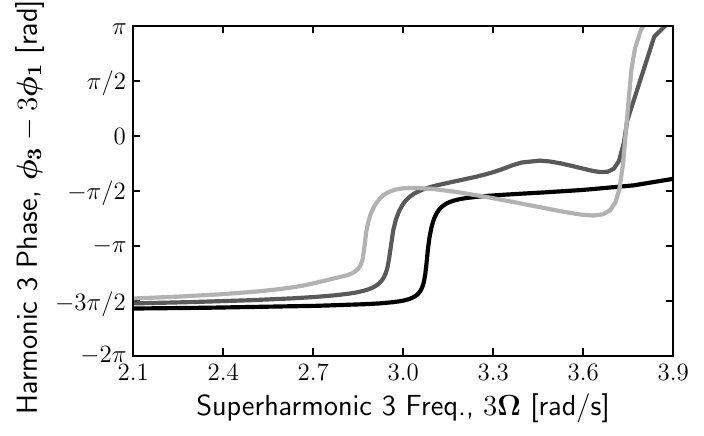}
		\caption{}
		\label{fig:motiv_constantA_phase}
	\end{subfigure}		\caption{
	Frequency responses for DOF 1 of the 3 DOF system described in \Cref{sec:3dof_sys} for controlled first harmonic response amplitude levels of 20 m, 30 m, and 40 m from darkest to lightest. The subplots are (a) peak response amplitude for a cycle normalized by force magnitude, (b) magnitude of the first harmonic response, (c) magnitude of the third harmonic response, and (d) phase difference between the third and first harmonics. Note that the forcing frequency ranges are identical for all plots, but plots (c) and (d) plot the superharmonic response frequency, which is three times the forcing frequency. Only plot (a) is normalized by force.
	} \label{fig:motiv_constantAmp}
\end{figure}

The complicated shape of the superharmonic peak in \Cref{fig:motiv_constantF_h3} is difficult to interpret from the perspective of single nonlinear mode theory. 
Taking inspiration from experimental approaches \cite{karaagacliExperimentalModalAnalysis2021, chenMeasurementIdentificationNonlinear2022}, responses can also be calculated for a constant response amplitude of the first harmonic while allowing the magnitude of the forcing to vary as shown in \Cref{fig:motiv_constantAmp}. 
The superharmonic resonances in \Cref{fig:motiv_constantA_h3} are easier to interpret than those shown in \Cref{fig:motiv_constantF_h3} and look similar to a single nonlinear mode response when taken in isolation. This case motivates the use of amplitude control throughout this paper and is more thoroughly discussed in \Cref{sec:theory}. 
It is noted that an interpolation problem can be solved to go between multiple curves calculated at constant amplitude to curves at constant force as is done in \cite{karaagacliExperimentalModalAnalysis2021}. 
Thus, the change from force to amplitude control is merely for convenience of understanding the same problem. 
As with the constant force case, the phase difference between the first and third harmonics (see \Cref{fig:motiv_constantA_phase}) can be used to identify the superharmonic resonance with VPRNM \cite{porterTrackingSuperharmonic2024}.
This system is more thoroughly investigated and discussed in  \Cref{sec:3dof_sys}.

\FloatBarrier

\section{Theory} \label{sec:theory}

This work focuses on harmonically forced systems with multiple DOFs. For a system with $ N $ displacement DOFs, the equations of motion are
\begin{equation} \label{eq:eom}
	\boldsymbol{M} \boldsymbol{\ddot{x}}
	+
	\boldsymbol{C} \boldsymbol{\dot{x}}
	+
	\boldsymbol{K} \boldsymbol{x}
	+
	\boldsymbol{T} \boldsymbol{f_{nl}}(\boldsymbol{Q} \boldsymbol{x})
	=
	\boldsymbol{F_{ext,0}}
	+
	\boldsymbol{F_{ext}} f_{mag,c} \cos (\Omega t)
	+
	\boldsymbol{F_{ext}} f_{mag,s} \sin (\Omega t)
\end{equation}
for $ N \times 1 $ vectors $ \boldsymbol{x} $, $ \boldsymbol{\dot{x}} $, and $ \boldsymbol{\ddot{x}} $ corresponding to displacements, velocities, and accelerations respectively. 
The linear portions of the system are described by the $ N \times N $ mass $ \boldsymbol{M} $, damping $ \boldsymbol{C} $, and stiffness $ \boldsymbol{K} $ matrices.
The nonlinear forces are calculated for $ N_{nl} $ local displacements defined by the $ N_{nl} \times N $ matrix transformation $ \boldsymbol{Q} $, and the nonlinear forces are distributed to the global degrees of freedom with the $ N \times N_{nl} $ matrix $ \boldsymbol{T} $.
External forcing is applied to the system is the shape of the $ N \times 1 $ vector $ \boldsymbol{F_{ext}} $. This forcing vector is scaled by scalar coefficients $ f_{mag,c} $ and $ f_{mag,s} $ for the sine and cosine components at a frequency of $ \Omega $. For later discussion, the total force magnitude is taken as
\begin{equation}
	f_{mag} = \sqrt{ f_{mag,c}^2 + f_{mag,s}^2 }.
\end{equation}
A $ N \times 1 $ static external force of $ \boldsymbol{F_{ext,0}} $ is applied for the large system in \Cref{sec:hbrb_system}. The solution to the static problem is briefly discussed in \Cref{sec:prestress_theory}.

\subsection{Harmonic Balance Method and Continuation} \label{sec:theory_hbm}

Throughout this paper, HBM solutions are utilized as truth solutions and form a basis for other methods under steady-state harmonic excitation \cite{krackHarmonicBalanceNonlinear2019}. HBM assumes the harmonic form of the solution
\begin{equation} \label{eq:x_t}
	\boldsymbol{x}(t) = 
	\boldsymbol{X_0}
	+
	\sum_{n=1}^H \left[ \boldsymbol{X_{nc}} \cos(n \Omega t) + \boldsymbol{X_{ns}} \sin(n \Omega t) \right] 
	=
	\boldsymbol{X_0} + \sum_{n=1}^H \boldsymbol{X_{n}} \cos(n \Omega t - \boldsymbol{\phi_n})
\end{equation}
where $ H $ harmonics are used, $ \boldsymbol{X_0} $ is the static displacements, $ \boldsymbol{X_{nc}} $ is the harmonic $ n $ cosine displacements, and $ \boldsymbol{X_{ns}} $ is the harmonic $ n $ sine displacements.
For later plotting, it is convenient to define vectors for the total displacement $ \boldsymbol{X_n} $ and phase $ \boldsymbol{\phi_n} $ of each DOF at each harmonic.
HBM gives the discrete equations
\begin{equation} \label{eq:hbm}
	\begin{split}
		\boldsymbol{K} \boldsymbol{X_0} + \boldsymbol{F_{nl, 0}} - \boldsymbol{F_{ext, 0}} &= \boldsymbol{0}
		\\
		(-n \Omega^2 \boldsymbol{M} + \boldsymbol{K}) \boldsymbol{X_{nc}} + (n \Omega \boldsymbol{C}) \boldsymbol{X_{ns}} + \boldsymbol{F_{nl, nc}} - \boldsymbol{F_{ext, nc}} &= \boldsymbol{0} \ \ 
		\forall n \in \{1, \dots, H\}
		\\
		(-n \Omega^2 \boldsymbol{M} + \boldsymbol{K}) \boldsymbol{X_{ns}} - (n \Omega \boldsymbol{C}) \boldsymbol{X_{nc}} + \boldsymbol{F_{nl, ns}} - \boldsymbol{F_{ext, ns}} &= \boldsymbol{0} \ \ 
		\forall n \in \{1, \dots, H\},
	\end{split}
\end{equation}
where $ \boldsymbol{F_{nl, 0}} $ is a $ N\times 1 $ vector of the constant component of the nonlinear forces. For harmonic $ n $ and $ q = c $ or $s $ (corresponding to cosine and sine respectively), the $ N\times 1 $ vectors $ \boldsymbol{F_{ext, nq}} $ and $ \boldsymbol{F_{nl, nq}} $ are the applied external forcing and internal nonlinear forces respectively. 

The nonlinear forces for the HBM solution are calculated with the Alternating Frequency-Time (AFT) method \cite{krackHarmonicBalanceNonlinear2019}.
The AFT method transforms the harmonic displacement coefficients to a discrete time series utilizing \eqref{eq:x_t}. The nonlinear forces $ \boldsymbol{f_{nl}}(\boldsymbol{Q x}(t)) $ are calculated in the time domain, and a Fourier transform is used to return to frequency domain nonlinear forces.
Nonlinear forces in the frequency domain are mapped back to the global coordinates using the matrix $ \boldsymbol{T} $.
To obtain steady-state hysteretic behavior, the nonlinear forces are calculated in the time domain over two cycles and the nonlinear forces from the second cycle are used~\cite{krackHarmonicBalanceNonlinear2019}.

This work utilizes continuation to calculate solutions to systems of equations while allowing a control variable (e.g., frequency for HBM) to vary. Continuation adds an additional constraint equation on a given step and relaxes the control variable to be an unknown.
Each step of continuation makes a prediction tangent to the previous solution and then applies an orthogonal corrector to find a new solution point similar to \cite{peetersNonlinearNormalModes2009, rensonNumericalComputationNonlinear2016}.
The present work utilizes the implementation from \cite{porterTrackingSuperharmonic2024}, which is available at \cite{porterTMDSimPy}.
The output of continuation is several solutions to the system of equations (e.g., HBM) for a range of the control variable (e.g., frequency for HBM).

For HBM solutions at a constant forcing amplitude, the $ (2H + 1) N $ equations of HBM are utilized as formulated in \eqref{eq:hbm} for the $ (2H + 1) N $ unknowns corresponding to the displacement coefficients.
For constant forcing amplitude, $ f_{mag,s} $ is set to zero and $ f_{mag,c} $ is set to a constant value. 
As discussed in \Cref{sec:motiv_example}, it is useful to consider the system under amplitude control. In this work, the magnitude of the first harmonic is controlled to have a constant amplitude of $ A_1 $ for the $ k $th derivative of displacement\footnote{For displacement control, $ k=0 $, and for acceleration control, $ k=2 $.}
by adding the equation
\begin{equation} \label{eq:amp_constraint}
	\Omega^{(2k)} \left[ \boldsymbol{R_1 X_{1c}} \right]^2 + \Omega^{(2k)} \left[ \boldsymbol{R_1 X_{1s}} \right]^2 - A_1^2 = 0
\end{equation}
to HBM where $ \boldsymbol{R_1} $ is a $ 1 \times N $ vector that defines what amplitude should be controlled to be constant. Since an additional equation is added to HBM, $ f_{mag,c} $ is relaxed to be an unknown while $ f_{mag,s} $ is still fixed to be zero. Thus, HBM with amplitude control results in $ (2H + 1) N + 1 $ equations and unknowns.
Utilizing continuation to vary the forcing frequency adds an additional equation and unknown resulting in $ (2H + 1) N + 2$ equations and unknowns.

\subsection{Variable Phase Resonance Nonlinear Modes (VPRNM)} \label{sec:vprnm_theory}

This work investigates VPRNM for MDOF systems to track superharmonic resonances.
Previous work proposed VPRNM for SDOF systems \cite{porterTrackingSuperharmonic2024} as an extension of PRNM for superharmonic resonances \cite{volvertPhaseResonanceNonlinear2021}.
The variable phase criteria of VPRNM is required for tracking superharmonic resonances in frictional systems \cite{porterTrackingSuperharmonic2024}, but previous work did not consider MDOF systems. 

To formulate VPRNM, the displacement of the superharmonic $ n $ of interest is
\begin{equation}
	\boldsymbol{x_n}(t) = \boldsymbol{X_{nc}} \cos(n \Omega t) + \boldsymbol{X_{ns}} \sin(n \Omega t),
\end{equation}
and the displacement considering only the first $ n-1 $ harmonics is 
\begin{equation}
	\boldsymbol{x_{0:(n-1)}}(t) = \boldsymbol{X_0} + \sum_{k=1}^{n-1} \boldsymbol{X_{kc}} \cos(k \Omega t) + \boldsymbol{X_{ks}} \sin(k \Omega t).
\end{equation}
Utilizing AFT, VPRNM calculates the excitation of the $ n $th superharmonic as 
\begin{equation} \label{eq:fbroad}
	\boldsymbol{F_{nq,broad}} = - \boldsymbol{T} \mathcal{F}_{nq} \{ \boldsymbol{f_{nl}} [ \boldsymbol{Q} \boldsymbol{x_{0:(n-1)}}(t) ] \} ,
\end{equation}
where $ \mathcal{F}_{nq} $ represents a Fourier transform for harmonic number $n$. Throughout this work, subscripts $q = c $ or $s$ represent cosine or sine terms respectively for frequency domain quantities. 
In addition, nonlinear forces related to violations of superposition are defined as
\begin{equation} \label{eq:Fsup_def}
	\boldsymbol{F_{kq,sup,n}} = 
	-\boldsymbol{T}\mathcal{F}_{kq} \{ \boldsymbol{f_{nl}}[\boldsymbol{Q}\boldsymbol{x}(t)] - \boldsymbol{f_{nl}}[\boldsymbol{Q}\boldsymbol{x_n}(t)] - \boldsymbol{f_{nl}}[\boldsymbol{Q}\boldsymbol{x_{0:(n-1)}}(t)] \}.
\end{equation}
Therefore, the total nonlinear forces on the $ n $th harmonic can be rewritten as
\begin{equation}
	\boldsymbol{T}\mathcal{F}_{nq}\{\boldsymbol{f_{nl}}\left[\boldsymbol{Qx}(t)\right]\}
	=
	\boldsymbol{T} \mathcal{F}_{nq} \{ \boldsymbol{f_{nl}} [\boldsymbol{Q}\boldsymbol{x_n}(t)] \}
	-
	\boldsymbol{F_{nq,broad}}
	-
	\boldsymbol{F_{nq,sup,n}}.
\end{equation}
Using this rearrangement of nonlinear forces, the HBM equations for the superharmonic can then be rearranged as
\begin{subequations} \label{eq:vprnm_rewrite_hbm}
\begin{equation}
	(-n \Omega^2 \boldsymbol{M} + \boldsymbol{K}) \boldsymbol{X_{nc}} + (n \Omega \boldsymbol{C}) \boldsymbol{X_{ns}} +
	\boldsymbol{T} \mathcal{F}_{nc} \{ \boldsymbol{f_{nl}} [\boldsymbol{Q}\boldsymbol{x_n}(t)] \}
	 = 
	 \boldsymbol{F_{nc,broad}} + \boldsymbol{F_{nc,sup,n}}
\end{equation}
\begin{equation}
	(-n \Omega^2 \boldsymbol{M} + \boldsymbol{K}) \boldsymbol{X_{ns}} - (n \Omega \boldsymbol{C}) \boldsymbol{X_{nc}} 
	+
	\boldsymbol{T} \mathcal{F}_{ns} \{ \boldsymbol{f_{nl}}[\boldsymbol{Q}\boldsymbol{x_n}(t)] \}  = 
	\boldsymbol{F_{ns,broad}} + \boldsymbol{F_{ns,sup,n}},
\end{equation}
\end{subequations}
which is exactly equivalent to the HBM equations for the $ n $th harmonic.
In \eqref{eq:vprnm_rewrite_hbm}, $ \boldsymbol{F_{nq,broad}} $ is interpreted as an external forcing acting on the superharmonic $ n $ resulting in a superharmonic resonance.
In practice, the terms $ \boldsymbol{F_{kq,sup,n}} $ are neglected for defining the VPRNM equation and represent the error in the phase constraint for VPRNM.
Therefore, the VPRNM phase constraint assumes that lower harmonic motion excites the superharmonic resonance, but that influences of the superharmonic resonance on the lower harmonics can be neglected. VPRNM shows errors where this assumption breaks down \cite{porterTrackingSuperharmonic2024}.

To capture phase resonance with respect to $ \boldsymbol{F_{nq,broad}} $ by analogy to primary resonances \cite{peetersDynamicTestingNonlinear2011}, VPRNM adds a constraint equation that the response be orthogonal to $ \boldsymbol{F_{nq,broad}} $ by adding the equation
\begin{equation} \label{eq:vprnm_eqn}
	\dfrac{1}{
		\left|\left |\begin{bmatrix}
			\boldsymbol{F_{nc,broad}}^T & \boldsymbol{F_{ns,broad}}^T
		\end{bmatrix}\right|\right|_2 
	\left|\left |\begin{bmatrix}
	\boldsymbol{X_{nc}}^T & \boldsymbol{X_{ns}}^T
		\end{bmatrix}\right|\right|_2 
	}
	\begin{bmatrix}
		\boldsymbol{F_{nc,broad}}^T & \boldsymbol{F_{ns,broad}}^T
	\end{bmatrix}
	\begin{bmatrix}
		\boldsymbol{X_{nc}} \\ \boldsymbol{X_{ns}}
	\end{bmatrix}
	=
	0
\end{equation}
to the set of equations for HBM \eqref{eq:hbm}.
When the left hand side of \eqref{eq:vprnm_eqn} is positive one, the superharmonic responds completely in phase and in the same direction as the excitation of the superharmonic. When the left side is negative one, then the response is completely out of phase, but in the same shape. 
When \eqref{eq:vprnm_eqn} is satisfied, then the response is considered to be in phase resonance. 
Additional interpretation of \eqref{eq:vprnm_eqn} as modal phase resonance is discussed in \Cref{sec:theory_high_dim}.
Since VPRNM adds an additional constraint equation to HBM, an additional unknown is also necessary. Here, the frequency is allowed to vary at a given amplitude level (e.g., force magnitude) to satisfy the VPRNM constraint.

For the 3 DOF system, the same approach as  \cite{porterTrackingSuperharmonic2024} is utilized where the VPRNM constraint is added to the HBM equations and continuation is conducted with respect to the scaling of the external force $ f_{mag,c} $.
For larger systems, a modified approach is needed to consistently obtain converged solutions, and is presented in the next section.

\subsection{Extension to High Dimensional Systems} \label{sec:theory_high_dim}

In addition to the 3 DOF system previously discussed in \Cref{sec:motiv_example} and fully explored in \Cref{sec:3dof_sys}, this work considers a model of a jointed structure detailed in \Cref{sec:hbrb_system}. The model of the jointed structure has 860 degrees of freedom and 232 frictional nonlinear elements making it significantly more challenging to solve than the 3 DOF system. To make the problem easier to solve, HBM is modified to use an amplitude and phase constraint rather than just constraining the amplitude. Thus, \eqref{eq:amp_constraint} is replaced with the two equations
\begin{subequations} \label{eq:amp_phase_constraint}
\begin{equation} 
	\Omega^k  \boldsymbol{R_1 X_{1c}} - A_1 = 0
\end{equation}
\begin{equation} 
	\Omega^k \boldsymbol{R_1 X_{1s}}  = 0 .
\end{equation}
\end{subequations} 
This requires two additional unknowns compared to HBM as both $ f_{mag,c} $ and $ f_{mag,s} $ are allowed to vary. 
With this modification, HBM with amplitude and phase control requires solutions to $ (2H + 1) N + 2 $ equations and unknowns. Utilizing continuation adds an additional equation and unknown resulting in $ (2H + 1) N + 3  $ equations and unknowns. 
Given this modification, the first harmonic response remains close to constant near a primary resonance rather than undergoing a phase shift and only the two unknowns $ f_{mag,c} $ and $ f_{mag,s} $ experience significant changes near the primary resonance allowing for large continuation steps. The superharmonic resonance still requires small continuation steps to resolve, but is expected to be easier to calculate when the first harmonic solution remains at constant amplitude and phase resulting in a more constant excitation $ \boldsymbol{F_{nq,broad}} $ of the superharmonic resonance.

To make VPRNM easier to solve, the HBM equations of \eqref{eq:hbm} are augmented with \eqref{eq:amp_phase_constraint} in addition to a modified VPRNM constraint (discussed next). Thus, the system contains $ (2H + 1) N + 3 $ equations and an equal number of unknowns comprised of the $ (2H + 1) N $ harmonic displacements, the force scaling for both the cosine and sine terms ($ f_{mag,c} $ and $ f_{mag,s} $ respectively), and the frequency at a given amplitude level.
For the large problem, continuation is conducted with respect to the response amplitude rather than the force level for VPRNM. This is also convenient since the experiments are conducted with acceleration control, so VPRNM can easily be initialized at and run between experimental acceleration amplitudes.
Further solver considerations for the large models are detailed in \Cref{sec:solver_consider}.

Initial simulations with VPRNM for the jointed structure showed that VPRNM, as formulated by \eqref{eq:vprnm_eqn}, sometimes resulted in two nearby solutions or no solutions (see \Cref{sec:results_modal_filter} for details). This challenge is likely related to conditions for isolating a single mode with velocity feedback \cite{scheelNonlinearModalTesting2022}. To solve this issue, a modal filter is applied to the VPRNM equation. 
The superharmonic resonance corresponds to the nonlinear excitation of a higher mode, and the linear mode shape of that mode $ \boldsymbol{\psi_{lin,n}} $ is utilized for the modal filter. 
Projecting the motion and force onto this mode shape yields
\begin{subequations}
\begin{equation}
	\boldsymbol{q_{vprnm}}^T = 
	\boldsymbol{\psi_{lin,n}}^T \boldsymbol{M}
	\begin{bmatrix}
		\boldsymbol{X_{nc}} & \boldsymbol{X_{ns}}
	\end{bmatrix}
\end{equation}
\begin{equation} \label{eq:modal_force_vprnm}
	\boldsymbol{\eta_{vprnm}}^T = 
	\boldsymbol{\psi_{lin,n}}^T 
	\begin{bmatrix}
		\boldsymbol{F_{nc,broad}} & \boldsymbol{F_{ns,broad}}
	\end{bmatrix}
\end{equation}
\end{subequations}
These $ 2 \times 1 $ modal quantities (corresponding to the cosine and sine terms) are then utilized in VPRNM as 
\begin{equation} \label{eq:vprnm_eqn_modal_filter}
	\dfrac{\boldsymbol{\eta_{vprnm}}^T
		\boldsymbol{q_{vprnm}}}{
		\left|\left|\boldsymbol{\eta_{vprnm}}\right|\right|_2 
		\left|\left|\boldsymbol{q_{vprnm}}\right|\right|_2 
	}
	=
	0.
\end{equation}
This representation in the modal space corresponds to modal phase resonance with respect to the modal forcing as expressed by orthogonality in the complex plane and is consistent with the SDOF formulation of VPRNM \cite{porterTrackingSuperharmonic2024}.
If the motion of the superharmonic resonance is purely in the shape of $ \boldsymbol{\psi_{lin,n}} $, the numerator of \eqref{eq:vprnm_eqn_modal_filter} is exactly consistent with that of \eqref{eq:vprnm_eqn}.
If the mode shape of the superharmonic mode changes significantly with amplitude, \eqref{eq:vprnm_eqn_modal_filter} could be extended to include projections onto multiple constant basis vectors (e.g., by including modal derivatives \cite{vizzaccaroDirectComputationNonlinear2021}) while still filtering out much of the motion from other modes.

\section{VPRNM Based Reduced Order Model} \label{sec:full_rom_theory}

ROMs utilizing nonlinear modal solutions can significantly reduce computational costs. 
For a single isolated mode, EPMC can be used to characterize the response \cite{krackNonlinearModalAnalysis2015}. EPMC is detailed further in \Cref{sec:epmc_theory} and the implementation is available at \cite{porterTMDSimPy}.
The EPMC solution can then be utilized to reconstruct FRCs for isolated resonances as described in \cite{krackReducedOrderModeling2013, krackMethodNonlinearModal2013, schwarzValidationTurbineBlade2020}.
The details for two forms of a single mode ROM are briefly detailed here based on EPMC (Sections \ref{sec:phase_amp_rom} and \ref{sec:epmc_control_amp_rom}). Then the novel VPRNM ROM utilizing a combination of EPMC for individual modes and VPRNM to characterize the interaction between modes with a superharmonic resonance is presented in \Cref{sec:vprnm_rom}. 
The steps for implementing the VPRNM ROM with minimal discussion are summarized in \Cref{sec:rom_summary_steps} as an additional reference and high level overview.

\subsection{Constant Force EPMC ROM} \label{sec:phase_amp_rom}

EPMC (summarized in \Cref{sec:epmc_theory}) calculates amplitude dependent modal frequency $ \omega_{i} $, damping factor $ \zeta_{i} $, modal displacements $ \boldsymbol{X_{i,EPMC,nc}} $ and $ \boldsymbol{X_{i,EPMC,ns}} $ (harmonic $ n $ cosine and sine displacements respectively) for mode $ i $ and modal amplitude $ q_{i} $. The first harmonic contributions to the EPMC mode shape are represented in mass normalized and complex form as
\begin{equation} \label{eq:complex_mode_shape}
	\boldsymbol{\psi_{i}} = 
	\dfrac{1}{q_{i}} 
	\left[ \boldsymbol{X_{i,EPMC,1c}} - \mathrm{j} \boldsymbol{X_{i,EPMC,1s}} \right]
\end{equation}
where $ \mathrm{j} = \sqrt{-1} $. This is consistent with \eqref{eq:x_t} in that
\begin{equation}
	\Re \{ q_i \boldsymbol{\psi_{i}} e^{\mathrm{j} \Omega t}\}  
	=
	\boldsymbol{X_{i,EPMC,1c}} \cos(\Omega t) 
	+ \boldsymbol{X_{i,EPMC,1s}} \sin(\Omega t) 
	.
\end{equation}
The equations of motion \eqref{eq:eom} can then be projected onto a single degree of freedom $ q_{i} $ similar to \cite{krackReducedOrderModeling2013}, yielding
\begin{equation} \label{eq:mode_project}
	\left[-\Omega^2 + \mathrm{j} 2 \Omega \zeta_{i} \omega_{i} + \omega_{i}^2\right] q_{i} e^{\mathrm{j} (\Omega t - \varphi)} = 
	\boldsymbol{\psi_{i}}^H \boldsymbol{F_{ext,1}} e^{\mathrm{j} \Omega t},
\end{equation}
where $ (\cdot)^H $ denotes the Hermitian (conjugate) transpose, $ \varphi $ is the phase of the mode $ i $, and
\begin{equation}
	\boldsymbol{F_{ext,1}}
	=
	f_{mag,c} \boldsymbol{F_{ext}} - \mathrm{j} f_{mag,s} \boldsymbol{F_{ext}}.
\end{equation}

First the relationship between forcing frequency and response amplitude is determined. The forcing frequency to achieve a given modal response amplitude $ q_i $ of a mode under constant force excitation of the system can be calculated as \cite{schwarzValidationTurbineBlade2020}
\begin{equation} \label{eq:epmc_amp_only_rom}
	\Omega^2 = p_2 
	\pm 
	\sqrt{
	p_2^2 - \omega_i^4 + \dfrac{\left| \boldsymbol{\psi_i}^H \boldsymbol{F_{ext,1}} \right|^2}{q_i^2}
	}
\end{equation}
for
\begin{equation}
	p_2 = \omega_{i}^2 - 2 (\omega_{i} \zeta_{i})^2 .
\end{equation}
From \eqref{eq:epmc_amp_only_rom}, only real positive values of $ \Omega $ are retained with complex values indicating that a response amplitude is not achieved for a given forcing level.

For many cases with a single nonlinear mode, \eqref{eq:epmc_amp_only_rom} is sufficient for understanding the system. However for the later development of the VPRNM ROM, it is necessary to calculate the phase of the responses.
Given an assumed response amplitude and a forcing frequency from \eqref{eq:epmc_amp_only_rom}, the phase $ \varphi $ of the modal response in \eqref{eq:mode_project} can be calculated by multiplying both sides of \eqref{eq:mode_project} by the conjugate of $ \left[\boldsymbol{\psi_{i}}^H \boldsymbol{F_{ext,1}} \right]$ and simplifying,
yielding
\begin{equation}
	\left[\boldsymbol{\psi_{i}}^H \boldsymbol{F_{ext,1}} \right]^H 
	\left[-\Omega^2 + \mathrm{j} 2 \Omega \zeta_{i} \omega_{i} + \omega_{i}^2\right] q_{i} e^{-\mathrm{j} \varphi} = 
	\left| \boldsymbol{\psi_{i}}^H \boldsymbol{F_{ext,1}} \right|^2 .
\end{equation}
Here, the right hand side is purely real, thus a four-quadrant arctangent operator can be applied to the product of the first two bracketed terms on the left to determine the phase $ \varphi $ as
\begin{equation} \label{eq:epmc_phase_rom}
	\varphi = \text{angle}\left[ \left[\boldsymbol{\psi_{i}}^H \boldsymbol{F_{ext,1}} \right]^H 
	\left[-\Omega^2 + \mathrm{j} 2 \Omega \zeta_{i} \omega_{i} + \omega_{i}^2\right] \right]
\end{equation}

Given the phase information, it is useful to format the EPMC ROM with phase information to match the format of HBM solutions. To that end, response vectors for each harmonic $ n $ and cosine $ c $ and sine $ s $ components are phase rotated as
\begin{equation} \label{eq:epmc_phase_rotate}
	\begin{bmatrix}
		\boldsymbol{X_{i,EPMC,ROM,nc}}^T
		\\
		\boldsymbol{X_{i,EPMC,ROM,ns}}^T
	\end{bmatrix}
	=
	\begin{bmatrix}
		\cos (n \varphi) & - \sin (n \varphi)
		\\
		\sin (n \varphi) & \cos (n \varphi)
	\end{bmatrix}
	\begin{bmatrix}
		\boldsymbol{X_{i,EPMC,nc}}^T
		\\
		\boldsymbol{X_{i,EPMC,ns}}^T
	\end{bmatrix}
\end{equation}
Here, the left side of the equation is the response vector for the ROM at a given phase $ \varphi $ corresponding to a specific forcing frequency $ \Omega $. 
The static contributions $ X_0 $ for the EPMC ROM can be taken directly from the EPMC solution without modification.

\subsection{Constant Amplitude EPMC ROM}
\label{sec:epmc_control_amp_rom}

As motivated in \Cref{sec:motiv_example}, it is beneficial to consider the superharmonic and internal resonance under amplitude control of the first harmonic motion. 
Therefore an EPMC based ROM for amplitude control of an isolated mode is formulated by rearranging \eqref{eq:epmc_amp_only_rom} as
\begin{equation} \label{eq:epmc_amp_control}
		q_i \sqrt{\Omega^4 - 2 \Omega^2 p_2  + \omega_{i}^4 } =  f_{mag} \left| \boldsymbol{\psi_{i}}^H  \left(  \boldsymbol{F_{ext}} e^{\mathrm{j} \varphi_f} \right)  \right| ,
\end{equation}
where $ \boldsymbol{F_{ext}} $ is the previously defined forcing direction, $ \varphi_f $ is the phase of the external forcing, $ f_{mag} $ is the total magnitude of the external forcing, and $ \Omega $ is a chosen forcing frequency of interest. Here, phase $ \varphi_f $ can be chosen arbitrarily because only the magnitude of the modal force is taken.
To calculate the response at a given displacement amplitude, the correct EPMC modal amplitude $ q_i $ is selected to give the desired displacement amplitude (generally requiring linear interpolation between calculated points). Then \eqref{eq:epmc_amp_control} can be solved for $ f_{mag} $ since everything else is already known utilizing the nonlinear modal properties at amplitude $ q_i $.
This EPMC ROM only considers a single nonlinear mode; consequently, it is only expected to be accurate near an isolated resonance (e.g., where a single mode dominates the first harmonic response). 
Note that for this EPMC ROM, the forcing frequency $ \Omega $ can be chosen as any value of interest without requiring the solution to a set of nonlinear equations unlike the EPMC ROM presented in \Cref{sec:phase_amp_rom}, which calculates forcing frequencies at the input modal amplitudes.

\subsection{Superharmonic Resonance VPRNM ROM}\label{sec:vprnm_rom}

The superharmonic resonance VPRNM ROM is based on a controlled displacement amplitude response of the first harmonic as motivated by \Cref{sec:motiv_example}. 
Since the present work considers hysteretic nonlinearities that are purely a function of the displacement history, a constant displacement amplitude approximates a constant nonlinear force $ \boldsymbol{F_{nq,broad}} $ from VPRNM (see \Cref{sec:vprnm_theory}). 
This approximation is most accurate near an isolated resonance where a single mode dominates the vibration shape of the fundamental motion over the forcing frequency range of interest.

While VPRNM calculates $ \boldsymbol{F_{nq,broad}} $ for the phase constraint, the force $ \boldsymbol{F_{nq,broad}} $ is not necessarily a good representation of the magnitude of the excitation of the superharmonic resonance due to $ \boldsymbol{F_{nq,sup,n}} $ in \eqref{eq:vprnm_rewrite_hbm}. 
Therefore, the VPRNM ROM approximates the magnitude of the excitation of the superharmonic resonance utilizing single nonlinear mode theory. 
First, the VPRNM continuation solution is interpolated to a single point
\begin{equation} \label{eq:interp_vprnm}
	(\boldsymbol{R_1} \boldsymbol{X_{VPRNM,1c}})^2 + (\boldsymbol{R_1} \boldsymbol{X_{VPRNM,1s}})^2 = A_{VPRNM\ ROM,1}^2,
\end{equation}
where $ \boldsymbol{R_1} $ is an $ 1 \times N $ vector that defines the DOF that is controlled to have a constant amplitude of $ A_{VPRNM\ ROM,1} $.
At the interpolated point, the harmonic displacements are $ \boldsymbol{X_{VPRNM,nq}} $, the forcing frequency is $ \Omega_{VPRNM} $, and the external forcing magnitude is $ f_{mag,VPRNM} $.

Next, the EPMC solution for the nonlinear mode corresponding to the $ n $th superharmonic response is interpolated to match the superharmonic response amplitude of VPRNM, extracted with the $ 1 \times N $ vector $ \boldsymbol{R_n} $, by satisfying the equation
\begin{equation} \label{eq:interp_super_EPMC}
	(\boldsymbol{R_n} \boldsymbol{X_{VPRNM,nc}})^2 + (\boldsymbol{R_n} \boldsymbol{X_{VPRNM,ns}})^2
	=
	(\boldsymbol{R_n} \boldsymbol{X_{S,EPMC,1c}})^2 + (\boldsymbol{R_n} \boldsymbol{X_{S,EPMC,1s}})^2
	.
\end{equation}
In this case, $ \boldsymbol{X_{S,EPMC,1q}} $ represents the fundamental motion of the EPMC solution for the nonlinear superharmonic mode.
At the interpolated EPMC solution, the modal frequency, damping, and amplitude are $ \omega_{S,EPMC} $, $ \zeta_{S,EPMC} $, and $ q_{S,EPMC} $ respectively.
Note that the fundamental motion of the superharmonic EPMC mode is at approximately the same frequency as the $ n $th harmonic of the VPRNM solution because the superharmonic resonance occurs at $ \omega_{S,EPMC}\approx n \Omega_{VPRNM} $.

Utilizing these interpolated points, the modal forcing on the superharmonic resonance is approximated based on \eqref{eq:epmc_amp_only_rom} with forcing frequency and superharmonic natural frequency both equal to $ n \Omega_{VPRNM} $ and modal properties from the interpolated superharmonic EPMC solution resulting in
\begin{equation} \label{eq:superharmonic_force}
	\left| \boldsymbol{\psi_{S}}^H \boldsymbol{F_{S}} \right|
	=
	2 q_{S,EPMC} \left( n \Omega_{VPRNM} \right)^2  \zeta_{S,EPMC}.
\end{equation}
It is assumed that the mode shape $ \boldsymbol{\psi_S} $ is constant for calculating the modal forcing of the superharmonic so that the quantity $ \left| \boldsymbol{\psi_{S}}^H \boldsymbol{F_{S}} \right| $ can be directly calculated a single time. 
For linear systems, phase resonance occurs at the natural frequency, so this calculation is consistent. 
Furthermore for a single linear mode, the derivation of \eqref{eq:superharmonic_force} and \Cref{sec:phase_amp_rom} exactly give the frequency response taking input of the response amplitude of the mode at phase resonance.

Utilizing $ \left| \boldsymbol{\psi_{S}}^H \boldsymbol{F_{S}} \right| $, the response of the superharmonic can be calculated for a range of response frequencies utilizing \eqref{eq:epmc_amp_only_rom}, \eqref{eq:epmc_phase_rom}, and \eqref{eq:epmc_phase_rotate}. For this calculation, the original EPMC points are used plus the interpolated point that satisfies \eqref{eq:interp_super_EPMC}.
The modal forcing is assumed to be purely real at this stage. 
The nonlinear modal frequency of the superharmonic mode is scaled as 
\begin{equation} \label{eq:scale_superharmonic_mode_freq}
	\omega_{modified}(q_i) = \omega_i(q_i) \dfrac{n \Omega_{VPRNM}}{\omega_{S,EPMC}}
\end{equation}
so that the superharmonic resonance peak matches the forcing frequency found with VPRNM.
Note that $ \zeta_i(q_i) $ from EPMC is not modified by this frequency shift and is calculated from $ \xi_i $ based on the unmodified frequency (see \Cref{sec:epmc_theory}).
The result of this EPMC ROM is a set of response vectors $ \boldsymbol{X_{S,EPMC,ROM,kq}}(\Omega_S) $ for the harmonic $ k $ and $ q $ denoting cosine or sine at superharmonic frequencies $ \Omega_S = n \Omega $ where $ \Omega_S $ is the result from \eqref{eq:epmc_amp_only_rom}.
Additionally, the modal amplitude for each response point $ q_S(\Omega_S) $ is saved for later.

Next, the phase of the superharmonic resonance needs to be matched to the VPRNM solution since the modal forcing from \eqref{eq:superharmonic_force} is assumed to be purely real.
The phase difference between the fundamental and superharmonic motions is critical to capture with the VPRNM ROM because different relative phases could result in adding or subtracting the fundamental and superharmonic response amplitudes to determine the maximum amplitude over a cycle.
Taking the ROM response at $ \Omega_S = n \Omega_{VPRNM} $ as 
$ \boldsymbol{X_{S,EPMC,ROM,1q}}(n \Omega_{VPRNM}) $, the phase difference between the ROM and the VPRNM solution can be calculated as
\begin{equation}
\begin{split} \label{eq:super_phase_diff}
	\varphi = &
	arctan2(\boldsymbol{R_n} \boldsymbol{X_{VPRNM,ns}}, \boldsymbol{R_n}  \boldsymbol{X_{VPRNM,nc}})
	\\
	& -
	arctan2(\boldsymbol{R_n} \boldsymbol{X_{S,EPMC,ROM,1s}}(n\Omega_{VPRNM}), \boldsymbol{R_n}  \boldsymbol{X_{S,EPMC,ROM,1c}}(n\Omega_{VPRNM}))
\end{split}
\end{equation}
where $ arctan2 $ is a four-quadrant arctangent operator. 
The response vector at all forcing frequencies and harmonics $ k $ is rotated as 
\begin{equation} \label{eq:super_phase_rot}
	\begin{bmatrix}
		\boldsymbol{X_{S,ROM, kc}}^T \\
		\boldsymbol{X_{S,ROM, ks}}^T
	\end{bmatrix} 
	=
	\begin{bmatrix}
		\cos (k\varphi) & -\sin(k\varphi) \\
		\sin(k\varphi) & \cos (k\varphi)
	\end{bmatrix}
	\begin{bmatrix}
		\boldsymbol{X_{S,EPMC,ROM,kc}}^T \\
		\boldsymbol{X_{S,EPMC,ROM,ks}}^T
	\end{bmatrix}
\end{equation}
giving the VPRNM ROM solution for the superharmonic components.

Next, the mode at the fundamental forcing frequency is considered. The amplitude of the fundamental EPMC mode denoted with subscript $ F $ is interpolated to the controlled amplitude as
\begin{equation}\label{eq:interp_epmc_fund}
	(\boldsymbol{R_1} \boldsymbol{X_{F,EPMC,1c}})^2 + (\boldsymbol{R_1} \boldsymbol{X_{F,EPMC,1s}})^2 = A_{VPRNM\ ROM,1}^2
	,
\end{equation}
with the modal frequency, damping, and amplitude of $ \omega_{F,EPMC} $, $ \zeta_{F,EPMC} $, and $ q_{F,EPMC} $ respectively.
Then, the fundamental EPMC solution (only at the single interpolated point) is phase rotated to match the VPRNM phase with
\begin{subequations}
\begin{equation}\label{eq:fund_phase_diff}
		\varphi = 
		arctan2(\boldsymbol{R_1} \boldsymbol{X_{VPRNM,ns}}, \boldsymbol{R_1}  \boldsymbol{X_{VPRNM,nc}})
		-
		arctan2(\boldsymbol{R_1} \boldsymbol{X_{F,EPMC,1s}}, \boldsymbol{R_1}  \boldsymbol{X_{F,EPMC,1c}})
\end{equation}
\begin{equation} \label{eq:fund_phase_rot}
	\begin{bmatrix}
		\boldsymbol{X_{F,ROM, kc}}^T \\
		\boldsymbol{X_{F,ROM, ks}}^T
	\end{bmatrix} 
	=
	\begin{bmatrix}
		\cos (k\varphi) & -\sin(k\varphi) \\
		\sin(k\varphi) & \cos (k\varphi)
	\end{bmatrix}
	\begin{bmatrix}
		\boldsymbol{X_{F,EPMC,kc}}^T \\
		\boldsymbol{X_{F,EPMC,ks}}^T
	\end{bmatrix}
\end{equation}
\end{subequations}
for all harmonics $ k $. 
For this work, it is important to note that the EPMC backbone corresponding to the mode at the fundamental forcing frequency is calculated without including the harmonic corresponding to the superharmonic resonance. This prevents the EPMC backbone from displaying tongues commonly observed in NNMs near internal resonances \cite{kerschenNonlinearNormalModes2009, krackMethodNonlinearModal2013} since the present work uses EPMC to characterize the two modes separately and VPRNM to capture the internal resonance behavior.

At this point a full response vector can be constructed at each forcing frequency $ \Omega $ consisting of the static solution $ \boldsymbol{X_{VPRNM,0}} $ from the interpolated VPRNM solution, the rotated fundamental EPMC ROM from \eqref{eq:fund_phase_rot} (repeated for all frequencies), and the rotated superharmonic resonance ROM from \eqref{eq:super_phase_rot} (at the frequencies $ \Omega_S = n \Omega $). Here the harmonic indices $ k $ for the superharmonic response in \eqref{eq:super_phase_rot} need to be multiplied by $ n $ when the solutions are combined (e.g., harmonic 1 of the superharmonic response becomes harmonic $ n $, harmonic 2 becomes $ 2n $, etc.).

With the interpolated fundamental EPMC properties, the force to achieve constant amplitude can be calculated with \eqref{eq:epmc_amp_control} at frequencies\footnote{The superharmonic response was previously calculated at superharmonic frequencies $\Omega_S$.} $ \Omega = \Omega_S / n $ yielding $ f_{mag,fund,rom}(\Omega) $.
If higher resolution of the force is desired, additional forcing frequencies $ \Omega $ can be considered in the force calculation while linearly interpolating the harmonic displacements to the new frequencies.
This force can be distorted by $ \boldsymbol{F_{1q,sup,n}} $ and thus a correction is optionally applied. The difference in the force at the VPRNM solution is calculated as 
\begin{equation} \label{eq:rom_deltaF}
	\Delta F = f_{mag, VPRNM} - f_{mag,fund,rom}(\Omega_{VPRNM}).
\end{equation}
The necessary correction is assumed to be proportional to the superharmonic response modal amplitude $ q_S(n\Omega) $ based on \eqref{eq:epmc_amp_only_rom} yielding the corrected force
\begin{equation} \label{eq:vprnm_rom_force}
	f_{VPRNM,ROM}(\Omega) = 
	f_{mag,fund,rom}(\Omega) 
	+ \Delta F \dfrac{q_S(n\Omega)}{\max_{\Omega_S} q_S(\Omega_S)}.
\end{equation}
The steps of the VPRNM ROM are summarized in \Cref{sec:rom_summary_steps} to help in understanding the implementation.

\subsection{Discussion on VPRNM ROM Formulation}

Many steps require matching EPMC and VPRNM amplitude and phase, but exactly matching all DOFs is a highly overdetermined problem. 
Thus a scalar function is utilized to calculate the amplitude or phase. 
Mass normalized modal amplitudes and phases based on the the complex product $ \psi_i^H \boldsymbol{M} \psi_j$ could be utilized, but this would require the full mass matrix for the VPRNM ROM instead of just the backbone information.
Therefore, amplitude and phase information in the VPRNM ROM are extracted at a linear combination of DOFs utilizing $ 1\times N $ vectors $ \boldsymbol{R_1} $ and $ \boldsymbol{R_n} $ for the fundamental and superharmonic motions respectively. 
Vector products define the scalar function for the interpolation in \eqref{eq:interp_vprnm}, \eqref{eq:interp_super_EPMC}, \eqref{eq:interp_epmc_fund} and for matching phases in \eqref{eq:super_phase_diff} and \eqref{eq:fund_phase_diff}.
These vectors may generally be different, and different vectors are used in \Cref{sec:hbrb_results} since the control point for the first harmonic has much smaller displacement for the superharmonic resonance than other points.

The present work only considers VPRNM to track the superharmonic resonance. An alternative approach would be to use extrema tracking \cite{razeTrackingAmplitudeExtrema2024}, but that is not applicable in this case because of the smoothness requirements for the nonlinearities and computational costs.
It is expected that the present VPRNM ROM should be mostly applicable to extrema tracking of a superharmonic resonance.
If used with extrema tracking,  \eqref{eq:superharmonic_force} should be reformulated based on the amplitude resonance frequency of
\begin{equation}
	n \Omega_{VPRNM} = \omega_{EPMC,S} \sqrt{1 - \zeta_{S,EPMC}^2}
\end{equation}
to remain consistent with linear systems
and the frequency scaling in \eqref{eq:scale_superharmonic_mode_freq} should be consistently updated.

\section{Three Degree of Freedom System} \label{sec:3dof_sys}

\subsection{System Description} \label{sec:3dof_description}

\begin{figure}[h!]
	\centering
	\includegraphics[width=0.55\linewidth]{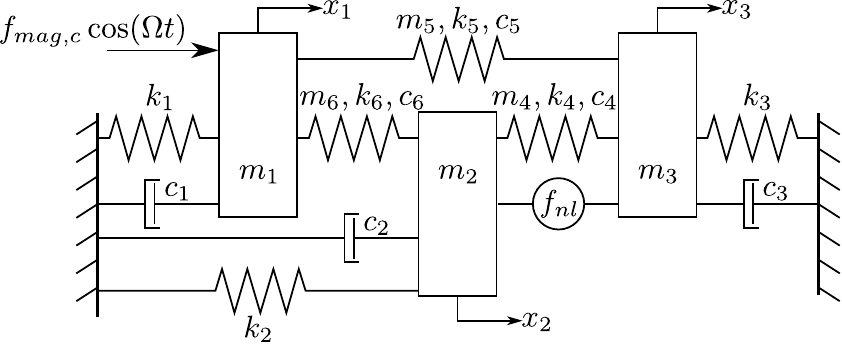}
	\caption{Schematic of 3 DOF system with nonlinear element.}
	\label{fig:3dof_schematic}
\end{figure}

\begin{figure}[!h]
	\centering
	\begin{subfigure}[c]{0.45\linewidth } 
		\centering
		\includegraphics[width=\linewidth]{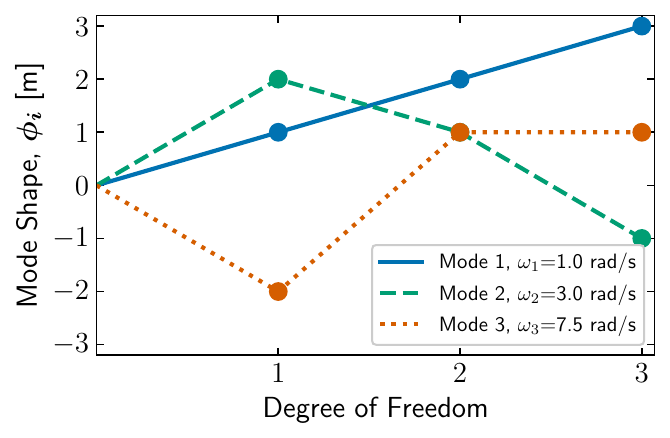}
		\caption{}
		\label{fig:3dof_overview_shapes}
	\end{subfigure}	\quad
	\begin{subfigure}[c]{0.45\linewidth } 
		\centering
		\includegraphics[width=\linewidth]{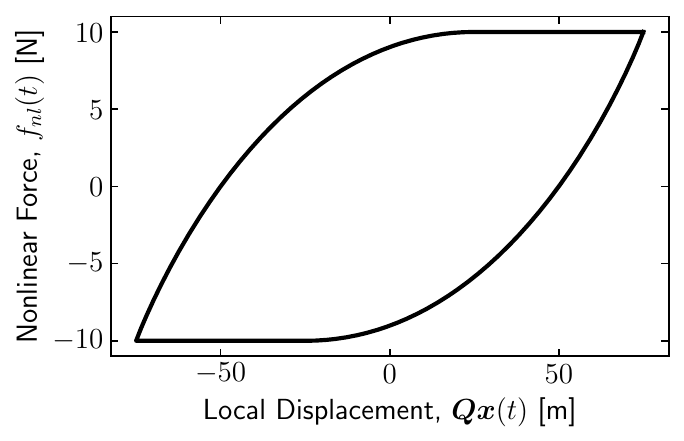}
		\caption{}
		\label{fig:3dof_iwan_hyst}
	\end{subfigure}	\\
	\begin{subfigure}[c]{0.45\linewidth } 
		\centering
		\includegraphics[width=\linewidth]{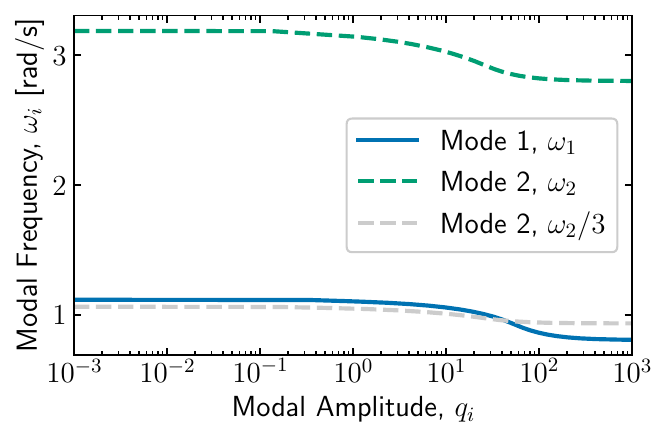}
		\caption{}
		\label{fig:3dof_epmc_freq}
	\end{subfigure}	\quad
	\begin{subfigure}[c]{0.45\linewidth } 
		\centering
		\includegraphics[width=\linewidth]{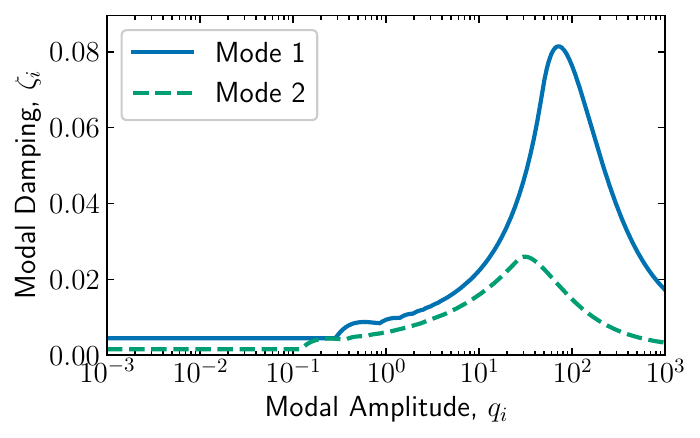}
		\caption{}
		\label{fig:3dof_epmc_damp}
	\end{subfigure}		\caption{Overview of the 3 DOF system: (a) mass-normalized mode shapes around a linear state including half of the linearized frictional force (representative of partial slip) and (b) example hysteresis for the four-parameter Iwan model used in this study. EPMC is used to calculate backbones of the system for (c) frequency and (d) damping of the first two modes.} \label{fig:3dof_overview}
\end{figure}

A 3 DOF system (see \Cref{fig:3dof_schematic}) is constructed to show a superharmonic and internal resonances between the first and second modes of the structure with an Iwan element nonlinearity.
Mode shapes shown in \Cref{fig:3dof_overview_shapes} are prescribed when the Iwan element is partially slipping resulting in dense $ 3 \times 3 $ mass and stiffness matrices (see full details in \Cref{sec:3dof_deriv}). 
Additionally, the modal frequencies are chosen to be $ \omega = $ 1.0 rad/s, 3.0 rad/s, and 7.5 rad/s by construction when the system is partially slipping.
This gives an expected 3:1 superharmonic resonance of the second mode near the primary resonance of the first mode resulting in an internal resonance. 
Mass proportional damping of $ \boldsymbol{C} = 0.01 \boldsymbol{M} $ is utilized.
The system is externally forced only at the first degree of freedom defining the vector 
\begin{equation}
	\boldsymbol{F_{ext}} = \begin{bmatrix}
		1 \\ 0 \\ 0
	\end{bmatrix},
\end{equation}
and only cosine forcing is considered (i.e., $ f_{mag} = f_{mag,c} \neq 0 $ and $ f_{mag,s} = 0 $ in \eqref{eq:eom}).
A four-parameter Iwan element \cite{segalmanFourParameterIwanModel2005} is utilized to contribute nonlinearity and acts between the second and third DOFs with transform matrices
\begin{equation}
	\boldsymbol{Q} = \boldsymbol{T}^T = \begin{bmatrix}
		0 & 1 & -1
	\end{bmatrix}
\end{equation}
corresponding to \eqref{eq:eom}.
An example hysteresis loop for the four-parameter Iwan model is shown in \Cref{fig:3dof_iwan_hyst}. The exact formulation for the four-parameter Iwan model is identical to \cite{porterTrackingSuperharmonic2024} (with $k_t=0.6$ N/s, $F_s=10.0$ N, $\chi=-0.5$, $\beta=0$); the implementation is available at \cite{porterTMDSimPy}, and details are reviewed in \Cref{sec:iwan_details}.

To profile the system and provide inputs for the VPRNM ROM, EPMC backbones for the first two modes for frequency and damping are shown in \Cref{fig:3dof_epmc_freq} and \Cref{fig:3dof_epmc_damp} respectively.
The EPMC backbones show behavior typical of frictional systems with both modes decreasing in frequency (1.12 rad/s to 0.81 rad/s and 3.19 rad/s to 2.80 rad/s for the first and second modes respectively). 
As the system slips, both modes show strong damping nonlinearities in \Cref{fig:3dof_epmc_damp}.

For the first mode backbone, only harmonics 0, 1, and 2 are used so that the EPMC backbone does not include the superharmonic resonance with the second mode. 
This prevents tongues in the backbone due to the internal resonance \cite{kerschenNonlinearNormalModes2009, krackMethodNonlinearModal2013} since the VPRNM ROM captures the internal resonance based on the VPRNM backbone and assumes only single modes contribute to both EPMC backbones.
For the second mode backbone, harmonics 0, 1, 2, and 3 are utilized since there is no superharmonic resonance at a multiple of the second modal frequency.
Harmonics 0, 1, 2, and 3 are utilized for subsequent simulations with HBM and VPRNM. 
Including harmonics 0 and 1-8 did not visibly change any of the results for the 3 DOF system (the third harmonic was deliberately omitted from the fundamental EPMC backbone in this case as well).
For all simulations (HBM, EPMC, and VPRNM) with the 3 DOF system, 1024 time points are used in AFT evaluations. 
Next, the superharmonic resonance behavior of this system is analyzed utilizing HBM and VPRNM.

\FloatBarrier

\subsection{Superharmonic Resonance Case} \label{sec:3dof_SR_results}

To analyze the superharmonic resonance behavior, the system described in \Cref{sec:3dof_description} is excited to have amplitudes of 10 m, 20 m, 30 m, 40 m, 50 m, and 70 m at DOF 1. In HBM, the external force scaling is allowed to vary and the constraint of \eqref{eq:amp_constraint} is applied.
Furthermore, the VPRNM ROM described in \Cref{sec:vprnm_rom} is compared to the results throughout this section.
For the VPRNM ROM, the amplitude calculation vectors are chosen as
\begin{equation}
	\boldsymbol{R_1} = \boldsymbol{R_n} = \begin{bmatrix}
		1& 0 & 0
	\end{bmatrix}
\end{equation}
resulting in the VPRNM ROM being most accurate for DOF 1.

\begin{figure}[h!]
	\centering
	\includegraphics[width=0.9\linewidth]{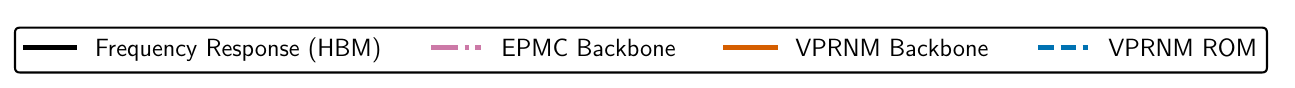}
	\\
	\includegraphics[width=\linewidth]{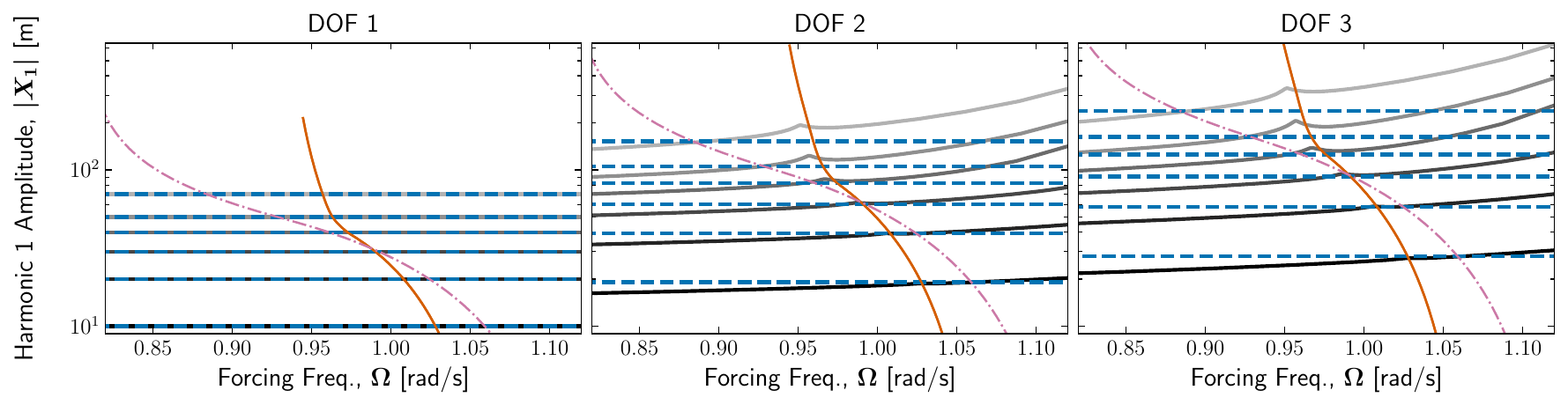}
	\caption{
		First harmonic response magnitude of the 3 DOF system. 
		Lighter gray lines correspond to higher controlled amplitudes at DOF 1 calculated with HBM.
		The proposed VPRNM ROM (dashed) lines approximates HBM truth solution lines (solid gray).
		The y-axes are shared between all three plots.
		The EPMC and VPRNM backbones track the primary and superharmonic resonance frequencies respectively.
	}
	\label{fig:3dof_h1}
\end{figure}

In \Cref{fig:3dof_h1}, DOF 1 shows the controlled amplitudes while DOFs 2 and 3 show some changes in amplitude because multiple modes contribute to the response, especially further from resonance. 
For DOFs 2 and 3, there are slight increases for the higher amplitude curves near the VPRNM backbone due to the superharmonic resonance affecting the first harmonic vibration shape.
The VPRNM ROM does not capture the shifts in the vibration response shape since it is based on a single EPMC mode, but it does match all three DOFs at the primary resonance marked by the EPMC backbone in \Cref{fig:3dof_h1}. 
Other studies \cite{krackReducedOrderModeling2013, krackMethodNonlinearModal2013, schwarzValidationTurbineBlade2020} have utilized nonlinear modes to calculate constant force FRCs and in some cases superimposed responses from multiple modes. Multiple EPMC modes are not utilized here for amplitude control because that would require a solution to a nonlinear system of equations and significantly increase computation time for the VPRNM ROM. Thus the controlled amplitude EPMC ROM (and the fundamental motion of the VPRNM ROM) is only accurate near an isolated primary resonance.

\begin{figure}[h!]
	\centering
	\includegraphics[width=0.75\linewidth]{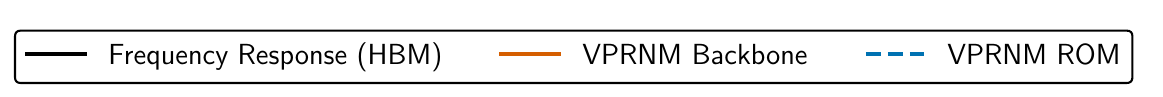}
	\\
	\includegraphics[width=\linewidth]{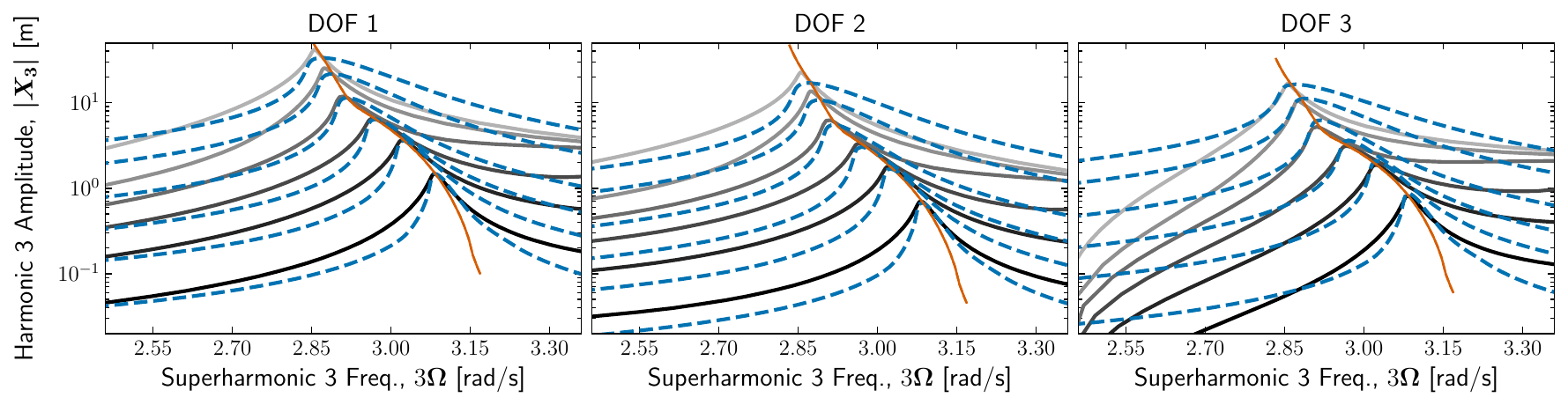}
	\caption{
		Third harmonic response magnitude of the 3 DOF system. 
		Lighter gray lines correspond to higher controlled amplitudes at DOF 1 calculated with HBM (see \Cref{fig:3dof_h1}).
		The proposed VPRNM ROM (dashed) lines approximates HBM truth solution lines (solid gray).
		The y-axes are shared between all three plots.
		The x-axes are the superharmonic resonance frequency and correspond to the same forcing frequency bounds as \Cref{fig:3dof_h1}.
	}
	\label{fig:3dof_h3}
\end{figure}

The vibration of the first harmonic excites the third harmonic through the nonlinear force resulting in superharmonic resonances shown in \Cref{fig:3dof_h3}.
These superharmonic resonances are well tracked by the orange VPRNM line.
For the lower amplitudes, the VPRNM ROM does a good job of matching the superharmonic resonance responses. 
This highlights how amplitude control of the first harmonic response allows for the application of single nonlinear mode theory to understand the superharmonic resonance through the VPRNM ROM.
For the highest two amplitudes, the VPRNM ROM has some clear errors for DOF 1 and 2 due to VPRNM not exactly capturing the peak response amplitude of the superharmonic resonance, likely due to neglecting $ \boldsymbol{F_{kq,sup,n}} $ in the phase constraint (see \Cref{sec:vprnm_theory}). 
Furthermore, at the higher amplitudes the VPRNM ROM shows wider peaks, seemingly overestimating the damping of the superharmonic response. This is a limitation of the proposed VPRNM ROM in that it relies on the damping of the superharmonic mode responding independently to approximate the superharmonic response of two modes acting simultaneously and influencing the damping of each other. If approaches for calculating the modal damping as a function of two modal amplitudes are matured, then such approaches could be used with the VPRNM ROM to more accurately capture the superharmonic resonance peak. 
It is important to note here that alternative efficient modeling techniques cannot capture superharmonic resonances at all. Thus despite the errors in the VPRNM ROM, it is significant improvement over alternatives.

\begin{figure}[h!]
	\centering
	\includegraphics[width=0.7\linewidth]{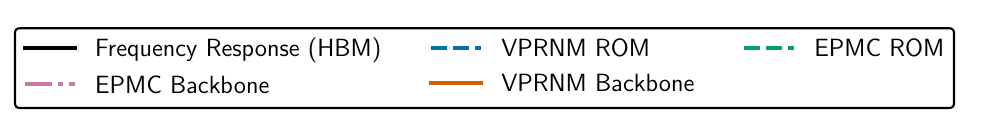}
	\begin{subfigure}[c]{0.45\linewidth } 
		\centering
		\includegraphics[width=\linewidth]{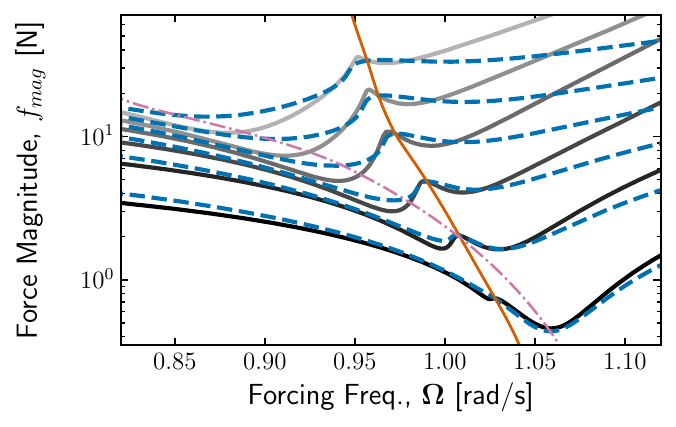}
		\caption{}
		\label{fig:3dof_VPRNMROM_F}
	\end{subfigure}
\quad
	\begin{subfigure}[c]{0.45\linewidth } 
		\centering
		\includegraphics[width=\linewidth]{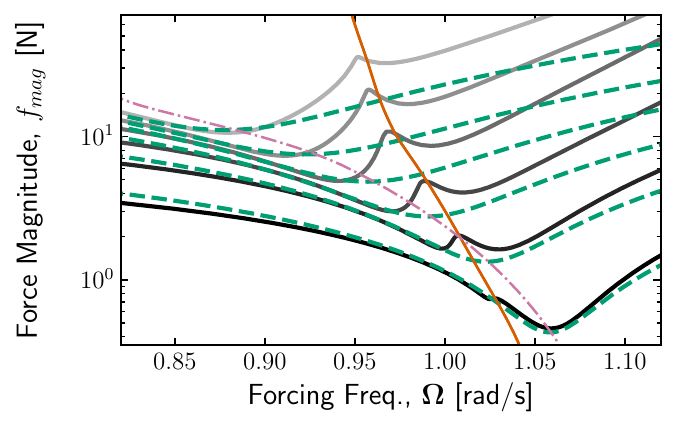}
		\caption{}
	\label{fig:3dof_EPMC_F}
	\end{subfigure}
	\caption{
	External force scaling for the 3 DOF system. 
	Lighter gray lines correspond to higher controlled amplitudes at DOF 1 calculated with HBM (see \Cref{fig:3dof_h1}).
	The HBM truth solution (solid gray) is better matched by (a) the VPRNM ROM (dashed) than (b) the EPMC ROM (dashed) without superharmonic resonance correction. The EPMC and VPRNM backbones track the primary and superharmonic resonance frequencies respectively.
	}
	\label{fig:3dof_force}
\end{figure}

Since the response amplitude is controlled, the required external force scaling is shown in \Cref{fig:3dof_force}.
The reconstructed force for the VPRNM ROM in \Cref{fig:3dof_VPRNMROM_F} is reasonably accurate near the superharmonic resonance (the orange VPRNM line) and for wider frequency ranges at lower amplitude and force levels. 
On the other hand, the EPMC ROM\footnote{As previously mentioned, the EPMC backbone is calculated without the third harmonic as superharmonic resonances introduce artificial effects in EPMC \cite{krackNonlinearModalAnalysis2015}.} based on amplitude control in \Cref{fig:3dof_EPMC_F} misses the peaks in force at the superharmonic resonance.
At higher frequencies, both ROMs are consistent with each other and differ significantly from the HBM truth solutions as the second mode (not considered for the force in either ROM) starts to contribute significantly to the fundamental motion. 

\FloatBarrier

\begin{figure}[h!]
	\centering
	\includegraphics[width=0.9\linewidth]{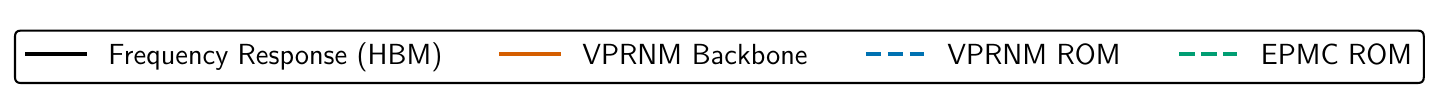}
	\\
	\begin{subfigure}[c]{\linewidth } 
		\centering
		\includegraphics[width=\linewidth]{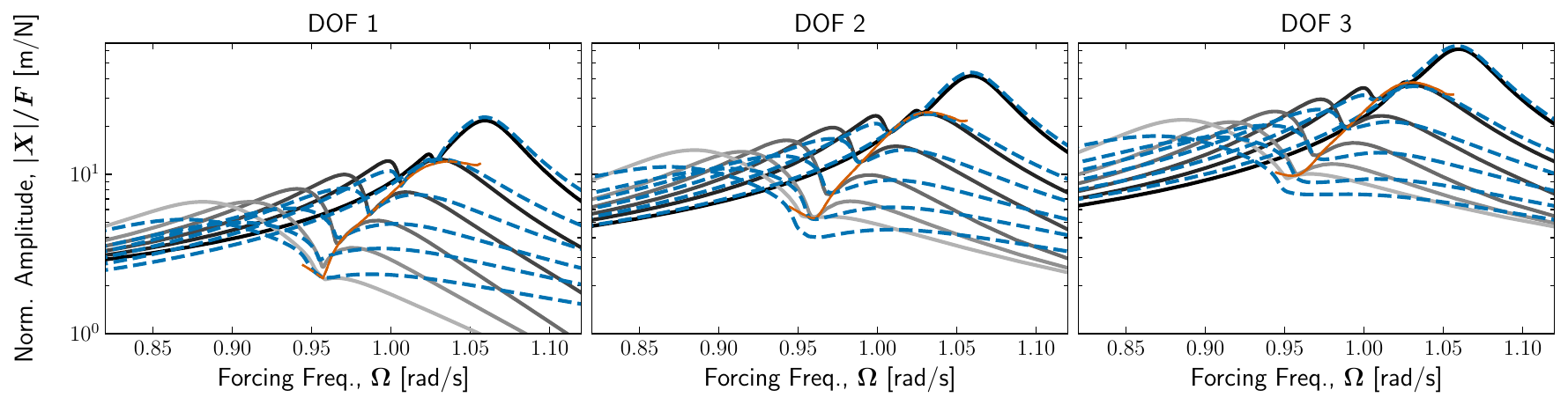}
		\caption{}
		\label{fig:3dof_VPRNMROM_tot}
	\end{subfigure}
\\
	\begin{subfigure}[c]{\linewidth } 
		\centering
		\includegraphics[width=\linewidth]{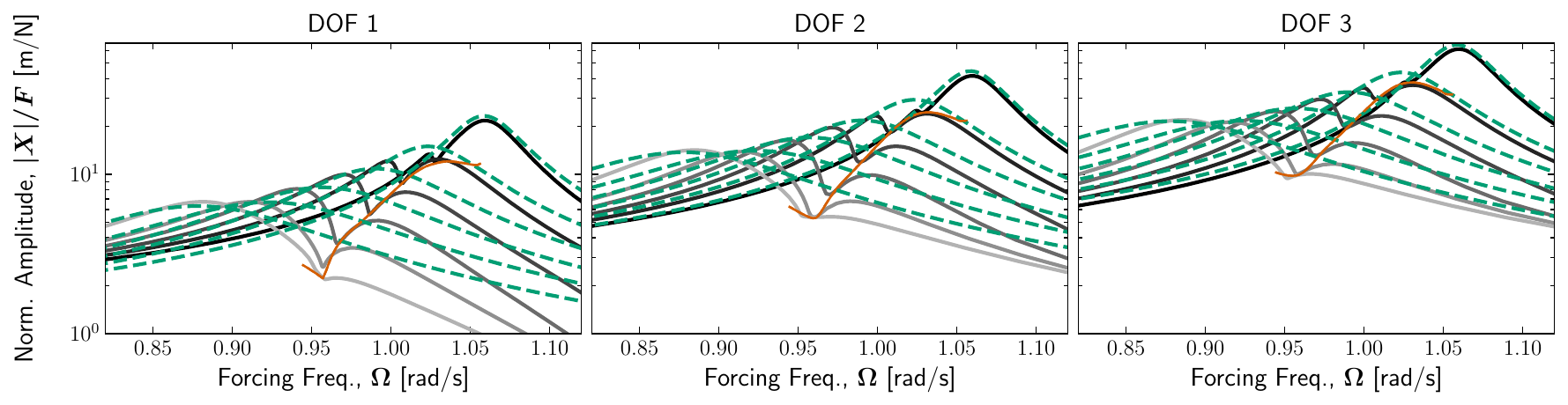}
		\caption{}
		\label{fig:3dof_EPMC_tot}
	\end{subfigure}
	\caption{
	Total response magnitude divided by force level of the 3 DOF system. 
	Lighter gray lines correspond to higher controlled amplitudes at DOF 1 calculated with HBM (see \Cref{fig:3dof_h1}).
	The HBM truth solution (solid gray) is better matched by (a) the VPRNM ROM (dashed) than (b) the EPMC ROM (dashed) without superharmonic resonance.
	The y-axes are shared between all three plots in both (a) and (b).
	}
	\label{fig:3dof_total}
\end{figure}

The total response considering the highest amplitude at any point in a cycle divided by the force scaling at the given forcing frequency is shown in \Cref{fig:3dof_total}. 
For the normalized plots, the highest control amplitudes (lightest colors) correspond to the most damping and the lowest normalized response amplitudes. 
As before, the VPRNM backbone does a good job of tracking the superharmonic resonance behavior, which corresponds to a local minimum due to the phase difference between the first and third harmonics. 
Likewise, the VPRNM ROM captures the FRC behavior of the superharmonic resonance in \Cref{fig:3dof_VPRNMROM_tot} especially for DOF 1. 
While qualitatively correct, there are clear errors for DOFs 2 and 3 at the highest control amplitudes (lightest curves) because the choice of the control location in the VPRNM ROM is DOF 1 ($\boldsymbol{R_1}$ and $\boldsymbol{R_n}$ in \Cref{sec:vprnm_rom}). If a different location was of primary interest, then the control location could be changed in the VPRNM ROM to provide more accurate results for a different DOF. 
The VPRNM ROM can also be compared to the state of the art presented by the EPMC ROM in \Cref{fig:3dof_EPMC_tot}. 
The EPMC ROM does not capture any of the superharmonic resonance behavior highlighting that the VPRNM ROM is a significant improvement for computationally efficient approaches despite some errors compared to the truth solution of HBM.
As previously mentioned, EPMC is calculated without the third harmonic (filtering out the superharmonic resonance) because including the superharmonic resonance in the EPMC backbone results in artificial modal interactions \cite{krackNonlinearModalAnalysis2015} and does not allow the EPMC ROM to capture the superharmonic resonance in FRCs.

Away from the primary resonance, the VPRNM ROM shows some clear errors relative to the HBM truth solution. The absolute errors are exaggerated in \Cref{fig:3dof_VPRNMROM_tot} due to the log scale, and the errors can be attributed to the assumption of a single mode for the primary response. Similar errors are seen in \Cref{fig:3dof_EPMC_tot}, which represents the current state of the art for nonlinear mode based ROMs.

\FloatBarrier

\begin{figure}[h!]
	\centering
	\includegraphics[width=0.75\linewidth]{hbm_vprnm_rom_legend.pdf}
	\\
	\includegraphics[width=\linewidth]{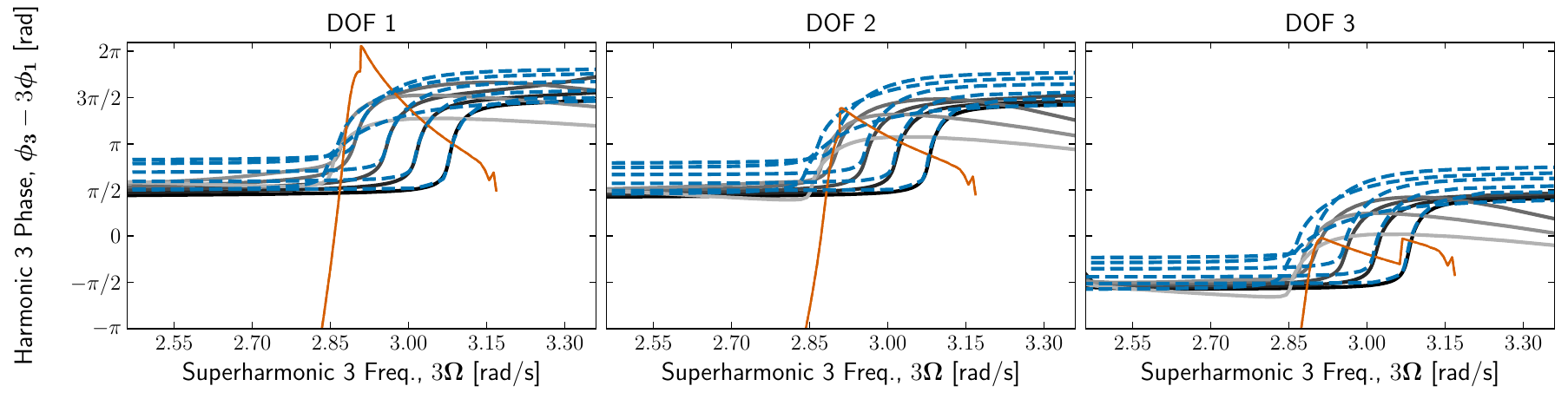}
	\caption{
	Difference in phase between third and first harmonics of the 3 DOF system. 
	Lighter gray lines correspond to higher controlled amplitudes at DOF 1 calculated with HBM (see \Cref{fig:3dof_h1}).
	The proposed VPRNM ROM (dashed) lines approximates HBM truth solution lines (solid gray).
	The y-axes are shared between all three plots.
	The x-axes are the superharmonic resonance frequency and correspond to the same forcing frequency bounds as \Cref{fig:3dof_h1}.
	}
	\label{fig:3dof_h3phase}
\end{figure}

The relative phase between the third and the first harmonics is shown in \Cref{fig:3dof_h3phase} since the VPRNM backbone is based on a phase resonance criteria. Here the VPRNM ROM shows a clear shift of $ \pi $ as the superharmonic resonance is crossed and matches the HBM phase at the VPRNM backbone for DOF 1.
Of note, the phase for the highest amplitude case (lightest curve, 70 m control amplitude) shows a phase shift of only $ \pi/2 $ rather than $ \pi $ for DOF 2. 
This can be attributed to other vibration motion contributing to the response making it harder to clearly identify the superharmonic resonance, and likely contributing to previously discussed errors.
Utilizing the modal filter with VPRNM that is proposed for the high DOF system in \Cref{sec:theory_high_dim} did not notably change the results in this section.
The sharp corners in the VPRNM backbones in \Cref{fig:3dof_h3phase} are attributed to the discrete implementation of the Iwan element, which shows non-smooth behavior at low amplitudes when only a few sliders transition between stick and slip behaviors.

\begin{figure}[h!]
	\centering
	\includegraphics[width=0.5\linewidth]{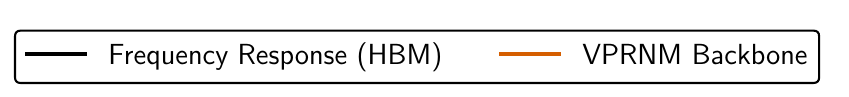}
	\\
	\includegraphics[width=\linewidth]{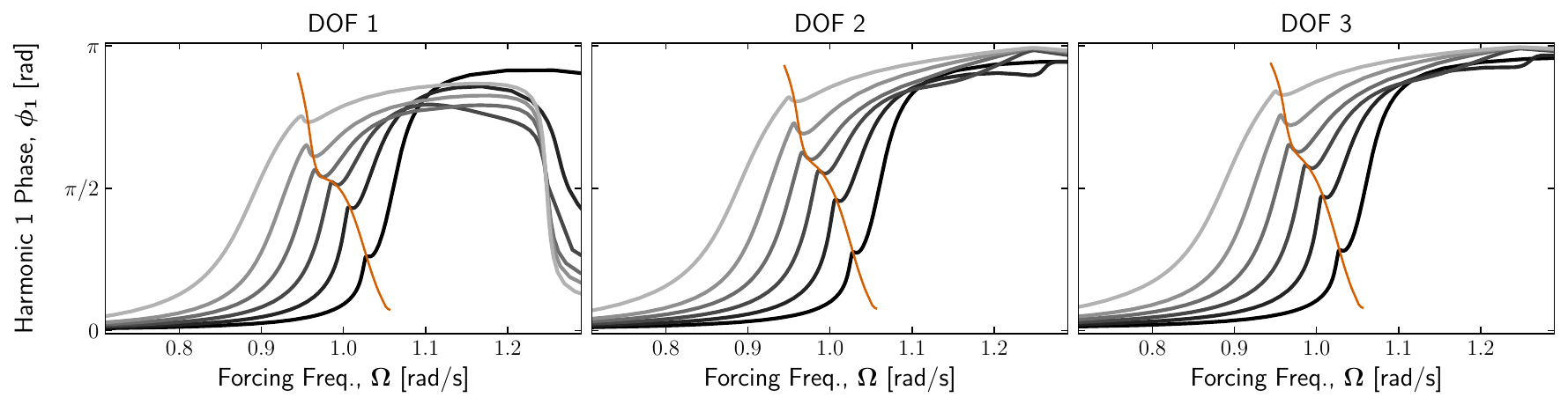}
	\caption{
	First harmonic phase of the 3 DOF system. 
	Lighter gray lines correspond to higher controlled amplitudes at DOF 1 calculated with HBM.
	The VPRNM backbone (orange) tracks distortion in the phase shift due to the superharmonic resonance.
	The y-axes are shared between all three plots.
	}
	\label{fig:3dof_h1phase}
\end{figure}

The phase of the first harmonic shows a typical phase shift of $ \pi $  in \Cref{fig:3dof_h1phase} as well as some phase distortions near the superharmonic resonance in the HBM results. The present formulations of the VPRNM and EPMC ROMs do not calculate a phase difference between the first harmonic response and the forcing magnitude, so no VPRNM or EPMC ROM results are presented in \Cref{fig:3dof_h1phase}.

\FloatBarrier

\subsection{No Superharmonic Resonance Case} \label{sec:3dof_noSR}

\FloatBarrier

\begin{figure}[!h]
	\centering
	\begin{subfigure}[c]{0.45\linewidth } 
		\centering
		\includegraphics[width=\linewidth]{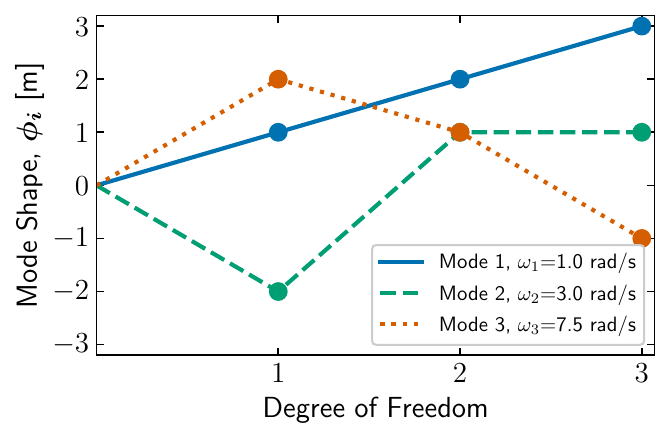}
		\caption{}
		\label{fig:noSR_modeshapes}
	\end{subfigure}	\\
	\begin{subfigure}[c]{0.45\linewidth } 
		\centering
		\includegraphics[width=\linewidth]{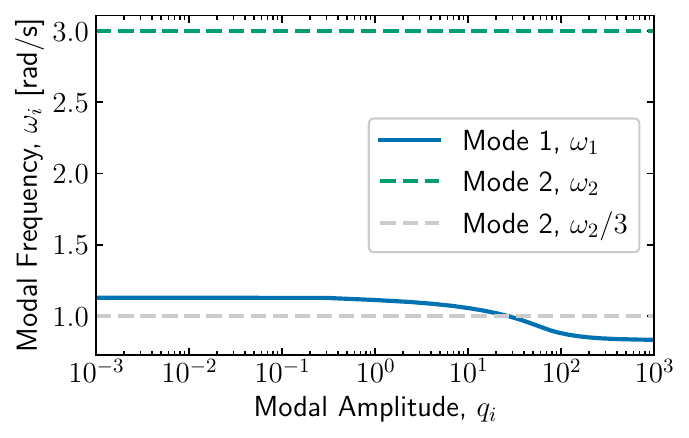}
		\caption{}
	\end{subfigure}	\quad
	\begin{subfigure}[c]{0.45\linewidth } 
		\centering
		\includegraphics[width=\linewidth]{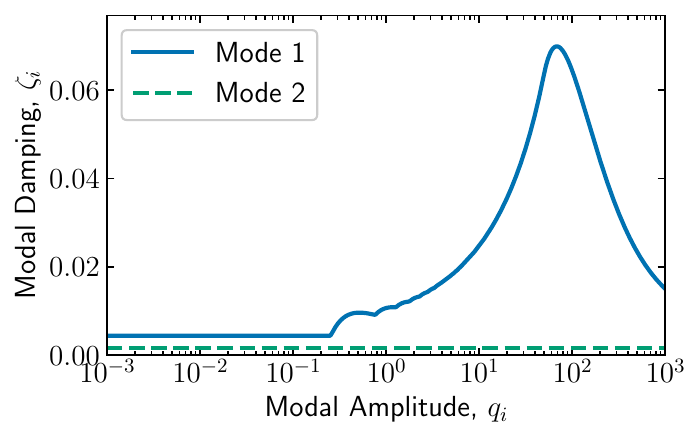}
		\caption{}
	\end{subfigure}		\caption{
	Overview of the modified 3 DOF system: (a) mass-normalized mode shapes around a linear state including half of the linearized frictional force (representative of partial slip), (b) modal frequency from EPMC, and (c) modal damping from EPMC.
	} 
	\label{fig:noSR_EPMC}
\end{figure}

A second case with a modified 3 DOF system is investigated without a superharmonic resonance. For this case, the three mode shapes are the same as in the previous section, but the second and third mode shapes are switched (see \Cref{fig:noSR_modeshapes}). Therefore, the second mode at a frequency of $3.0 $ rad/s has zero displacement over the nonlinear element (equal displacement of DOFs 2 and 3). 
The EPMC results for the first two modes of the modified system are shown in \Cref{fig:noSR_EPMC}. Here, since the nonlinearity does not act on the second mode, it remains linear at a constant frequency and damping value.

\begin{figure}[h!]
	\centering
	\begin{subfigure}[c]{\linewidth } 
		\centering
		\includegraphics[width=\linewidth]{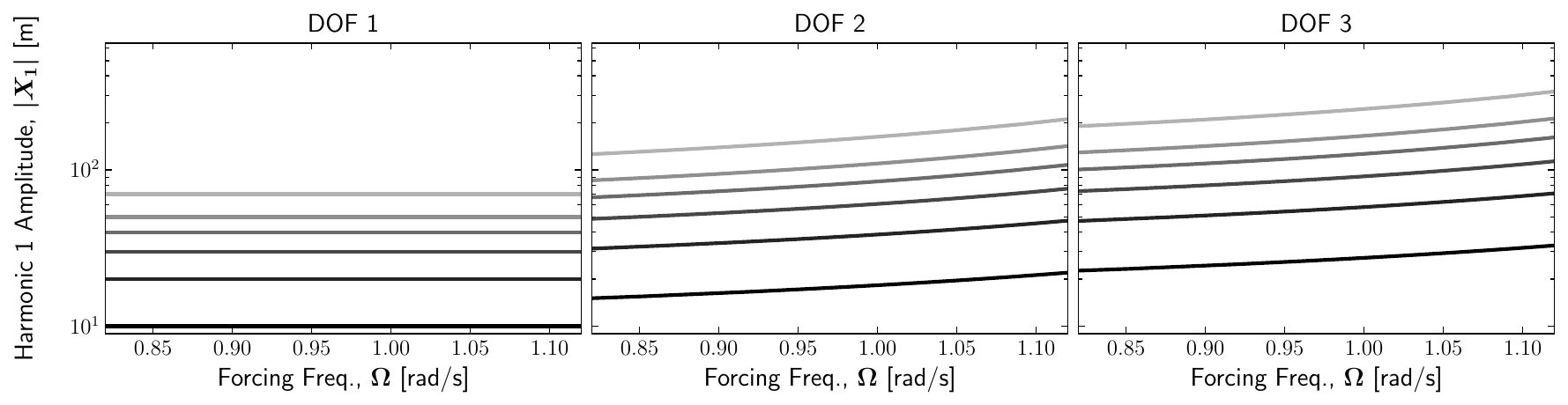}
		\caption{}
		\label{fig:noSR_H1}
	\end{subfigure}
	\\
	\begin{subfigure}[c]{\linewidth } 
		\centering
		\includegraphics[width=\linewidth]{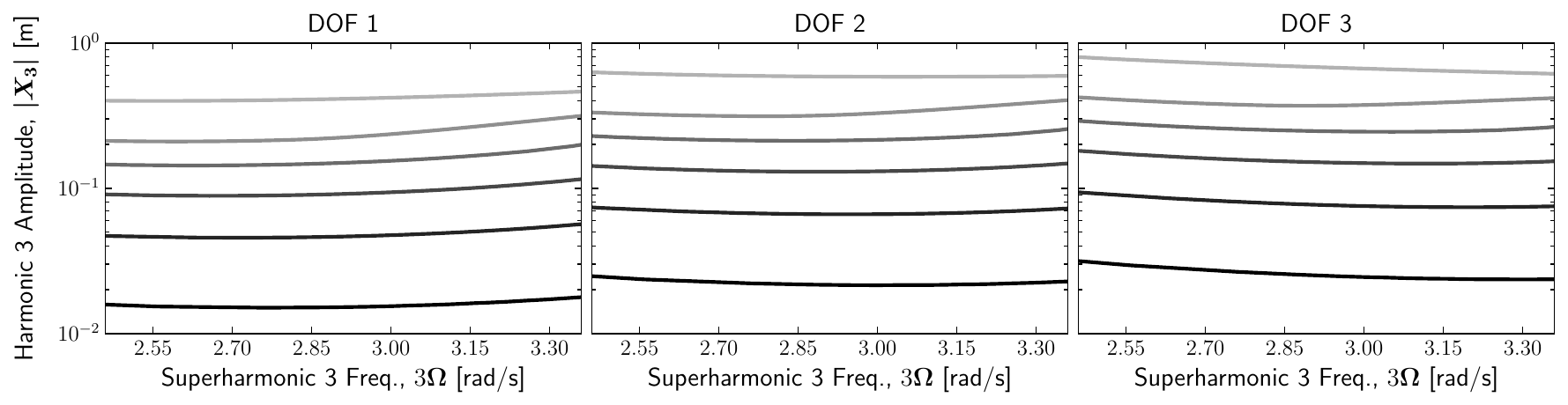}
		\caption{}
		\label{fig:noSR_H3}
	\end{subfigure}
	\\
	\begin{subfigure}[c]{\linewidth } 
		\centering
		\includegraphics[width=\linewidth]{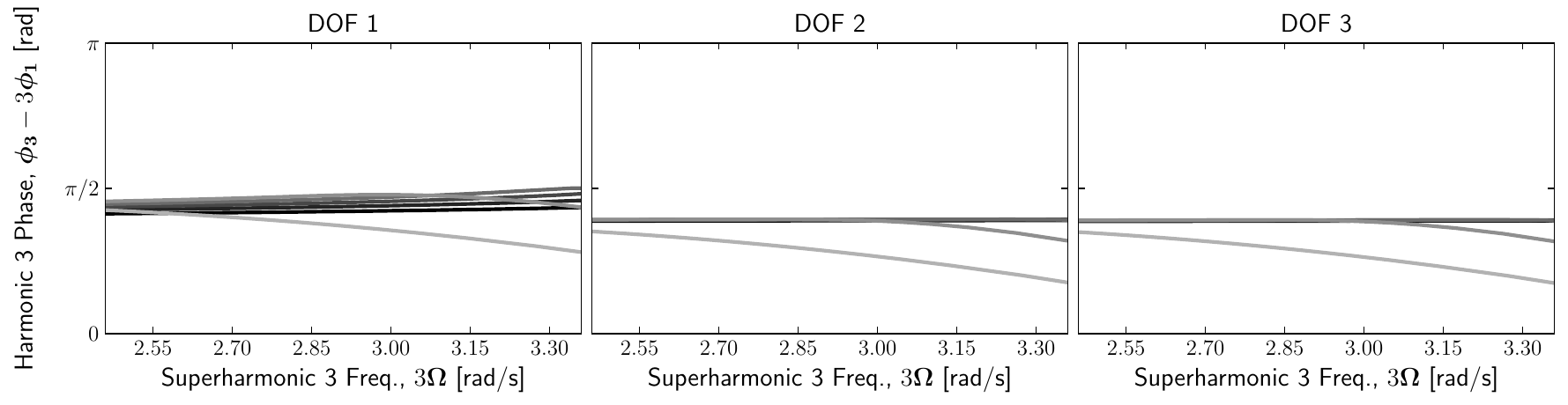}
		\caption{}
		\label{fig:noSR_H3Phase}
	\end{subfigure}
	\caption{
	Responses for the modified 3 DOF system without superharmonic resonance: (a) first harmonic response magnitude,
	(b) third harmonic response magnitude,
	and (c) difference in phase between third and first harmonics.
	Note that the plot bounds in (a) and (b) are different the the superharmonic response amplitude is consistently at least two orders of magnitude smaller than the fundamental response.
	Lighter gray lines correspond to higher controlled amplitudes at DOF 1 calculated with HBM as given by (a).
	The y-axes are shared between all three plots in each subplot.
	The forcing frequency bounds are the same for each subplot, but (b) and (c) plot the superharmonic response frequency on the x-axis corresponding to three times the forcing frequency.
	}
	\label{fig:noSR_amp}
\end{figure}

This system is excited harmonically and the amplitude of the first harmonic at DOF 1 is again controlled to be 10 m, 20 m, 30 m, 40 m, 50 m, and 70 m (see \Cref{fig:noSR_H1}).  In this case, no significant response or clear peak in the third superharmonic is observed in \Cref{fig:noSR_H3} (corresponding to approximately 3 rad/s and the natural frequency of the second mode).
The lack of a supharmonic resonance is confirmed by investigating the phase of the third harmonic in \Cref{fig:noSR_H3Phase}, which does not show a shift by $\pi$ that would indicate a resonance phenomena.
Since no superharmonic resonance is observed the VPRNM backbone and VPRNM ROM are not considered.

Conceptually, the lack of a superharmonic resonance fits with the formulation of VPRNM. For this system $\boldsymbol{F_{nq,broad}}$ from \eqref{eq:fbroad} acts equally and opposite on DOFs 2 and 3. This is interpreted as external excitation of a potential superharmonic resonance of mode 2. Converting $\boldsymbol{F_{nq,broad}}$ to modal forcing results in $ \boldsymbol{\psi_2}^T \boldsymbol{F_{nc,broad}} = 0$ and $\boldsymbol{\psi_2}^T \boldsymbol{F_{ns,broad}} = 0$ because the linear mode shape $\boldsymbol{\psi_2}$ has equal displacements at DOFs 2 and 3.
Therefore, there is no modal forcing for a superharmonic resonance of mode 2 and no superharmonic resonances of mode 2 are expected. 
This highlights that the framework of VPRNM can be utilized to understand the requirements for a superharmonic resonance beyond the necessary but not sufficient requirement of an integer ratio between the forcing and response frequencies.
Next, computation times for the 3 DOF system with a superharmonic resonance are discussed. 

\FloatBarrier

\subsection{Computation Time}
\FloatBarrier

The computation time for the VPRNM ROM is compared to HBM for the 3 DOF system to illustrate the computational efficiency of the VPRNM ROM in \Cref{tab:3dof_time}.
Considering both the VPRNM ROM construction (EPMC and VPRNM calculations) and evaluation, the VPRNM ROM represents a speedup of 1.75 times compared to the truth solution of HBM at six amplitude levels. 
Once the VPRNM ROM is constructed, responses at additional amplitude levels can be constructed 15,900 times faster than HBM solutions. 
Furthermore, the VPRNM ROM construction provides additional information not captured in the HBM solution such as allowing for EPMC ROM based construction of primary resonances responses near the second mode, which are not considered here.

For timing, the VPRNM backbone is only run between external forces of 0.4 N and 60 N (rather than 0.1 N and 100 N used for the plots) since this is sufficient to generate the VPRNM ROM over the amplitude regime considered with HBM. Likewise, the maximum mass normalized modal amplitude of the EPMC backbones for the first and second modes are reduced to $ q_1=100 $ and $ q_2=31.62$.
Updating these bounds does not affect the quality of the VPRNM ROM results shown here and gives a more fair comparison to HBM because only the necessary inputs to the VPRNM ROM are calculated.

\begin{table}[h!]
	\centering
	\caption{Simulation times for the 3 DOF system with superharmonic resonance on a desktop computer (6 core Intel i7-10710U CPU, 32 GB of RAM, 1.10 GHz processor). HBM, VPRNM, and EPMC computation times are averaged over 10 calculations. VPRNM ROM computation time is averaged over 10,000 calculations.}
	\label{tab:3dof_time}
	\begin{tabular}{ccc}
		\hline \hline 
		Method & Amplitude Levels & Computation Time (sec) \\ \hline \hline
		\multicolumn{3}{l}{Backbone Calculations} \\
		\hline
		VPRNM & NA & 15.1 \\
		EPMC Mode 1 & NA & 33.7 \\
		EPMC Mode 2 & NA & 27.0 \\ 
		VPRNM ROM Construction & NA & 75.8 \\ 
		\hline \hline
		\multicolumn{3}{l}{FRC Calculations} \\ \hline
		HBM & 6 & 133 \\ 
		VPRNM ROM Evaluation & 6 & 8.39e-03 \\
		\hline \hline 
		
	\end{tabular}
\end{table}

\FloatBarrier

\FloatBarrier

\section{Half Brake-Reu{\ss} Beam System} \label{sec:hbrb_system}

The Half Brake-Reu{\ss} Beam (HBRB) is a variant of the Brake-Reu{\ss} Beam (BRB) with half of the nominal beam thickness that was previously tested in \cite{chenMeasurementIdentificationNonlinear2022}. 
The HBRB has a three bolt lap joint with 7.94 mm (5/16 inches) diameter bolts and a contact interface that is nominally identical to the BRB, a common benchmark for studying nonlinear vibration of jointed structures \cite{brakeMechanicsJointedStructures2017, brakeObservationsVariabilityRepeatability2019}.
In a previous experimental study \cite{chenMeasurementIdentificationNonlinear2022}, the HBRB exhibited a superharmonic and internal resonance
when excited near the first bending mode resulting in large responses of the third bending mode at seven times the forcing frequency.
This section details the experimental test procedure (\Cref{sec:exper_procedure}) and modeling approach (\Cref{sec:model_procedure}).

\subsection{Experimental Procedure} \label{sec:exper_procedure}

This work considers similar HBRB tests to \cite{chenMeasurementIdentificationNonlinear2022} with the key improvement that instrumented bolts are used to more precisely control the bolt preload values. 
Relationships between bolt torque and bolt tension show significant variability \cite{ruan_variability_2019}, so recent modeling efforts have relied on instrumented bolts to more accurately characterize bolt preload \cite{balajiTractionbasedMultiscaleNonlinear2020, balajiDissipative2021, porterQuantitativeAssessmentModel2022, porterPredictivePhysicsbasedFriction2023}.
Three bolt tensions of 12006 N (high), 9069 N (medium), and 6135 N (low) are utilized for the experiments. All tests at each bolt tension are conducted without disassembling the structure on the same day and bolt tensions are checked after initial random excitation (to allow for the joint to settle) and after all tests are completed (to verify no significant loss after settling).
After settling, bolt tensions checked throughout the tests showed variations of approximately 100 N at the high preload, 60 N at the medium preload, and 10 N at the low preload.
Assembly of the HBRB follows best practices\footnote{The HBRB is clamped to a solid beam; a shim is used to align the bolt holes. Then the center bolt is tightened to 70\% of the preload; the two outer bolts are tightened to 70\% of the preload. Lastly, all bolts are tightened to approximately 105\% of the preload in the same order (to account for some preload being lost during settling).}
established in \cite{brakeObservationsVariabilityRepeatability2019}.
Three washers are used on each side of the bolts to protect the strain gauges as the instrumented bolts are designed for the thicker BRB.

\begin{figure}[h!]
	\centering
	\includegraphics[width=0.9\linewidth]{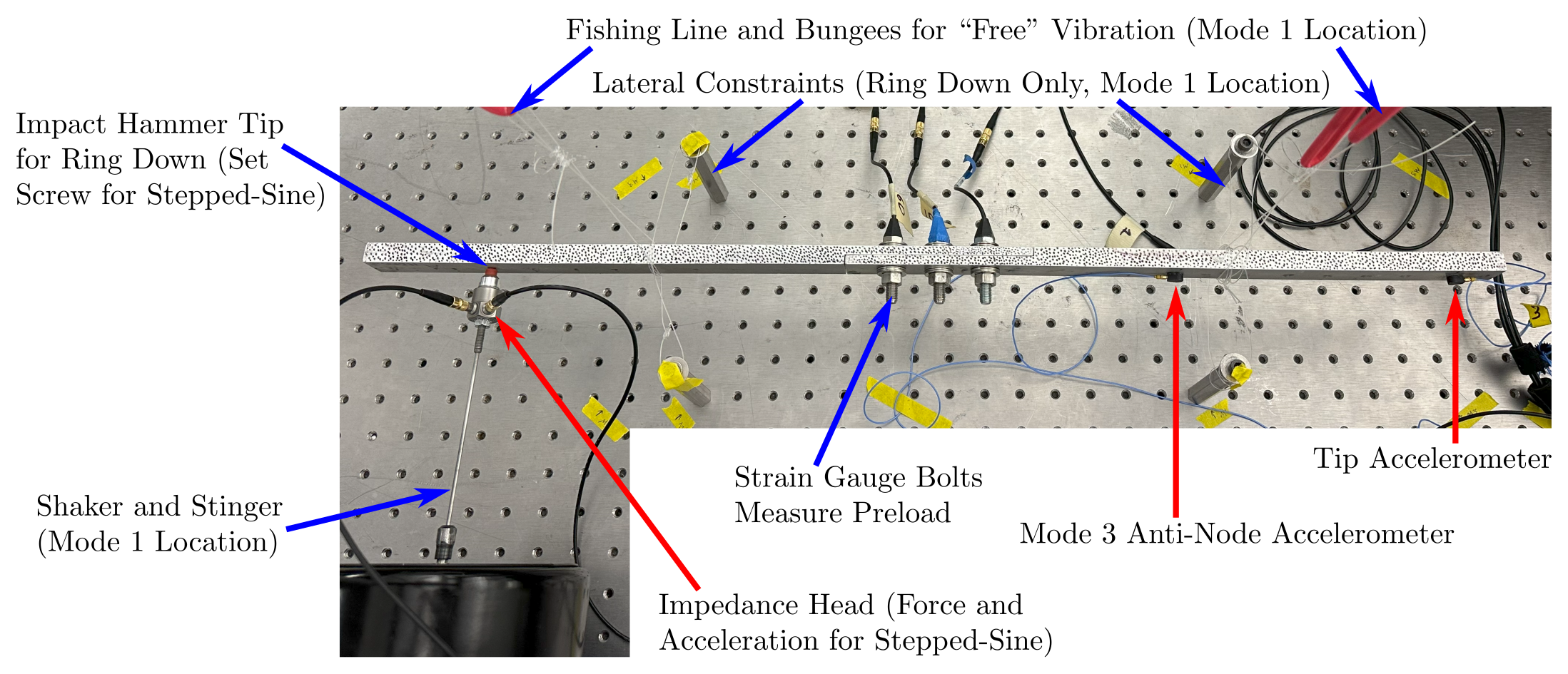}
	\caption{
		Experimental setup for shaker ring down of the first bending mode of the Half Brake-Reu{\ss} Beam. 
		For stepped-sine testing, no lateral constraints are used and the impedance head is screwed into to the structure. 
		The red arrows denote the locations of accelerometers used for reporting experimental results.
		For the third bending mode, the supports and lateral constraints are moved towards the ends of the beam and the shaker is moved towards the center based on \cite{chenMeasurementIdentificationNonlinear2022}.}
	\label{fig:ringdown_setup}
\end{figure}

Experimental tests are conducted to obtain the response to harmonic excitation (stepped-sine), the nonlinear modal properties of the first and third bending modes (shaker ring down), and linear properties of the unassembled beams (hammer impact). The setup for shaker ring down is shown in \Cref{fig:ringdown_setup}.
For all tests, uniaxial accelerometers (PZB Piezotronics Model 352A24) are used at 30 mm (1.2 inches) and 205 mm (8.1 inches) from the right end of the beam as shown in \Cref{fig:ringdown_setup}.
For subsequent discussion, these accelerometer locations are termed ``Tip'' and ``Mode 3 Anti-Node'' respectively.\footnote{Note, that the mode 3 anti-node is closer to 178 mm (7 inches) from the end \cite{chenMeasurementIdentificationNonlinear2022}, but the accelerometer cannot be placed there because of bungee supports.}
Additionally, for stepped-sine tests, the force and acceleration at the shaker is measured with an impedance head (PZB Piezotronics Model 288D01).
For all tests, the assembled HBRB is suspended from fishing line and bungee cords to create ``free'' boundary conditions at the node points of the mode of interest (178 mm (7 inches) or 76 mm (3 inches) from each end for modes 1 and 3 respectively) based on \cite{chenMeasurementIdentificationNonlinear2022}.

For stepped-sine tests, a Siemens LMS controller is utilized to force the structure at accelerations of 0.5 g, 1.0 g, 2.0 g, 4.0 g, 6.0 g, and 1.0 g (in order) at the shaker location. For each acceleration level, tests are conducted first sweeping from low to high and then high to low forcing frequency. For the highest preload value, tests are conducted from 82 to 85 Hz with steps of 0.05 Hz. For the medium and low preload levels, tests are conducted from 81.5 to 85 Hz with steps of 0.05 Hz to ensure that the resonance frequencies are captured.
For stepped sine testing, no lateral constraints are used and the shaker is screwed into the structure 76 mm (3 inches) from the left side of the HBRB, which corresponds to a node point for the third bending mode.

To obtain nonlinear modal frequency and damping trends, shaker ring down tests are conducted as described in \cite{chenMeasurementIdentificationNonlinear2022} and shown in \Cref{fig:ringdown_setup}. 
For shaker ring down, the HBRB is constrained from lateral rigid body motion by fishing line at approximately 178 mm (7 inches) and 76 mm (3 inches) from each end for the first and third bending modes respectively (corresponding to nodal points \cite{chenMeasurementIdentificationNonlinear2022}).
Then the shaker is preloaded against the structure using an impact hammer tip for the contact point. 
The shaker is located at 76 mm (3 inches) and 178 mm (7 inches) from the left end of the beam for the first and third bending modes respectively. 
During the test, the shaker initially applies single harmonic excitation to the structure at the natural frequency of the mode of interest (as identified by the stepped-sine tests) to achieve an amplitude of approximately 10 g at the tip accelerometer. 
Once a steady-state response has been reached, a DC voltage is sent to the shaker causing it to retract from the structure and the free decay vibration behavior is measured. 
Since the shaker is preloaded against the structure but not directly attached, disturbances caused by pulling the shaker back are small.
Shaker ring down provides better experimental data to compare against nonlinear modal analysis because a single mode dominates the response unlike hammer impact tests that excite all modes simultaneously.

\FloatBarrier

Lastly, hammer impact tests are conducted on the unassembled halves of the HBRB to obtain properties for the linear portions of the model. First modal frequencies were identified from the free decay of five hammer impact tests of each half as 305.0 Hz and 300.9 Hz for Beams A and B respectively.\footnote{Tests used a sampling frequency of 16 kHz and recorded 128,000 samples including the hammer impact. Using only the free decay portion of the signal, the frequency resolution of the Fourier transform was 0.126 Hz for all 10 identifications. The total range across all five hits of a given half of the structure was less than the 0.1 Hz.} 
The average density of the unassembled beams (using the measured masses and nominal volumes) was 7850 kg/m$^3$.
The average frequency between the two unassembled values of 302.9 Hz is used with the finite element model of the unassembled half beams (see \Cref{sec:model_procedure}) to obtain an elastic modulus of 173.4 GPa.
Other parameters are chosen to match  \cite{porterPredictivePhysicsbasedFriction2023}, including a Poisson's ratio of 0.3 \cite{balajiDissipative2021}, a yield strength of 330 MPa \cite{SNL_Materials304L} and a hardening modulus of 620 MPa \cite{SNL_Materials304L}.

\FloatBarrier

\subsection{Modeling Procedure} \label{sec:model_procedure}

\begin{figure}[!h]
	\centering
	\begin{subfigure}[c]{0.48\linewidth } 
		\centering
				\includegraphics[width=\linewidth, trim={0 1.5cm 0 0},clip]{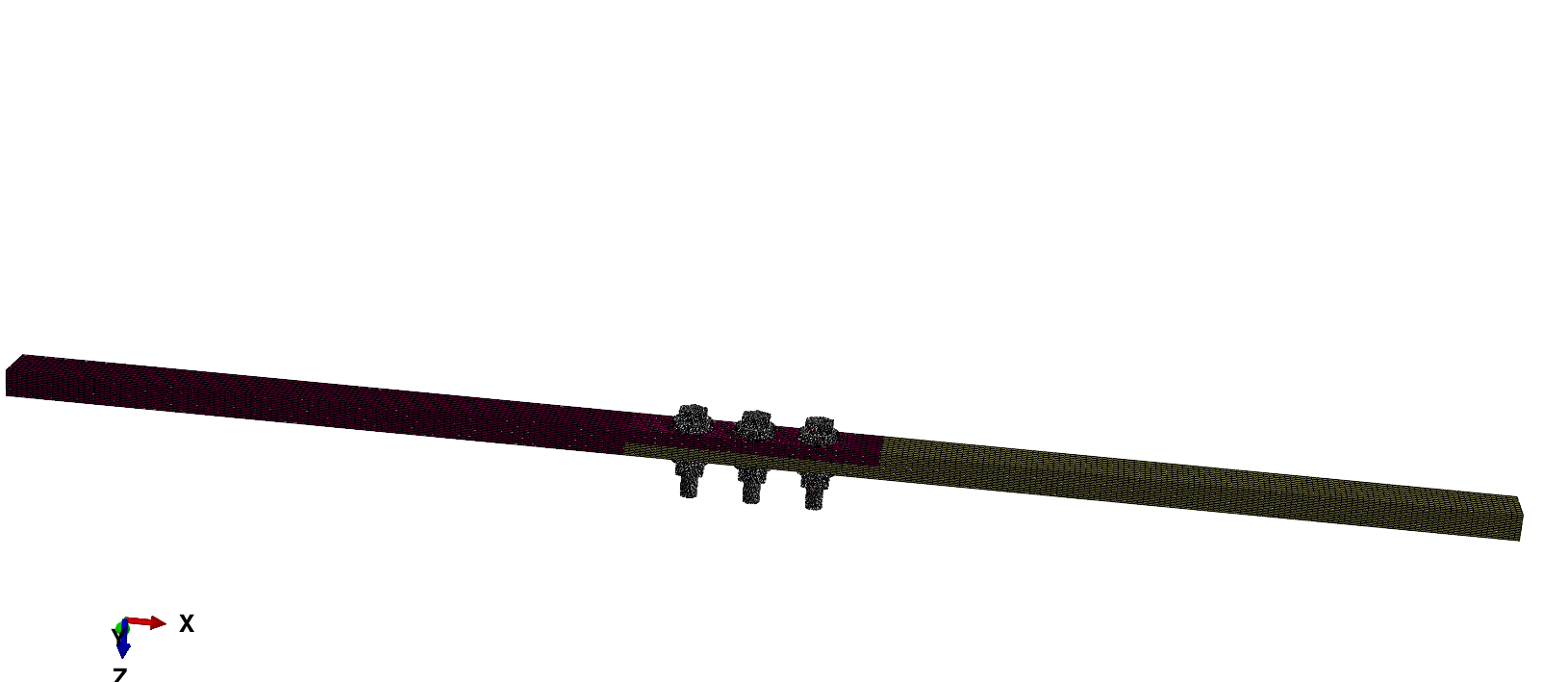}
		\caption{}
			\end{subfigure}	\hfill
	\begin{subfigure}[c]{0.48\linewidth } 
		\centering
		\includegraphics[width=\linewidth, trim={0 3cm 0 0},clip]{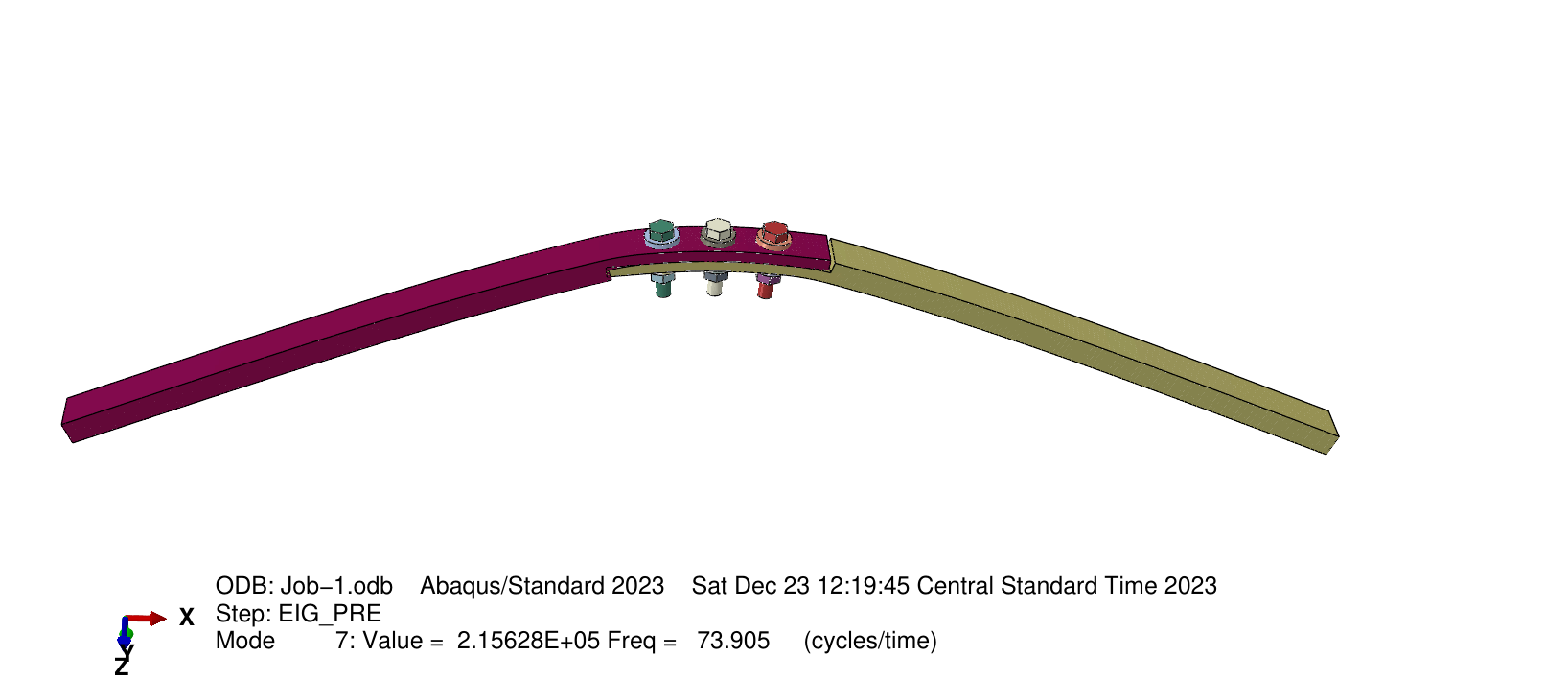}
		\caption{}
			\end{subfigure}	\\
	\begin{subfigure}[c]{0.48\linewidth } 
		\centering
		\includegraphics[width=\linewidth, trim={0 3cm 0 0},clip]{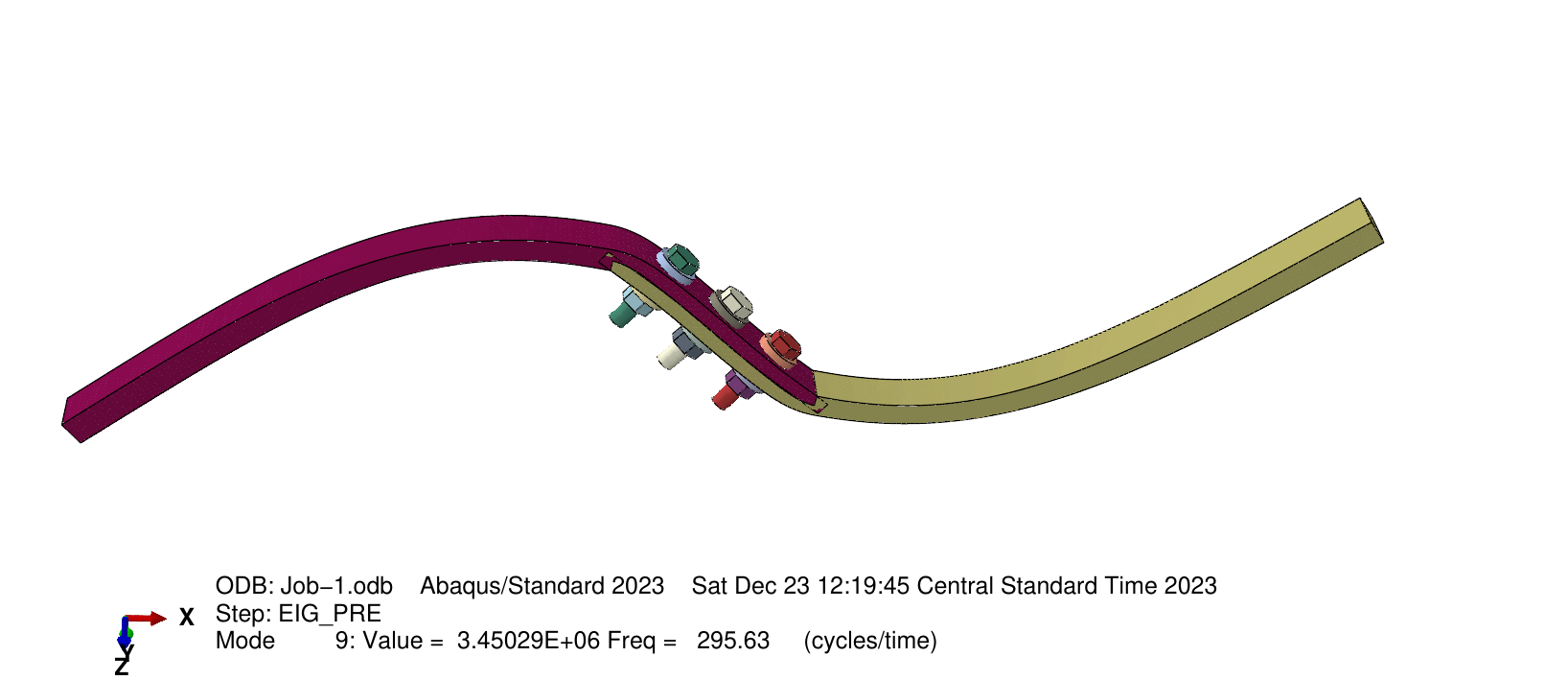}
		\caption{}
			\end{subfigure}	\hfill
	\begin{subfigure}[c]{0.48\linewidth } 
		\centering
		\includegraphics[width=\linewidth, trim={0 3cm 0 0},clip]{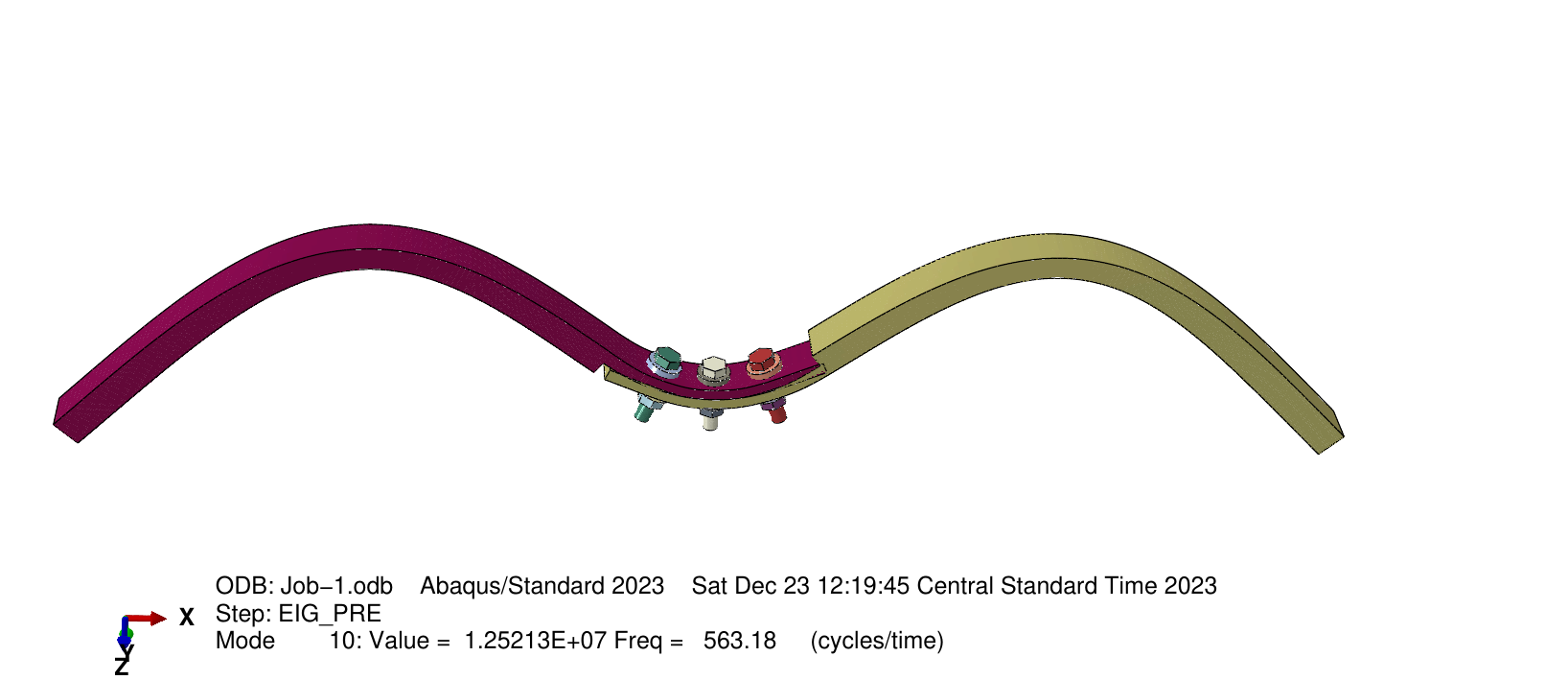}
		\caption{}
			\end{subfigure}		\caption{Abaqus HBRB model of (a) undeformed mesh, (a) first bending mode, (b) second bending mode, and (c) third bending mode.} \label{fig:hbrb_abaqus}
\end{figure}

The HBRB is modeled in Abaqus utilizing the procedure developed in \cite{balajiTractionbasedMultiscaleNonlinear2020} (also see \cite{balajiABAQUSMATLABTutorialJointed2024}). The Abaqus model (see \Cref{fig:hbrb_abaqus}) utilizes a structured interface mesh with the same arrangement of 588 elements as previous studies with the BRB \cite{balajiTractionbasedMultiscaleNonlinear2020, balajiReducedOrderModeling2021}. Based on a convergence study, 16 elements are used through the bending thickness direction of the beam (nominal 12.7 mm thickness) rather than 10 in previous studies with the BRB (nominal 25.4 mm thickness). Furthermore, smaller elements are used for the solid portions of the beam away from the interface and for the nuts, bolts, and washers. 
The full Abaqus model contains 124370 C3D8R elements and 124370 nodes, each with 3 degrees of freedom.
Average masses of the beams, bolts (including strain gauge cables), nuts, and washers are matched between the model and physical measurements. 
Two masses of 0.8 grams, corresponding to the two accelerometers, are added as point masses to the final model.

\begin{figure}[!h]
	\centering
	\begin{subfigure}[c]{0.49\linewidth } 
		\centering
		\includegraphics[width=\linewidth]{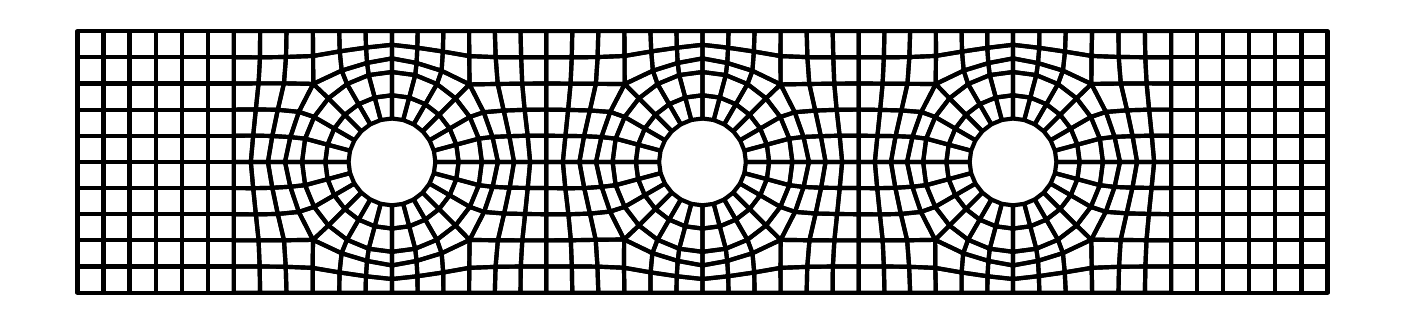}
		\caption{}
			\end{subfigure}	\hfill
	\begin{subfigure}[c]{0.49\linewidth } 
		\centering
		\includegraphics[width=\linewidth]{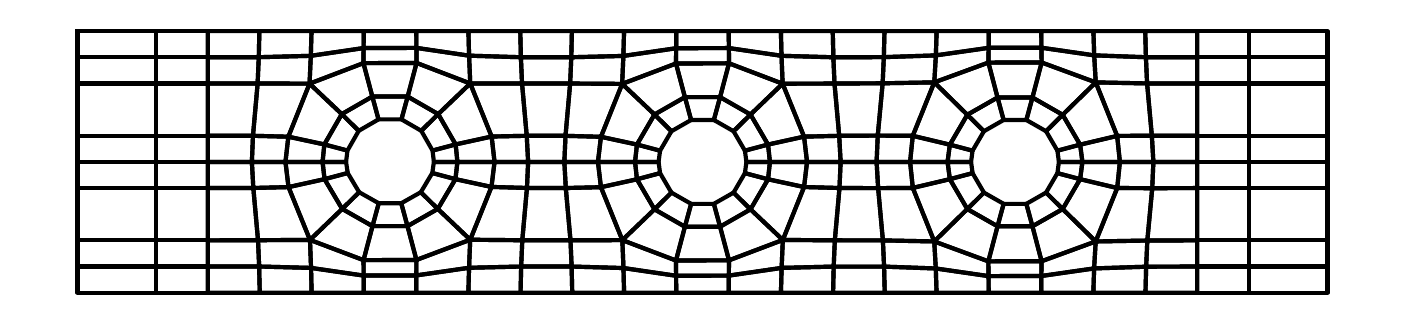}
		\caption{}
			\end{subfigure}		\caption{Interface meshes: (a) original 588 ZTE from Abaqus and (b) reduced to 232 ZTE based on \cite{balajiReducedOrderModeling2021}.} \label{fig:hbrb_interface_mesh}
\end{figure}

Utilizing Abaqus, an initial Hurty/Craig-Bampton (HCB) model reduction \cite{hurtyDynamicAnalysisStructural1965, craigCouplingSubstructuresDynamic1968} is conducted to the interface degrees of freedom plus 500 fixed interface modes.
Then, the model reduction procedure (including a second HCB) of \cite{balajiReducedOrderModeling2021} is applied to further reduce the model to 232 Zero-Thickness Elements (ZTEs), 282 nodes$ \times 3 $ DOFs, and 20 fixed interface modes. 
A null space transformation matrix described in \cite{balajiTractionbasedMultiscaleNonlinear2020} is also generated in this procedure to constrain out motion in the six rigid body modes of the assembled structure resulting in a final model with 860 DOFs. A seventh rigid body mode of the interface separating along the bolts is preserved \cite{balajiTractionbasedMultiscaleNonlinear2020}.
The resulting reduced interface mesh is the same as was used in \cite{porterPredictivePhysicsbasedFriction2023, balajiReducedOrderModeling2021}, but the mass and stiffness matrices correspond to the HBRB instead of the BRB. 
The original and reduced interface meshes are shown in \Cref{fig:hbrb_interface_mesh}.
For the reduced mesh, mapping matrices $\boldsymbol{Q}$ and $\boldsymbol{T}$ relate global displacements to displacements at one quadrature point per element and tractions at quadrature points to global forces respectively (see \eqref{eq:eom}).
Additionally, vectors are calculated for applying bolt tension and shaker excitation forces and for extracting accelerometer displacements.

The model also incorporates a linear damping matrix.
Linear viscous damping is not well understood \cite{akayResearchNeedsOpen2016}, so the present effort does not attempt to predict the linear damping component of the response.
Based on experimental shaker ring down tests, low amplitude damping values of the first
and third
bending modes of the assembled HBRB were identified as 1.35e-3 and 5.45e-4 respectively. These damping ratios were utilized to create a mass and stiffness proportional damping matrix around the prestressed state in the same way as~\cite{porterPredictivePhysicsbasedFriction2023}.

One small, but significant change to the model reduction code from \cite{balajiReducedOrderModeling2021} was necessary to achieve well conditioned and positive definite reduced mass and stiffness matrices. During the HCB step in the model reduction procedure, the unreduced stiffness matrix is represented as 
\begin{equation}
	\boldsymbol{K_{full}} = \begin{bmatrix}
		\boldsymbol{K_{bb}} & \boldsymbol{K_{bi}} \\
		\boldsymbol{K_{bi}}^T & \boldsymbol{K_{ii}}
	\end{bmatrix}
\end{equation}
and the solution of
\begin{equation}
	\boldsymbol{\Phi_C} = \boldsymbol{K_{ii}}^{-1} \boldsymbol{K_{bi}}^T
\end{equation}
is required. Since the HBRB model contains 6 rigid body modes with the interface degrees of freedom, this is poorly conditioned. While the previous work simply utilized the backslash operator in MATLAB to calculate a pseudo-inverse \cite{balajiReducedOrderModeling2021}, it resulted in reduced matrices that were not symmetric positive definite in this work. 
Instead, a null space transform is utilized for this matrix solution. Here, the singular value decomposition of $ \boldsymbol{K_{ii}} $ is calculated as
\begin{equation}
	\boldsymbol{K_{ii}} = \boldsymbol{U S V}^T.
\end{equation}
A null space matrix $ \boldsymbol{L_{ii}} $ is defined to be all but the last 6 columns of the matrix $ \boldsymbol{V} $. Finally, the matrix problem is solved as
\begin{equation}
	\boldsymbol{\Phi_C} = \boldsymbol{L_{ii}}
	\left[\boldsymbol{L_{ii}}^T  \boldsymbol{K_{ii}} \boldsymbol{L_{ii}} \right]^{-1} \left[\boldsymbol{L_{ii}}^T \boldsymbol{K_{bi}}^T \right].
\end{equation}
Next, the nonlinear friction forces applied in the reduced HBRB model are described.

\subsubsection{Nonlinear Frictional Forces}

The HBRB model utilizes a physics-based rough contact friction model developed in \cite{porterPredictivePhysicsbasedFriction2023}.
The friction model calculates nonlinear forces by integrating over a probability distribution of interactions between asperities (defined as locally maximal surface features); forces between 100 asperities with different initial gaps are utilized to numerically approximate this integral. 
Normal interactions between asperities are captured by elastic loading \cite{hertzUberdieBer1882}, elastic-plastic loading \cite{ghaednia_strain_2019}, and modified elastic unloading \cite{brakeAnalyticalElasticPlastic2015}.
Asperities are modeled with stick-slip tangential elements using the Cattaneo-Mindlin stiffness for contacting spheres \cite{cattaneoSulContatto1938, mindlinComplianceElasticBodies1949} and a Coulomb friction coefficient $\mu$ (proportional to normal force).
Capturing the mesoscale topology of the interface was critical in \cite{porterPredictivePhysicsbasedFriction2023}, so the present work deterministically models the mesoscale topology by offsetting the initial displacement for each element appropriately. 
For the present work, the friction model from \cite{porterPredictivePhysicsbasedFriction2023} is implemented in Python utilizing JAX \cite{JAXHighPerformanceArray} for automatic differentiation and just-in-time compilation \cite{porterTMDSimPy}. The JAX nonlinear force evaluations are effectively executed in parallel over available threads when running simulations. 

To obtain properties for the friction model, both interfaces of the HBRB are scanned. For the mesoscale topology, a KEYENCE LJ-V7080 in-line profilometer moved by a servomotor (at approximately 0.15 mm/s) is used to measure the full interface in one scan as was done in \cite{balajiTractionbasedMultiscaleNonlinear2020}. The scan data does not report coordinates with height measurements, so the edges of the interface were identified from scans and the x-y positions were scaled to fit the nominal interface size.
Once trimmed around the edges and holes, and leveled, the data was interpolated to a regular grid and a 2D Gaussian filter was applied with standard deviation of $\sigma=102$ points (5.1 mm) and a support of $2 \sigma +1$ points consistent with \cite{porterPredictivePhysicsbasedFriction2023}. 
While \cite{porterPredictivePhysicsbasedFriction2023} padded the scan data with border measurements, this work pads the scan data with the mirror reflection as this appeared to introduce less artificial effects when manually comparing to raw scan data.
The mesoscale topology is the total gap considering the filtered height data of both sides of the interface.

Both sides of the interface are also scanned in 12 patches with a KEYENCE VR-5100 White Light Interferometer at a vertical resolution of 4 $\mu$m and a horizontal resolution of 23.6 $\mu$m. Each patch of data is trimmed, leveled, and a mesoscale topology is calculated. The mesoscale topology is subtracted from the leveled height data and then the asperity processing procedure from \cite{porterPredictivePhysicsbasedFriction2023} is applied within each patch.\footnote{An erosion size of 2 points is utilized for identifying asperities consistent with the baseline identification used in \cite{porterPredictivePhysicsbasedFriction2023}.} 
The asperity processing uses a watershed algorithm based on \cite{wen_reconstruction_2020} to identify potential regions of asperities and a total ellipsoid fit to obtain asperity properties. 
The median of the asperity radii across all 24 scans (12 patches times 2 sides of the interface) is taken as a friction model input, and the heights of all asperities are utilized to calculate a probability distribution for the initial gap between two asperities.

Overall, the asperity identification procedure is nearly identical to \cite{porterPredictivePhysicsbasedFriction2023} with slight adjustments since the scans from the KEYENCE VR-5100 are not stitched together.
On the other hand, this work utilizes the KEYENCE LJ-V7080 to obtain the mesoscale topology instead of attempting to stitch scans from multiple patches as was done in \cite{balajiDissipative2021} and used in \cite{porterPredictivePhysicsbasedFriction2023}. 
Additionally, the process for calculating the mesoscale topology required manual updates (e.g., how much area should be trimmed around the edges and holes) and comparisons to the raw data to have confidence in the results.
However, this process is still easier than stitching together multiple scans from the KEYENCE VR-5100.
The calculated mesoscale topology of each patch from the KEYENCE VR-5100 is only used to get consistent heights of asperities relative to the smoothed surface and not in the model because the exact locations of each scan are not calculated. This does not limit the use of asperity statistics from the KEYENCE VR-5100, because those statistics are for the full surface not specific locations.

\section{Half Brake-Reu{\ss} Beam Results} \label{sec:hbrb_results}

\subsection{Surface Processing Results}

\begin{figure}[!h]
	\centering
	\begin{subfigure}[c]{0.59\linewidth } 
		\centering
						\includegraphics[width=\linewidth, trim={2.5cm 0 3.2cm 0},clip]{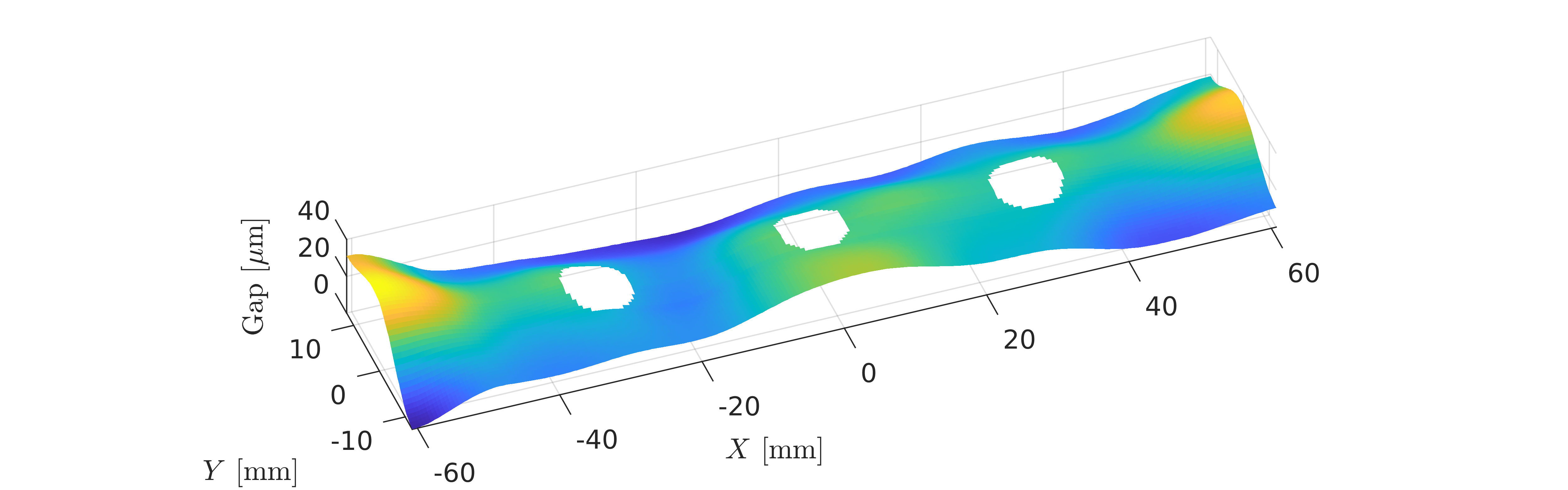}
		\caption{}
			\end{subfigure}	\hfill
	\begin{subfigure}[c]{0.39\linewidth } 
		\centering
				\includegraphics[width=\linewidth]{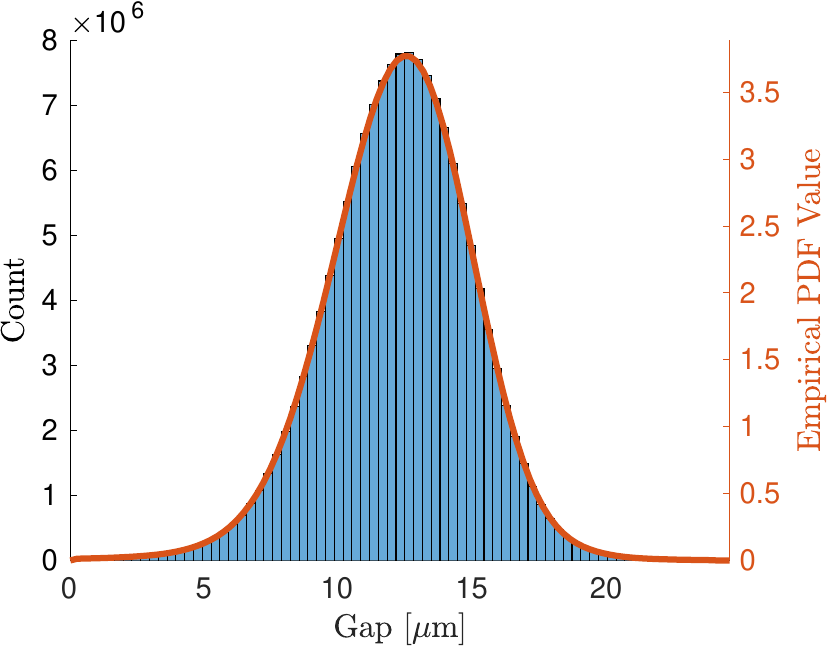}
		\caption{}
			\end{subfigure}		\caption{HBRB surface processing results: (a) mesoscale topology and (b) distribution of asperity gaps.} \label{fig:hbrb_surface}
\end{figure}

\Cref{fig:hbrb_surface} shows the mesoscale topology and the distribution of initial gaps between asperities. 
Here, the mesoscale topology has a maximum gap of 38.7 $ \mu $m, with many of the largest gaps appearing towards the edges.
\Cref{tab:hbrb_asp_results} presents the asperity statistics for the HBRB compared to those previously identified for the BRB \cite{porterPredictivePhysicsbasedFriction2023}. The asperity statistics are very similar to those reported for the BRB in \cite{porterPredictivePhysicsbasedFriction2023} suggesting that asperity statistics calculated based on the machining process (since the beams were manufactured in the same way) would be appropriate to use prior to initial prototyping. 
Next, dynamic results are presented for the HBRB.

\FloatBarrier

\begin{table}[h!]
	\caption{Asperity parameters for rough contact modeling including BRB results from \cite{porterPredictivePhysicsbasedFriction2023} for reference and newly calculated HBRB surface parameters.} \label{tab:hbrb_asp_results}
	\centering
	\begin{tabular}{cccc}
		\hline \hline
		Parameter & Description & BRB Value \cite{porterPredictivePhysicsbasedFriction2023} & HBRB Value \\ \hline
		$ R' $  & Large principal relative radius of curvature & 3.167 mm & 2.849 mm\\
		$ R'' $ & Small principal relative radius of curvature & 0.588 mm & 0.504 mm\\
		$ R_e $ & Equivalent radius & 1.365 mm & 1.199 mm \\
		$ z_{max} $ & Maximum gap between asperities &  28.05 $ \mu $m & 24.50 $ \mu $m \\
				$ \eta $ & Area density of asperities & 1.371 mm$ ^{-2} $ & 1.544 mm$ ^{-2}$ \\
		\hline \hline
	\end{tabular}
\end{table}

\FloatBarrier

\subsection{Single Mode Results}

\subsubsection{Experimental Results}

\begin{figure}[h!]
	\centering
	\includegraphics[width=0.55\linewidth]{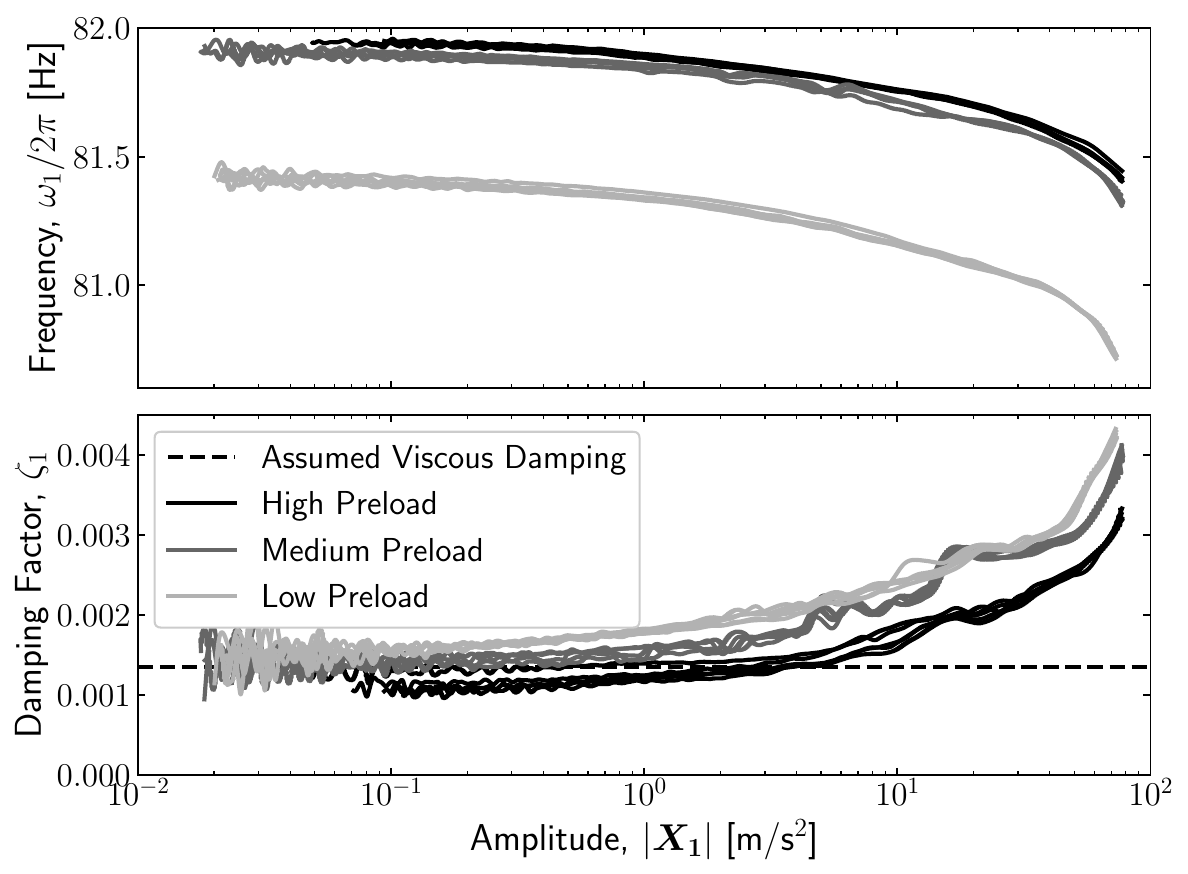}
	\caption{
	Experimental shaker ring down results for the first bending mode of the HBRB comparing the three preload values of 12006 (high), 9069 (medium), and 6135 (low) N per bolt. 
	The assumed viscous damping that is used in the models is plotted on the bottom for reference. The acceleration is at the tip accelerometer (see \Cref{fig:hbrb_ringdown_exp}). Six response backbones are plotted for each preload level. 
}
	\label{fig:hbrb_ringdown_exp}
\end{figure}

The shaker ring down tests for the first bending mode show nonlinear trends in frequency and damping in \Cref{fig:hbrb_ringdown_exp}. 
Time series data from shaker ring down is processed with the Peak Finding and Fitting (PFF) method \cite{jinIdentificationInstantaneousFrequency2020} to calculate the instantaneous frequency and damping as a function of amplitude.\footnote{For the first bending mode, a third order Butterworth bandpass filter is used with a center frequency of 82 Hz and a bandwidth of 10 Hz. 
End effects at high and low amplitude levels from PFF are removed in \Cref{fig:hbrb_ringdown_exp}.}
As the preload is decreased, the system exhibits a lower natural frequency, a larger shift in frequency, and increased damping. For all cases, the nonlinearity is most strongly observed in the damping with an increase by a factor of two to three.

\FloatBarrier

\FloatBarrier

Shaker ring down and PFF are also applied to the third bending mode (corresponding to the superharmonic resonance response in \Cref{sec:hbrb_sr_results}).\footnote{For the third bending mode, the third order Butterworth bandpass filter used a center frequency of 569 Hz with a bandwidth of 10 Hz.}
No visible nonlinear trends were observed utilizing PFF, so properties are calculated based on all amplitudes together and reported in \Cref{tab:hbrb_ringdown_m3}.
Minimal dependency on preload is observed for the third bending mode.

\begin{table}[h!]
	\centering
	\caption{Shaker ring down results for the third bending mode of the HBRB. No clear nonlinear trends were observed, so the properties are assumed constant over all amplitude levels. Amplitude levels are reported for the tip accelerometer.} 
	\label{tab:hbrb_ringdown_m3}
	\begin{tabular}{cccc}
		\hline \hline 
		& \multicolumn{3}{c}{Preload} \\ \cline{2-4} 
		& Low (6135 N) & Medium (9069 N) & High (12006 N) \\ \hline
		Lowest Amplitude, Tip [m/s$ ^2 $]  	& 1.93e-2 & 2.00e-2 & 1.53e-2 \\
		Highest Amplitude, Tip [m/s$ ^2 $] 	& 5.39e1  & 5.10e1  & 5.49e1 \\ \hline
		Mean Frequency, $ \omega_3 $ [Hz]               & 568.3 & 569.6 & 569.6 \\
		Standard Deviation Frequency, $ \sigma $ [Hz] & 0.3   & 0.3   & 0.3\\ \hline
		Mean Damping Factor, $ \zeta_3 $ & 5.85e-4 & 5.06e-4 & 5.43e-4  \\
		Standard Deviation Damping Factor, $ \sigma $ & 8.85e-5 & 8.02e-5 & 7.80e-5 \\
		\hline \hline
	\end{tabular}
\end{table}

\FloatBarrier

\subsubsection{Modeling Results} \label{sec:hbrb_epmc_results}

\begin{figure}[h!]
	\centering
	\includegraphics[width=0.55\linewidth]{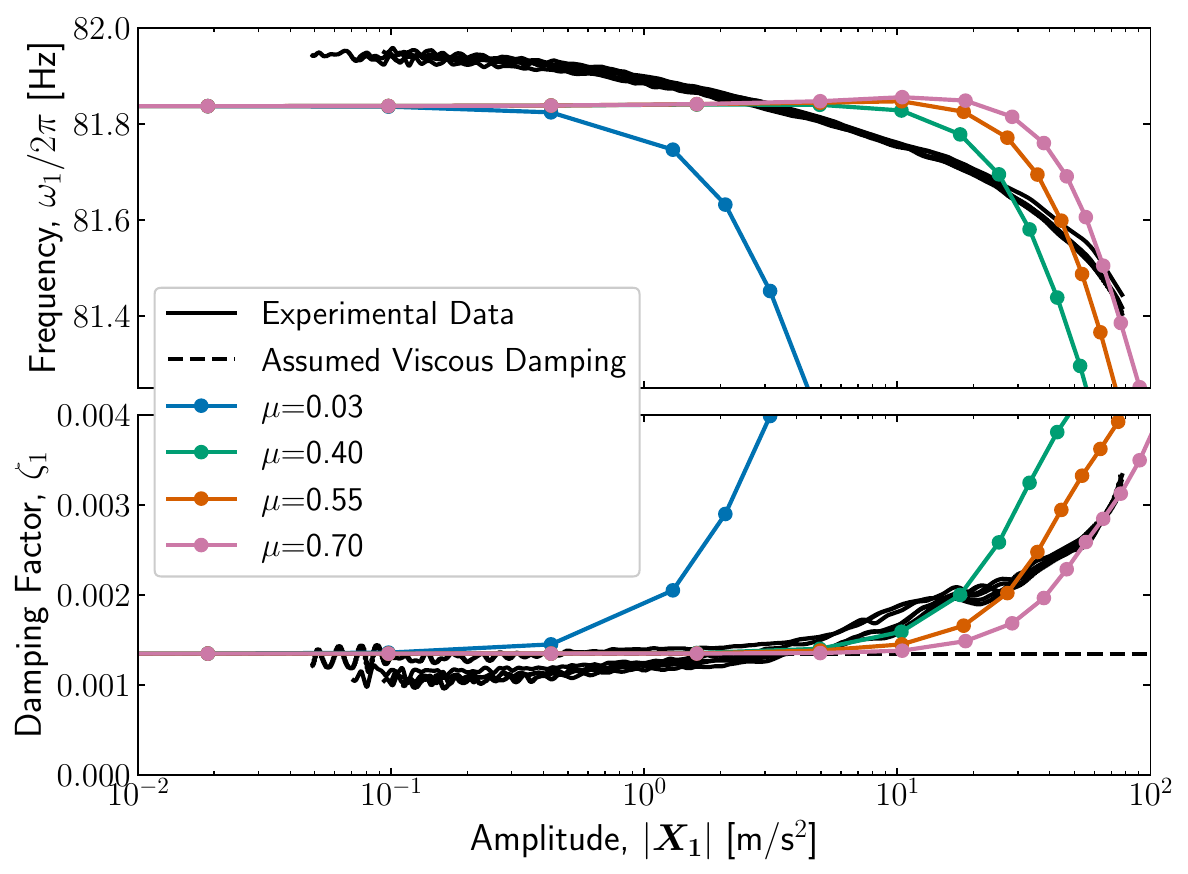}
	\caption{First bending mode EPMC results for different friction coefficients $ \mu $ compared to experiments using the highest preload value and the tip accelerometer. A friction coefficient of $ \mu=0.55 $ is chosen for subsequent models.}
	\label{fig:hbrb_epmc_high}
\end{figure}

The HBRB model (see \Cref{sec:model_procedure}) is applied to simulate the first bending mode of the HBRB at the highest preload (see  \Cref{fig:hbrb_epmc_high}) with EPMC \cite{krackNonlinearModalAnalysis2015} and the physics-based friction model \cite{porterPredictivePhysicsbasedFriction2023}. 
EPMC is utilized with harmonics 0, 1, 2, and 3 and 128 time steps for AFT, which is consistent with \cite{porterPredictivePhysicsbasedFriction2023}.
Here, the linear frequency is a prediction from this model by directly applying the modeling procedure from \cite{porterPredictivePhysicsbasedFriction2023} and achieves a close match to the experiments with an error of less than 0.2 Hz. 
The only manually tuned parameters were for processing the mesoscale topology (prior to looking at frequency results) and the friction coefficient based on the EPMC results (discussed next).
The linear damping is assumed based on the experiments, so it is not a feature of the nonlinear friction model (see \Cref{sec:model_procedure}).

The friction coefficient for the model cannot be predicted prior to comparing the model and experiments. 
The friction coefficient of $\mu=0.03$ from \cite{porterPredictivePhysicsbasedFriction2023} shows frequency and damping shifts at much lower frequencies than the experimental results in \Cref{fig:hbrb_epmc_high}. 
Based on running simulations with several friction coefficients, a value of $\mu=0.55$ is selected for the remainder of this work.
The higher friction coefficient for this work compared to \cite{porterPredictivePhysicsbasedFriction2023} could be attributed to normal load dependency of the friction coefficient (since the mesoscale topology and normal pressures differ), a different operating regime with less nonlinearity, modal interactions (since experiments in \cite{porterPredictivePhysicsbasedFriction2023} used hammer impact tests), or other unknown differences in physics.

Here, the qualitative trends of the models show a more abrupt increase in damping and decrease in frequency than the experimental results. 
This difference in qualitative trends between the models and experiments was not observed in \cite{porterPredictivePhysicsbasedFriction2023}, potentially because the experiments in that work showed larger frequency and damping shifts. 
Therefore, it is possible that the physics-based friction model is ill-suited to predicting the regime of the smaller frequency and damping shifts observed in this work. 
Results could be improved by the Mindlin-Iwan Fit model from \cite{porterPredictivePhysicsbasedFriction2023}, which further discretizes each asperity as a set of sliders based on the Mindlin partial slip solution \cite{mindlinComplianceElasticBodies1949, mindlinEffectsOscillatingTangential1952, mindlinElasticSpheresContact1953}, and results in greater damping and frequency shifts at low amplitudes \cite{porterPredictivePhysicsbasedFriction2023}.
This discretization of the asperities increased computation time by a factor of 6 in \cite{porterPredictivePhysicsbasedFriction2023}, so is not considered here.

Also of interest, \Cref{fig:hbrb_epmc_high_physics} shows the effect of simplifying assumptions in the friction model on the EPMC backbone at high preload. Utilizing the BRB asperity statistics from \cite{porterPredictivePhysicsbasedFriction2023} gives reasonable results, indicating that asperity statistics can be based on the manufacturing process because both the BRB and HBRB were manufactured in the same way.
Using elastic asperities or a flat mesoscale topology also decreased the natural frequency, which is consistent with \cite{porterPredictivePhysicsbasedFriction2023}. 
Lastly, excluding the mesoscale topology significantly changes how the system transitions from stick to slip with a smaller decrease in frequency, but a comparable increase in damping.
Additional results at low and medium preload levels for the first bending mode are included in \Cref{sec:extra_epmc_appendix}.

\newcommand{\EpmcPhysicsCaption}[1]{First bending mode EPMC results for different physics assumptions compared to experiments using the {#1} preload value, the tip accelerometer, and a friction coefficient of $ \mu=0.55 $.}

\FloatBarrier

\begin{figure}[h!]
	\centering
	\includegraphics[width=0.55\linewidth]{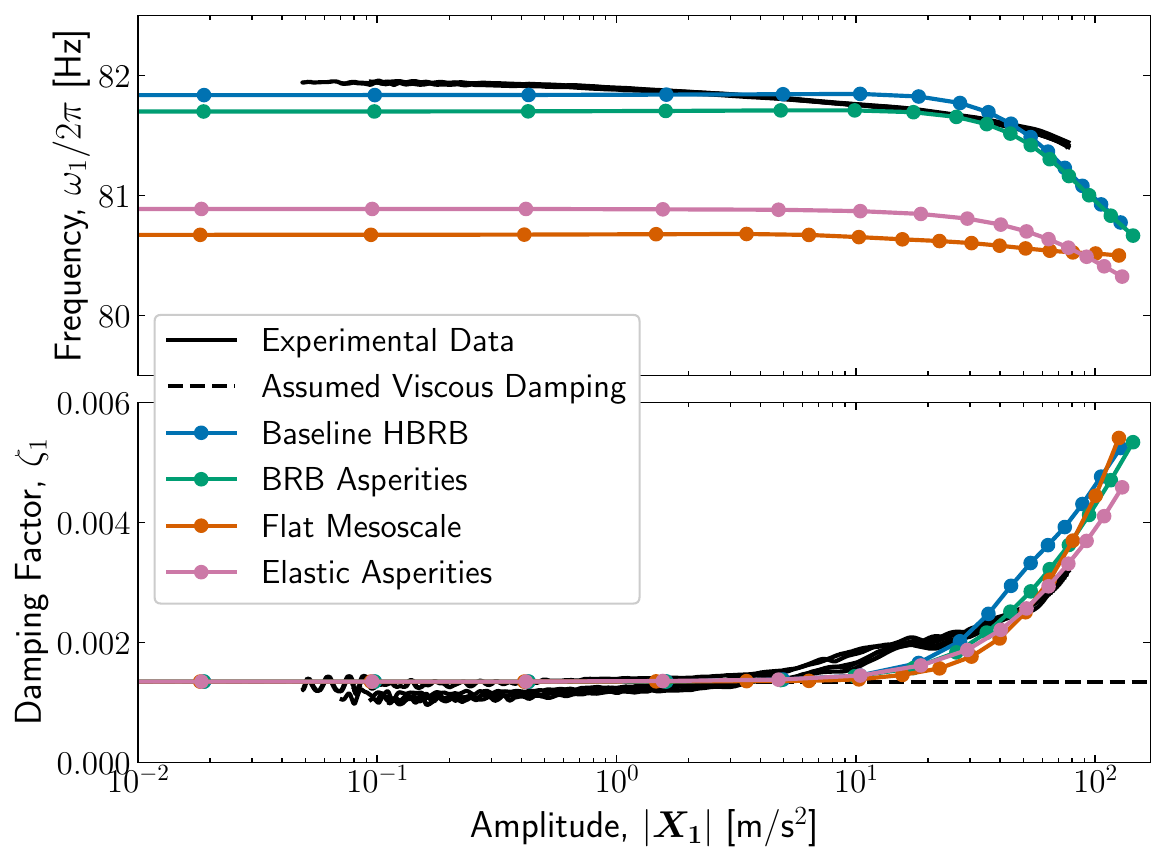}
	\caption{\EpmcPhysicsCaption{high}
	}
	\label{fig:hbrb_epmc_high_physics}
\end{figure}

\FloatBarrier

\begin{figure}[h!]
	\centering
	\includegraphics[width=0.55\linewidth]{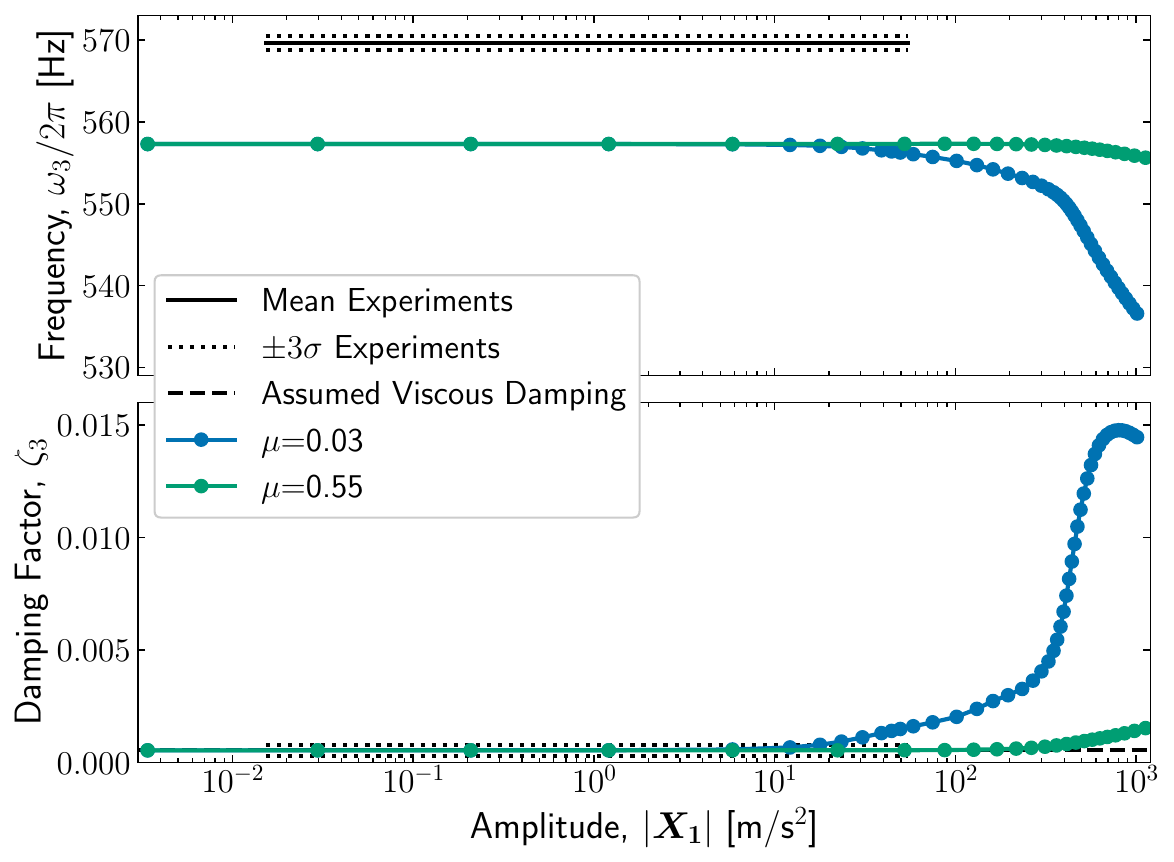}
	\caption{
	Third bending mode EPMC results for different friction coefficients $ \mu $ compared to experiments (with standard deviation $\sigma$) using the highest preload value and the tip accelerometer. A friction coefficient of $ \mu=0.55 $ is chosen for subsequent models.	
	Since the linear damping for the model is chosen as the mean experimental damping, the mean experimental damping falls directly below the plotted models. 
}
	\label{fig:hbrb_epmc_third_high}
\end{figure}

\Cref{fig:hbrb_epmc_third_high} shows the EPMC results for the third bending mode. Similar to the first bending mode, a friction coefficient of $\mu=0.03$ shows slip not observed in the experiments. The chosen friction coefficient of $\mu=0.55$ remains mostly linear over the regime tested in the experiments.
As with the first bending mode results, the model exactly matches the experimental damping because the model is set to have low amplitude viscous damping based on the experimental results, so this is not a prediction of the friction model.
On the other hand, the model linear frequency of 557 Hz is a prediction with a 2\% error compared to the experiments. 
EPMC results for the third bending mode at low and medium preload levels are included in \Cref{sec:extra_epmc_appendix}.

\FloatBarrier

Next, the superharmonic resonances of the HBRB are considered. 
Based on the EPMC results, exact matches between experiments and models for the superharmonic resonance are not expected since the nonlinear trends are not exactly captured. 
Furthermore, the lower frequency of the third bending mode in the model compared to the experiment is expected to result in a similar error in the forcing frequency giving rise to the superharmonic resonance.
However, the physics-based friction model still captures characteristics of the nonlinear trends making it relevant and useful for understanding the superharmonic resonances in the present work.

\subsection{Superharmonic Resonance Cases} \label{sec:hbrb_sr_results}

This section focuses on a 7:1 superharmonic resonance where the third bending mode of the HBRB responds at seven times the forcing frequency. 
This superharmonic resonance occurs near the primary resonance of the first bending mode, so the case is an internal resonance. 
Other superharmonic resonances are observed in the experimental data, but do not consistently appear because the experimental tests focused on the 7:1 superharmonic resonance of the third bending mode. 
Furthermore, other observed superharmonic resonances had smaller amplitudes in the collected data. 
Therefore, this section only considers the 7:1 superharmonic resonance of the third bending mode.

\subsubsection{Experimental Results} \label{sec:exp_stepped_sine}

\begin{figure}[h!]
	\centering
	\begin{subfigure}[c]{\linewidth } 
		\centering
		\includegraphics[width=\linewidth]{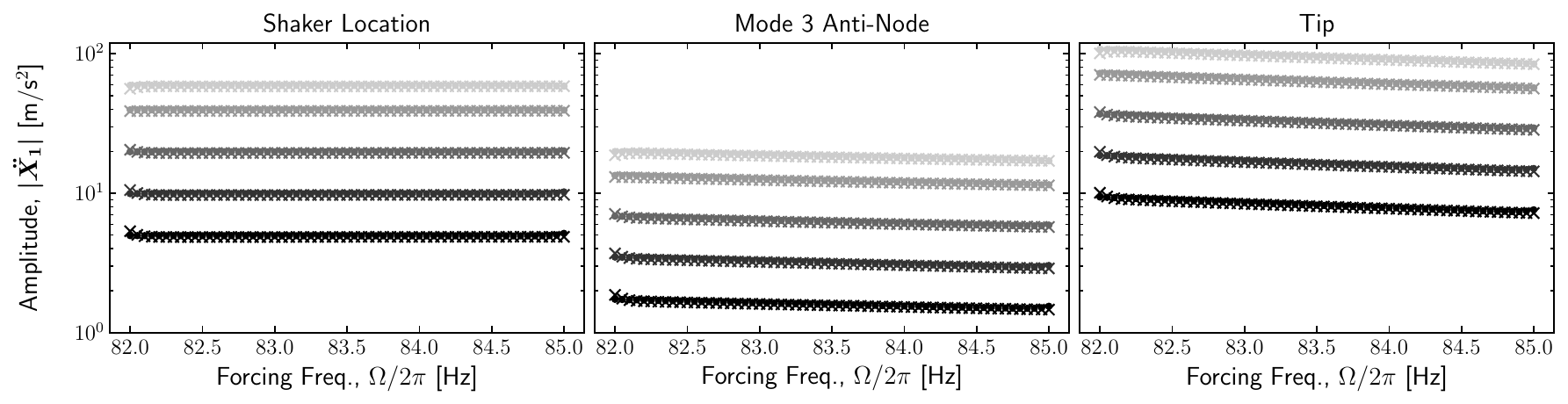}
		\caption{}
		\label{fig:exp_highPre_h1}
	\end{subfigure}
	\\
	\begin{subfigure}[c]{\linewidth } 
		\centering
		\includegraphics[width=\linewidth]{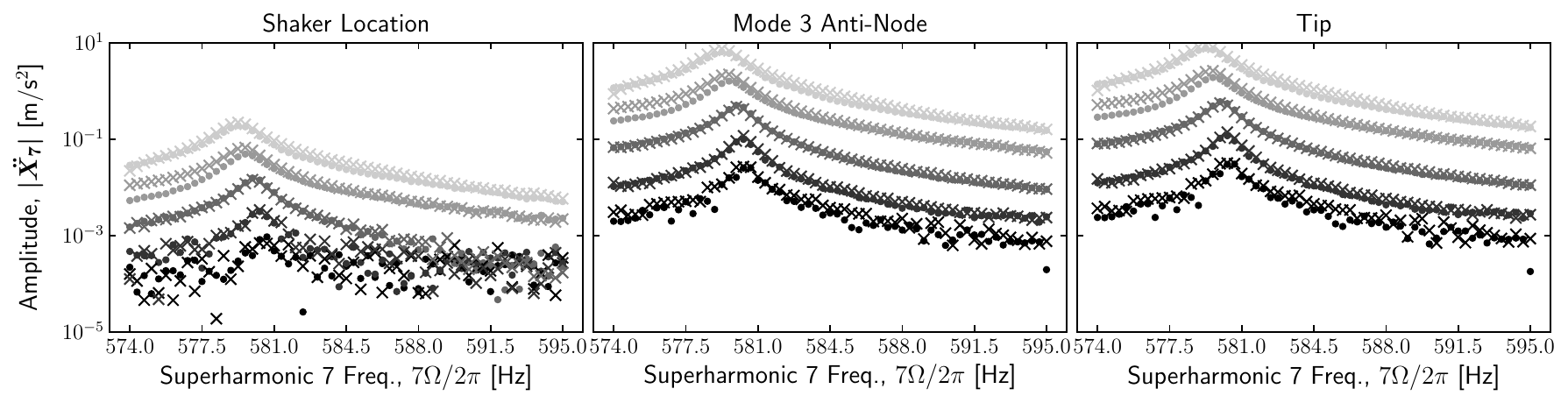}
		\caption{}
		\label{fig:exp_highPre_h7}
	\end{subfigure}
	\caption{
	HBRB experimental results for the high preload case stepped-sine response of (a) first harmonic and (b) seventh harmonic. 
	From dark to light, the tests are at amplitudes of 0.5 g, 1.0 g, 2.0 g, 4.0 g and 6.0 g for the first harmonic at the shaker.
	Tests are conducted sweeping from low to high frequency ($ \boldsymbol{\cdot} $) and high to low frequency ($ \boldsymbol{\times} $).
	The forcing frequency bounds are the same for (a) and (b), but (b) plots the superharmonic response frequency on the x-axis corresponding to seven times the forcing frequency.
}
	\label{fig:exp_highPre}
\end{figure}

Stepped-sine testing characterizes the superharmonic resonance as shown in \Cref{fig:exp_highPre} for the highest preload case. 
The first harmonic response is controlled at the shaker resulting in near constant response amplitudes for all three measurement locations (see \Cref{fig:exp_highPre_h1}).
The final tests at 1.0 g are not included, but showed consistent results with the initial 1.0 g results indicating repeatability after the 6.0 g tests.
The response of the seventh harmonic shows a clear peak corresponding to the superharmonic resonance. At the mode 3 anti-node, the superharmonic resonance achieves nearly half of the response amplitude of the first harmonic response for the 6.0 g (highest amplitude) case. 
If a sensitive component was located at the anti-node of mode 3 (and near the node of mode 1), it would be critical to capture acceleration due to the superharmonic resonance.
For both the first and seventh harmonics, the tests sweeping from low to high and high to low frequencies showed very repeatable results.
Experimental results for low and medium preload levels are included in \Cref{sec:extra_stepped_sine_appendix}.

\FloatBarrier

\begin{figure}[h!]
	\centering
	\includegraphics[width=\linewidth]{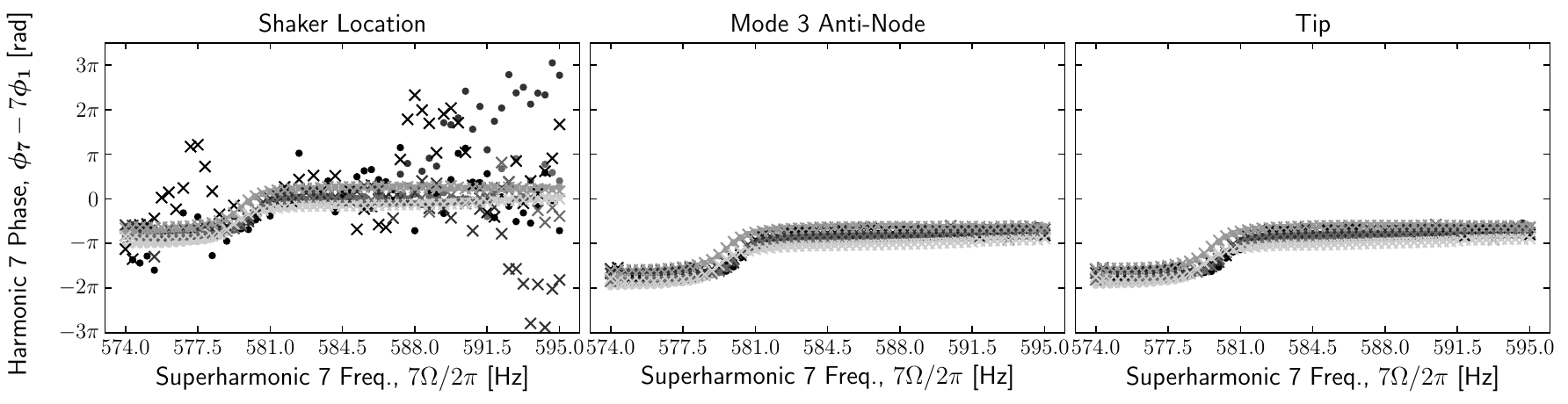}
	\\
	\includegraphics[width=\linewidth]{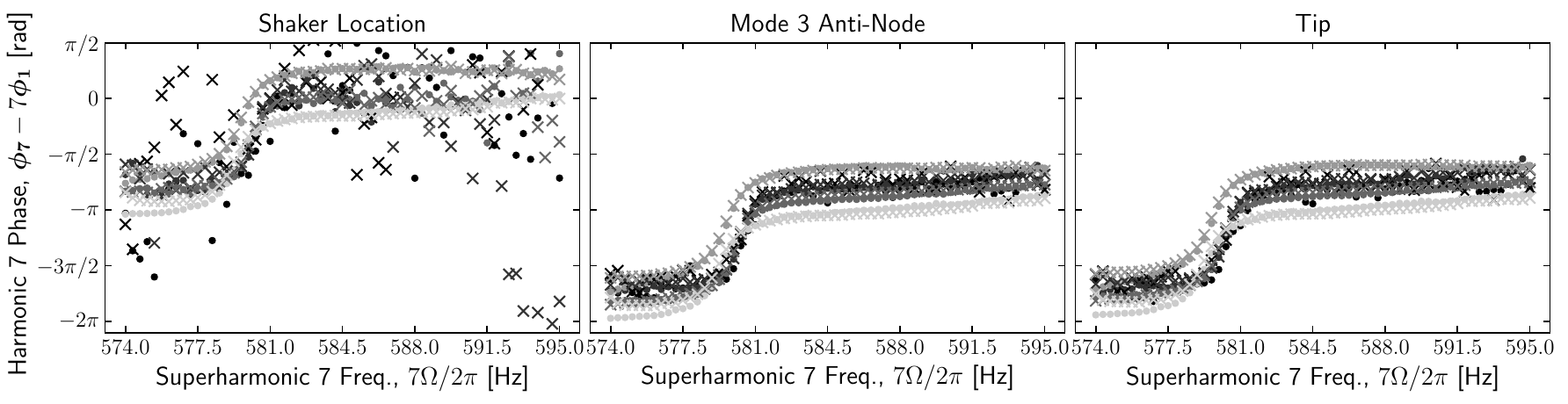}
	\caption{
	HBRB experimental phase of the seventh superharmonic resonance for the high preload case stepped-sine tests. 
	From dark to light, the tests are at amplitudes of 0.5 g, 1.0 g, 2.0 g, 4.0 g and 6.0 g for the first harmonic at the shaker.
	Tests are conducted sweeping from low to high frequency ($ \boldsymbol{\cdot} $) and high to low frequency ($ \boldsymbol{\times} $).
	Bottom is zoomed in version of top plot. There is significant noise for the phase at the shaker at low acceleration amplitudes (darker points) due to the small superharmonic displacements near the node of the superharmonic mode.
	}
	\label{fig:exp_highPre_h7Phase}
\end{figure}

The superharmonic resonance is also confirmed by the phase of the seventh harmonic in \Cref{fig:exp_highPre_h7Phase} that shows a clear increase by $\pi$ for all measurement locations.
At low amplitudes, the phase at the shaker exhibits significant noise since the shaker location is near a node for the third bending mode resulting in very small superharmonic amplitudes, especially at low control acceleration levels (see \Cref{fig:exp_highPre_h7}).

\FloatBarrier

\begin{figure}[h!]
	\centering
	\begin{subfigure}[b]{0.49\linewidth } 
		\centering
		\includegraphics[width=\linewidth]{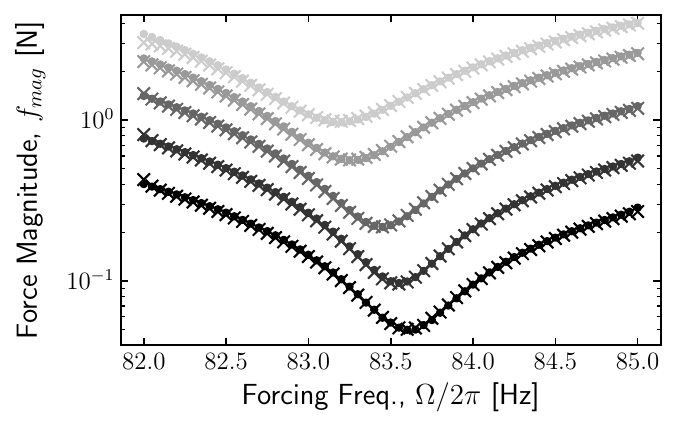}
		\caption{}
		\label{fig:exp_highPre_force_h1}
	\end{subfigure}
	\hfill
	\begin{subfigure}[b]{0.49\linewidth } 
		\centering
		\includegraphics[width=\linewidth]{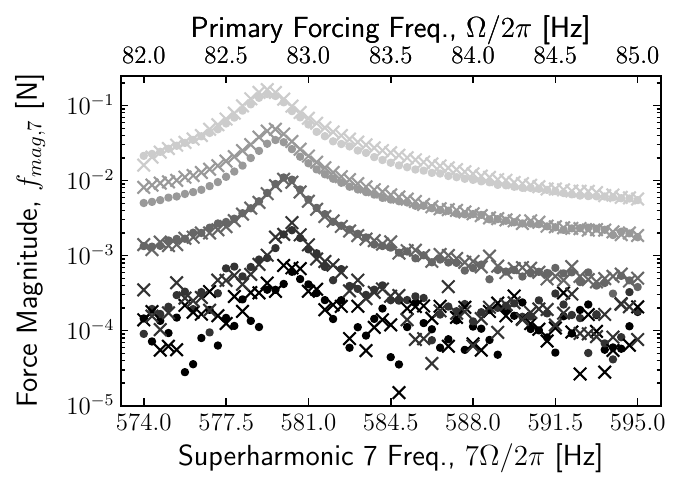}
		\caption{}
		\label{fig:exp_highPre_force_h7}
	\end{subfigure}
	\caption{
HBRB experimental force measurements of (a) first harmonic and (b) seventh harmonic for the high preload case stepped-sine tests.
From dark to light, the tests are at amplitudes of 0.5 g, 1.0 g, 2.0 g, 4.0 g and 6.0 g for the first harmonic at the shaker.
Tests are conducted sweeping from low to high frequency ($ \boldsymbol{\cdot} $) and high to low frequency ($ \boldsymbol{\times} $).
The forcing frequency bounds are the same for (a) and (b), but (b) plots the superharmonic response frequency on the x-axis corresponding to seven times the forcing frequency.
The first harmonic force in (a) is the desired input to the structure while the seventh harmonic force in (b) is undesired structure-shaker interaction.
	}
	\label{fig:exp_highPre_force}
\end{figure}

To achieve a constant amplitude, the controller adjusts the forcing magnitude shown in \Cref{fig:exp_highPre_force_h1}. Additionally, there are some shaker interaction effects resulting in non-zero force acting on the structure at the seventh harmonic corresponding to the superharmonic resonance as shown in \Cref{fig:exp_highPre_force_h7}.
The force on the superharmonic is approximately an order of magnitude smaller than that on the primary harmonic. As shown in the left most plot of \Cref{fig:exp_highPre_h7}, the response amplitude of the superharmonic at the shaker is much smaller than that of the first harmonic since it is near a node of the third bending mode of the HBRB. Therefore, the structure shaker interaction is expected to contribute minimal modal forcing to the superharmonic resonance, and it is believed that the superharmonic resonance is still primarily caused by the nonlinearity in the joint.
Advanced control algorithms have been demonstrated for eliminating higher harmonics of the shaker excitation \cite{hippoldShakerStructureInteraction2023}, but hardware to utilize such algorithms is unavailable for this study.

\FloatBarrier

\subsubsection{Modeling Results} \label{sec:hbrb_hbm_res}

HBM as formulated in \Cref{sec:theory} with the amplitude and phase constraint from \Cref{sec:theory_high_dim} is used to simulate the HBRB superharmonic resonance behavior shown in \Cref{fig:rom_highPre} and as a truth solution for comparing the VPRNM ROM.
HBM utilizes the zeroth and first nine harmonics. For acceleration levels of 0.5 g, 1.0 g and 2.0 g, results using harmonics 0 and 1-7 or 0 and 1-9 are consistent. 
For acceleration levels of 4.0 and 6.0 g, including harmonics 0 and 1-9 is significantly different from using only 0 and 1-7, but consistent with including the 0th and first 11, 13, or 15 harmonics. 
HBM convergence results are presented in  \Cref{sec:hbm_convg}.
The VPRNM ROM applies displacement amplitude control while the experiments and HBM simulations use acceleration amplitude control. The displacement for the VPRNM ROM is calculated by dividing the acceleration by the square of expected frequency for the superharmonic resonance (one seventh of the low amplitude EPMC frequency for the third bending mode, 79.6 Hz = 500 rad/s). Given the relatively narrow frequency band of interest, the difference between displacement and acceleration control is minor.
Overall, the controlled first harmonic responses are consistent between the experiments, HBM, and the VPRNM ROM in \Cref{fig:rom_highPre_h1}.

\FloatBarrier

\newcommand{\HbrbAmpCaption}[1]{HBRB experimental, HBM, and VPRNM ROM results for (a) first harmonic and (b) seventh harmonic at {#1} preload. 
From dark to light, the HBM solutions and experiments are at amplitudes of 0.5~g, 1.0~g, 2.0 g, 4.0 g and 6.0 g for the first harmonic at the shaker.
The forcing frequency bounds are the same for (a) and (b), but (b) plots the superharmonic response frequency on the x-axis corresponding to seven times the forcing frequency. The EPMC and VPRNM backbones track the primary and superharmonic resonance frequencies respectively. Neither resonances significantly changes the controlled first harmonic amplitude in (a).}

\begin{figure}[h!]
	\centering
	\includegraphics[width=0.85\linewidth]{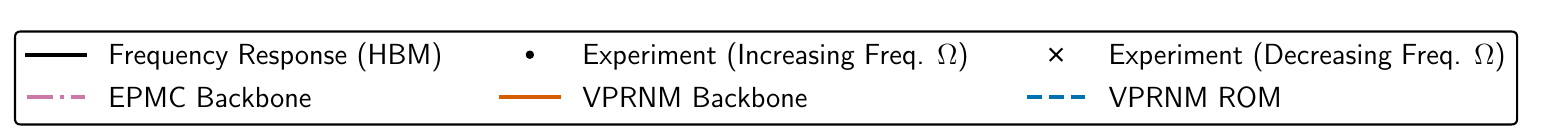}
	\begin{subfigure}[c]{\linewidth } 
		\centering
		\includegraphics[width=\linewidth]{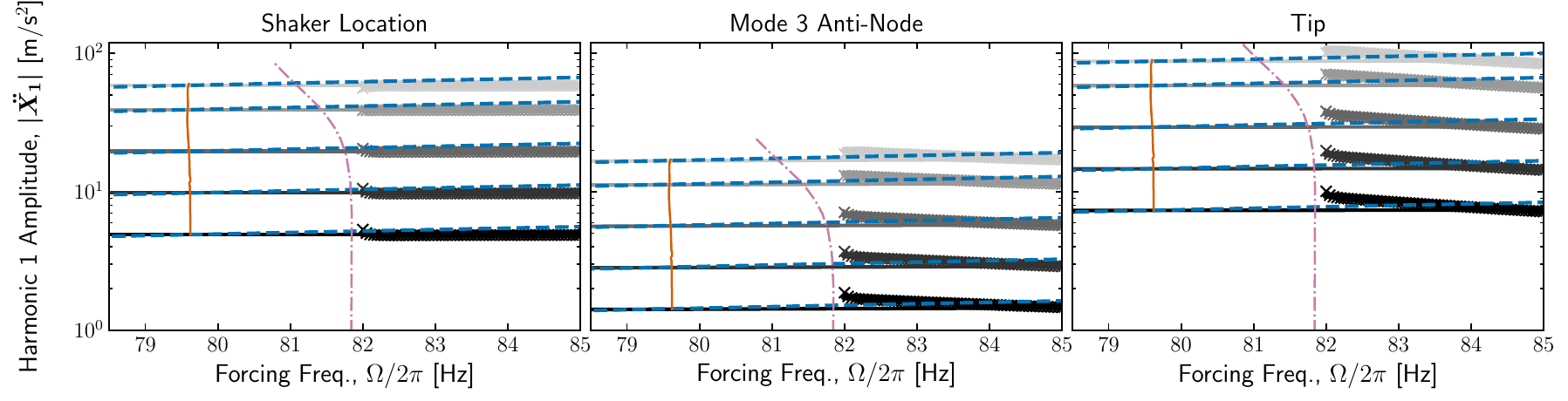}
		\caption{}
		\label{fig:rom_highPre_h1}
	\end{subfigure}
	\\
	\begin{subfigure}[c]{\linewidth } 
		\centering
		\includegraphics[width=\linewidth]{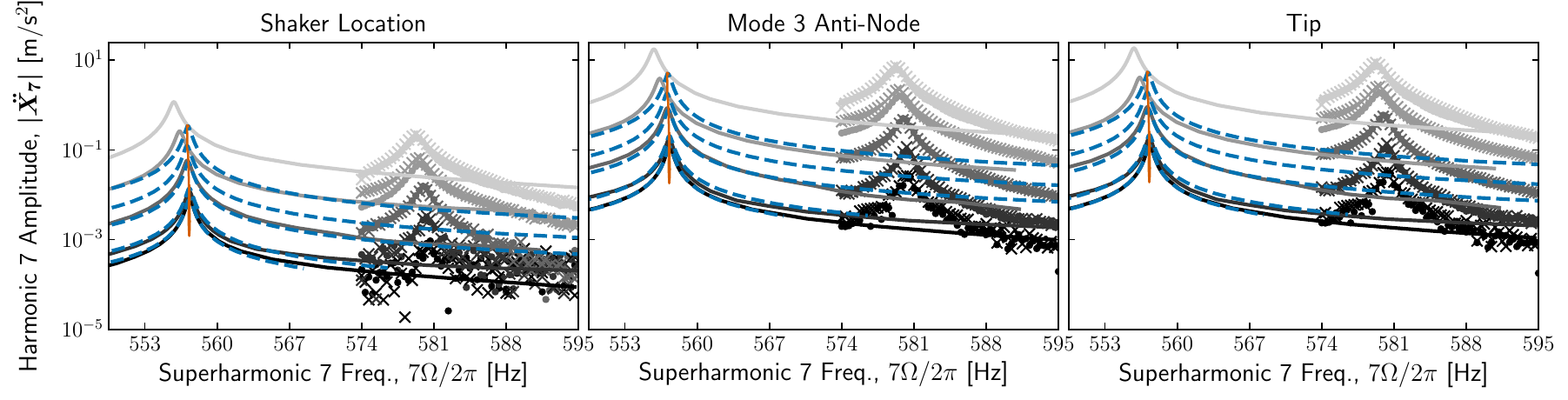}
		\caption{}
		\label{fig:rom_highPre_h7}
	\end{subfigure}
	\caption{\HbrbAmpCaption{high} }
	\label{fig:rom_highPre}
\end{figure}

The superharmonic resonance modeling and experimental results are qualitatively similar in \Cref{fig:rom_highPre_h7}. The differences in frequencies for the superharmonic resonances are attributed to errors in modeling the third bending mode of the HBRB and an increase in the frequency of the third bending mode when the shaker is attached (see \Cref{fig:hbrb_epmc_third_high}). 
Differences in response magnitude between models and experiments are attributed to not exactly capturing the nonlinear behavior of the HBRB (see \Cref{fig:hbrb_epmc_high}).
For the lowest three amplitudes, the VPRNM ROM precisely captures the superharmonic resonance behavior illustrating the utility of the approach.\footnote{The VPRNM ROM does not provide superharmonic solutions at all frequencies for the 0.5 g and 1.0 g cases because the response amplitude drops below the minimum amplitude captured in the EPMC simulations of the third bending mode. The EPMC backbone for the third bending mode is also up-sampled with linear interpolation to add an additional 100 amplitude levels to improve the resolution of the VPRNM ROM output.} 
For the 4.0 g and 6.0 g amplitudes, there are some clear errors in the VPRNM ROM because the VPRNM backbone does not track the peak response. 
These errors are likely caused by interactions between harmonics that are not accounted for in the VPRNM phase constraint (e.g., $\boldsymbol{F_{kq,sup,n}}$ from \eqref{eq:Fsup_def}). This is supported by the need for additional harmonics in HBM at the highest two amplitude levels (see \Cref{sec:hbm_convg}).
For additional investigations of how VPRNM breaks down for both unilateral spring contact (similar to normal direction) and hysteretic frictional forces (similar to tangent direction), see \cite{porterTrackingSuperharmonic2024}.
Analysis of the excitation mechanisms of the superharmonic resonance in \Cref{sec:hbrb_sr_causes} provides further insights into how VPRNM breaks down.
Despite some errors in the VPRNM ROM, it still provides useful characterizations of the superharmonic resonance that other ROMs (e.g., single nonlinear mode theory with EPMC) cannot capture.

To achieve constant amplitude, the forcing magnitude is varied in \Cref{fig:rom_highPre_force}. Again, the models and experiments show qualitative agreement.
The modeling errors in frequency and damping (i.e., minimum force levels) are attributed to attaching the shaker modifying the response of the beam, since the model is consistent with the low amplitude behavior from shaker ring down (see \Cref{fig:hbrb_epmc_high}).
Here, the HBM results show significant faceting due to the combination of adaptive continuation step sizing and a mostly constant solution with the amplitude and phase constraint of the first harmonic (see \Cref{sec:theory_high_dim}).
Given that the primary focus of this work is on the superharmonic resonance, which is well resolved by the adaptive step sizing, the lower computational costs for HBM to not resolve the primary resonance are accepted.
The VPRNM ROM matches the calculated HBM solutions well and provides higher resolution for comparison to experiments.
The VPRNM backbone tracks the frequencies of the superharmonic resonances and no notable disturbance of the external force is observed. Therefore, the force correction of \eqref{eq:vprnm_rom_force} is not applied here.
Alternatively, the EPMC backbone in \Cref{fig:rom_highPre_force}, corresponding to the forcing required for phase resonance, tracks the evolution of the primary resonance.

\newcommand{\HbrbForceCaption}[2]{HBRB experimental, HBM, and VPRNM ROM results for fundamental force magnitude at {#1} preload. 
From dark to light, the HBM solutions and experiments are at amplitudes of 0.5 g, 1.0 g, 2.0 g, 4.0 g and 6.0 g for the first harmonic at the shaker. The EPMC and VPRNM backbones track the primary and superharmonic resonance frequencies respectively. The superharmonic resonances, located at the VPRNM backbone (see {#2}), do not result in any notable distortion of the required external force for response control.}

\begin{figure}[h!]
	\centering
	\includegraphics[width=0.85\linewidth]{hbrb_all_legend.pdf}
	\includegraphics[width=0.49\linewidth]{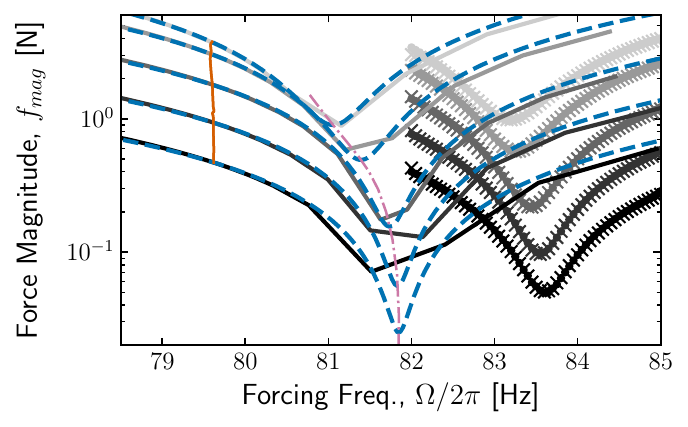}
	\caption{\HbrbForceCaption{high}{\Cref{fig:rom_highPre}} }
	\label{fig:rom_highPre_force}
\end{figure}

\FloatBarrier

Lastly, the phase of the superharmonic resonance is analyzed in \Cref{fig:rom_highPre_h7phase}. The VPRNM ROM provides a consistent phase to HBM for the 0.5 g and 1.0 g amplitude levels. However, the VPRNM ROM shows some errors for the higher amplitude levels, a behavior that is also characteristic of the neglected $\boldsymbol{F_{kq,sup,n}}$ forces \cite{porterTrackingSuperharmonic2024}.
The phase of the superharmonic resonances are also inconsistent between the modeling and the experimental results. 
This is attributed to a combination of model form error in the frictional contact nonlinearity (see \Cref{fig:hbrb_epmc_high}) and structure shaker interaction (see \Cref{fig:exp_highPre_force_h7}).
If structure shaker interactions are eliminated by advanced control \cite{hippoldShakerStructureInteraction2023}, analyzing the phase of superharmonic resonances could provide insight into the characteristics and model form errors of physics-based friction models for jointed structures.
However, such control algorithms and hardware are unavailable for this study.

Additional HBM and VPRNM results for low and medium preload levels are compared to experiments in \Cref{sec:extra_hbm_vprnm}.

\newcommand{\HbrbPhaseSevenCaption}[2]{HBRB experimental, HBM, and VPRNM ROM results for phase difference between seventh and first harmonics at {#1} preload. 
From dark to light, the HBM solutions and experiments are at amplitudes of 0.5 g (excluded for experiments), 1.0 g, 2.0 g, 4.0 g and 6.0 g for the first harmonic at the shaker.
The 0.5 g case for the experiments is excluded because of noise{#2}.}

\newcommand{\PhaseClipping}[0]{, and the plot bounds trim some experimental measurements}

\begin{figure}[h!]
	\centering
	\includegraphics[width=0.65\linewidth]{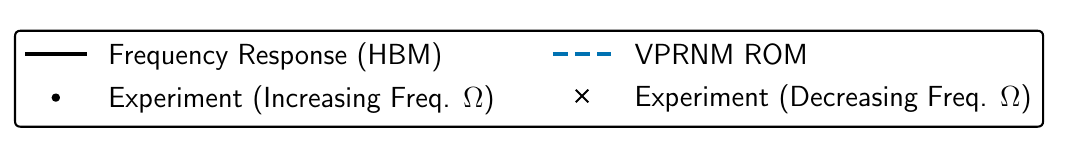}
	\includegraphics[width=\linewidth]{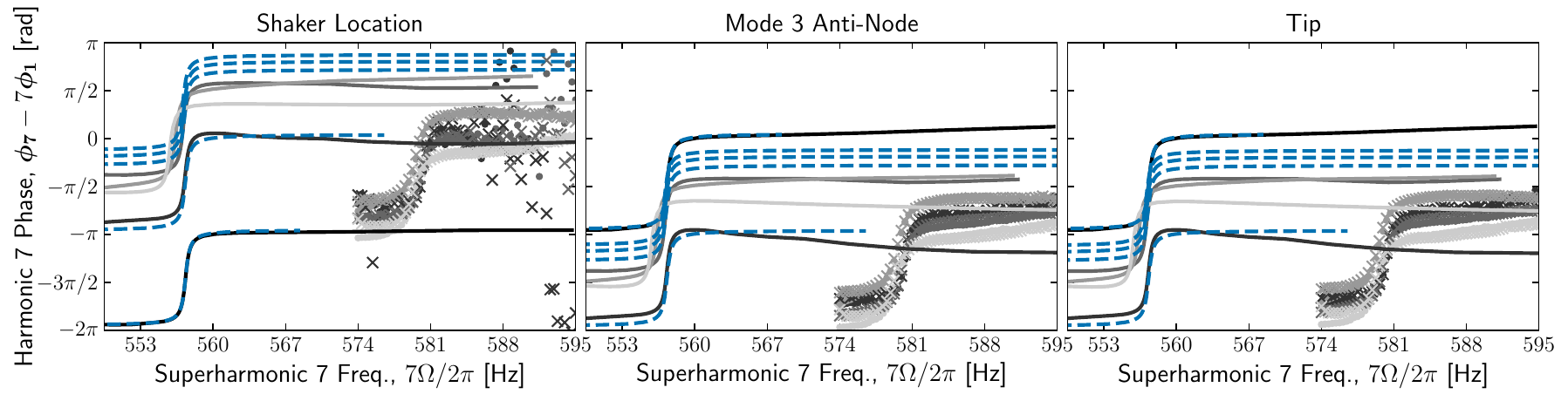}
	\caption{
		\HbrbPhaseSevenCaption{high}{\PhaseClipping}
	}
	\label{fig:rom_highPre_h7phase}
\end{figure}

\FloatBarrier

\subsubsection{Analysis of Causes of Superharmonic Resonance}
\label{sec:hbrb_sr_causes}

Based on the framework of VPRNM (see \Cref{sec:vprnm_theory}) and the insights from \Cref{sec:3dof_noSR}, the causes of the superharmonic resonance for the HBRB are analyzed. 
Specifically, this section plots the internal modal forcing due to $ \boldsymbol{F_{nq,broad}} $ from \eqref{eq:fbroad} as given by $ \boldsymbol{\eta_{vprnm}} $ in \eqref{eq:modal_force_vprnm}.
This modal forcing is calculated considering the linear mass-normalized mode shape of the third bending mode as is also done for the modal filtering of VPRNM (see \Cref{sec:theory_high_dim}). Consistent with \eqref{eq:x_t}, the phase of $ \boldsymbol{\eta_{vprnm}} $ is calculated as
\begin{equation}
	\phi_\eta = arctan2(\eta_{vprnm,s}, \eta_{vprnm,c})
\end{equation}
where $ \eta_{vprnm,c} $ and $ \eta_{vprnm,s} $ are the internal modal excitations for the cosine and sine terms respectively (corresponding to the first and second entries of the vector $ \boldsymbol{\eta_{vprnm}} $ from \eqref{eq:modal_force_vprnm}). 
In addition to calculating a total excitation of the superharmonic resonance, $ \boldsymbol{\eta_{vprnm}} $ can be calculated based on a subset of the nonlinear forces (e.g., the nonlinear forces from normal contact or within a single contact element).
In this section, the calculations are a post-processing step for previous HBM and VPRNM solutions and require no additional solutions to systems of equations.

\begin{figure}[h!]
	\centering
	\includegraphics[width=0.9\linewidth]{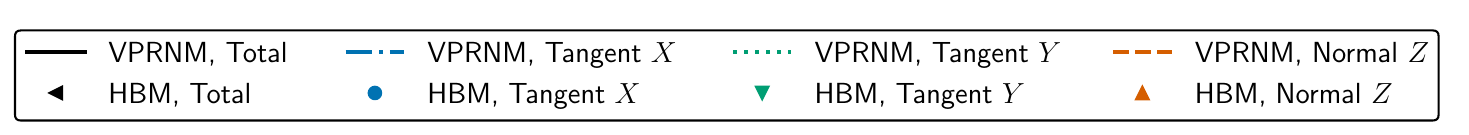}
	\begin{subfigure}[c]{0.49\linewidth } 
		\centering
		\includegraphics[width=\linewidth]{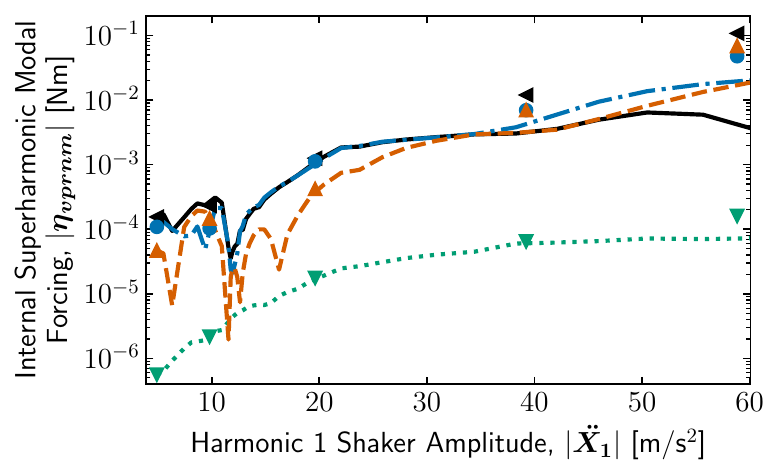}
		\caption{}
		\label{fig:cause_directions_amp}
	\end{subfigure}
	\hfill
	\begin{subfigure}[c]{0.49\linewidth } 
		\centering
		\includegraphics[width=\linewidth]{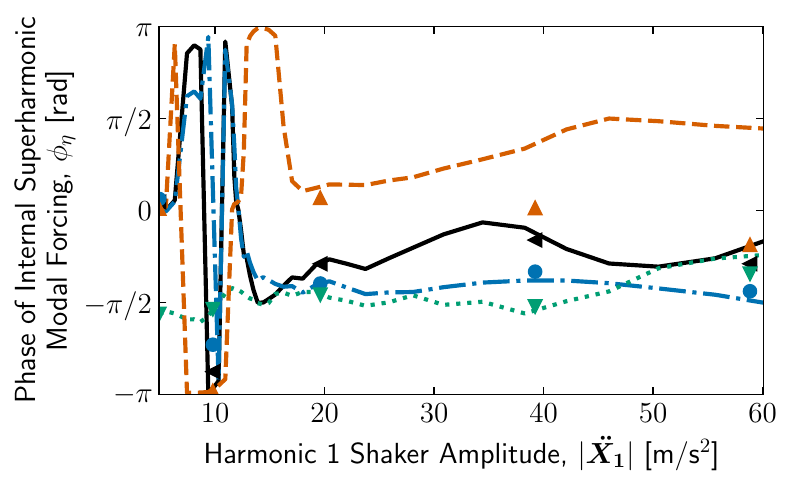}
		\caption{}
		\label{fig:cause_directions_phase}
	\end{subfigure}
		\caption{Contributions of nonlinear frictional forces in different local directions at the interface to the overall excitation of the superharmonic resonance for (a) amplitude and (b) phase. Solutions are calculated along the VPRNM backbone and for the HBM points with the highest superharmonic resonance amplitude at the tip accelerometer.}
	\label{fig:cause_directions}
\end{figure}

First, the total contributions of the nonlinear forces in each direction integrated across the interface are considered in \Cref{fig:cause_directions}.
Both the tangent $ X $ (aligned with the length of the beam) and normal $ Z $ directions provide significant contributions to the excitation of the superharmonic resonance. The VPRNM solutions are consistent with the HBM solutions at the peak superharmonic response for the 0.5~g, 1.0~g, and 2.0~g cases. 
For the 4.0 g and 6.0 g cases, there are clear differences, which are attributed to VPRNM not fully capturing the peak response point from HBM.
In \Cref{fig:cause_directions}, the rapid changes in excitation at low response amplitudes are not unexpected given that the model has significant complexity from the 232 nonlinear elements contributing to the response, each with a different mesoscale topology gap (see \Cref{fig:hbrb_surface}).
As expected, the contributions for the tangent $ Y $ direction are small because most of the motion is in the tangent $ X $ and normal $ Z $ directions for the first and third bending modes, which dominate the primary and superharmonic resonance responses respectively.

At 4.0 g and 6.0 g, the phases of the tangent $ X $ and normal $ Z $ excitations for HBM are close together, so the total excitation of the superharmonic resonance is larger than either component. Conversely, for VPRNM at 6.0 g, a phase difference of nearly $ \pi $ between the tangential $ X $ and normal $ Z $ excitations (see \Cref{fig:cause_directions_phase}) results in a decrease in the total excitation (see \Cref{fig:cause_directions_amp}).
The phase differences between VPRNM and HBM for the normal $ Z $ direction are similar to behavior observed for an SDOF system with a unilateral spring \cite{porterTrackingSuperharmonic2024}, which is a simplified representation of normal contact.
For the SDOF system, the phase of $ \boldsymbol{F_{nq,broad}} $ shifts at the superharmonic resonance causing VPRNM to poorly track the peak response amplitude. 
Likewise, the normal direction superharmonic excitation for the HBRB exhibits a phase shift along the HBM solution near the superharmonic resonance,\footnote{This phase shift is observed as a phase difference in $ \boldsymbol{\eta_{vprnm}} $ between the HBM and VPRNM solutions at a given amplitude.} likely causing the errors in VPRNM. 
In the future, this post-processing could be used to check if the normal nonlinearity significantly contributes to excitation of the superharmonic resonance and to identify amplitude levels where VPRNM may give less accurate results.

\begin{figure}[h!]
	\centering
	\includegraphics[width=0.6\linewidth]{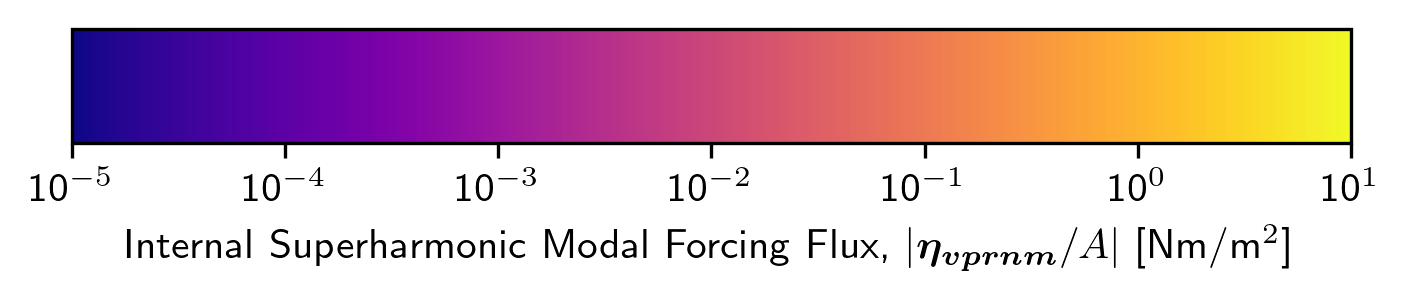}
	\\
	\begin{subfigure}[c]{0.49\linewidth } 
		\centering
		\includegraphics[width=\linewidth]{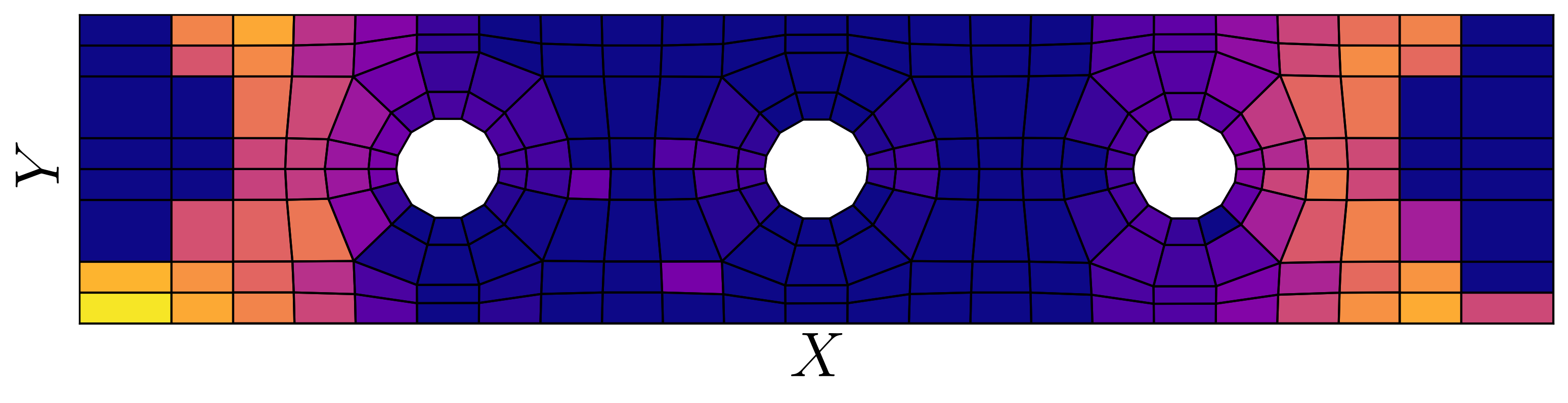}
		\caption{}
		\label{fig:phiF_mesh_lowAmpX}
	\end{subfigure}
	\hfill
	\begin{subfigure}[c]{0.49\linewidth } 
		\centering
		\includegraphics[width=\linewidth]{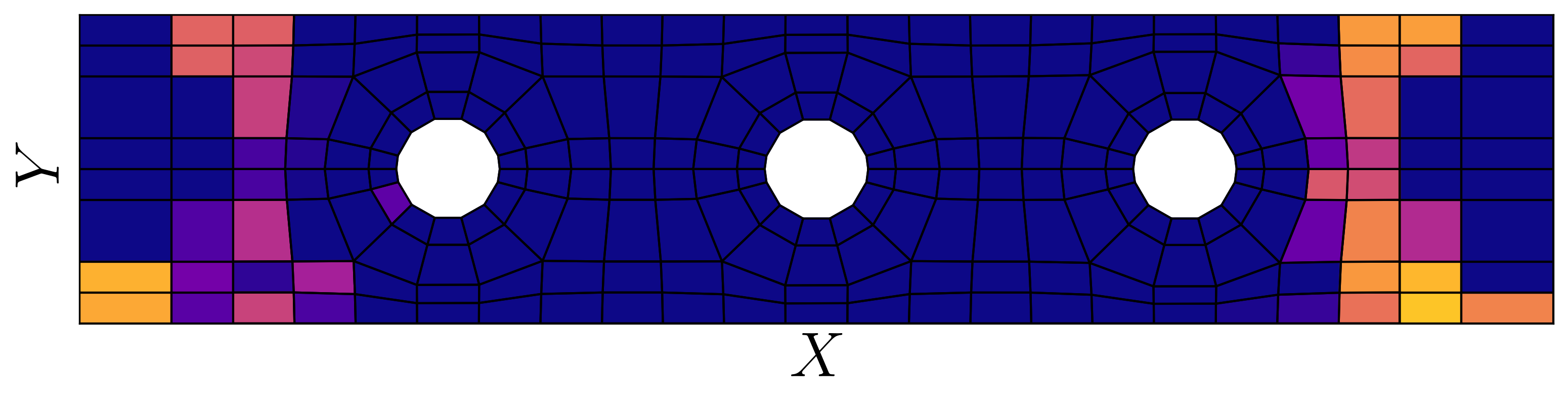}
		\caption{}
		\label{fig:phiF_mesh_lowAmpZ}
	\end{subfigure}
	\includegraphics[width=0.6\linewidth]{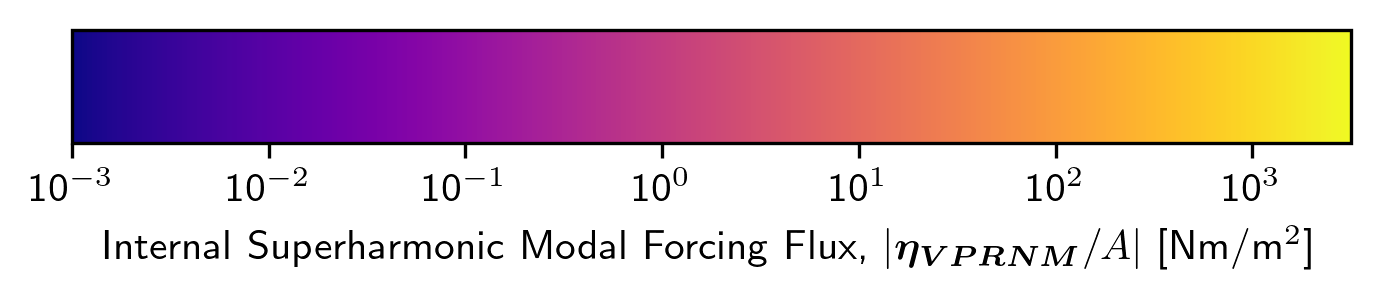}
	\\
	\begin{subfigure}[c]{0.49\linewidth } 
		\centering
		\includegraphics[width=\linewidth]{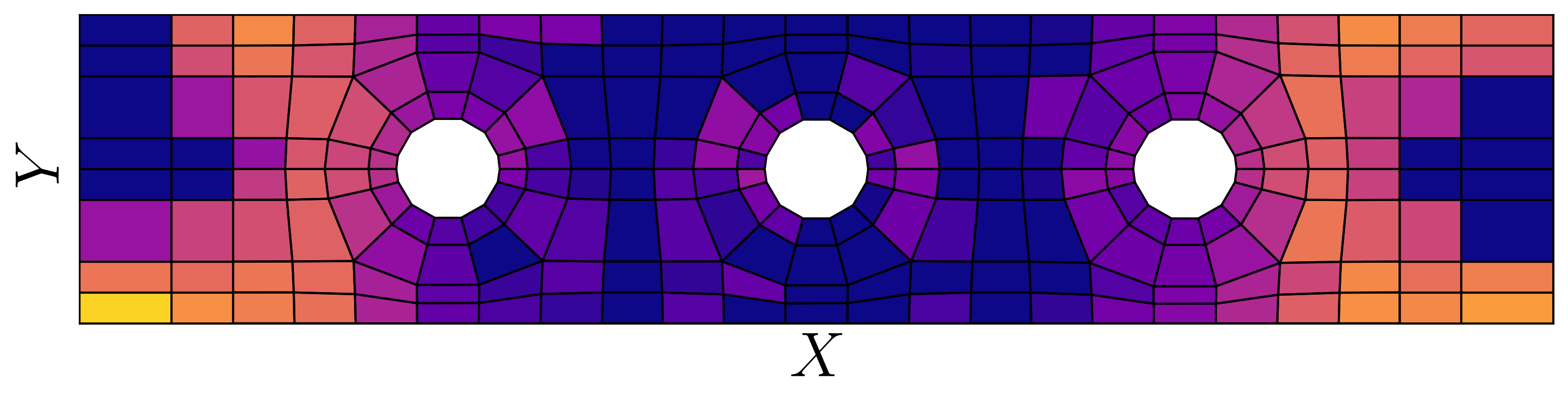}
		\caption{}
		\label{fig:phiF_mesh_highAmpX}
	\end{subfigure}
	\hfill
	\begin{subfigure}[c]{0.49\linewidth } 
		\centering
		\includegraphics[width=\linewidth]{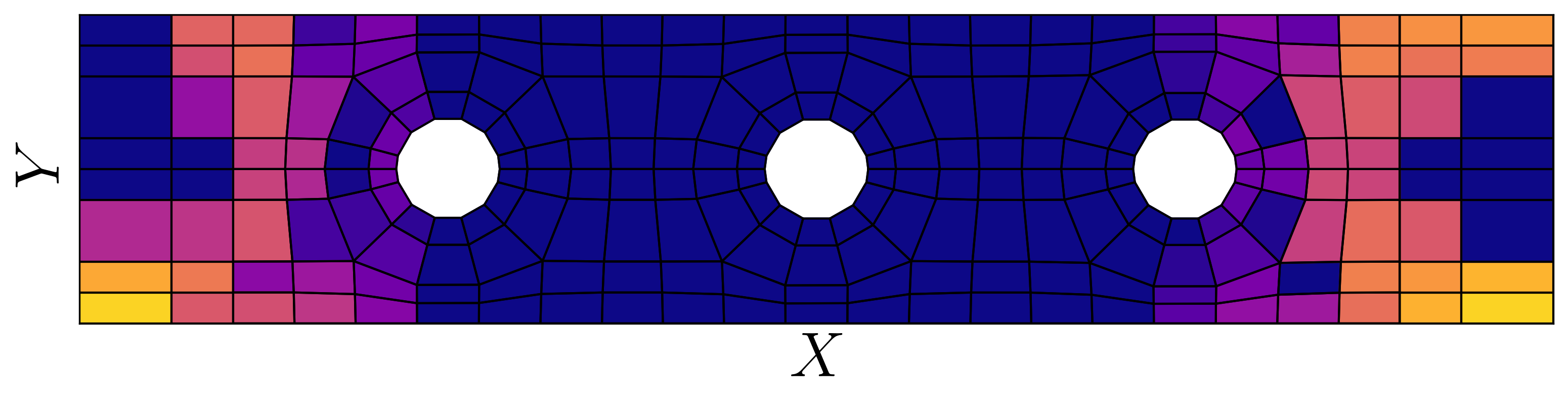}
		\caption{}
		\label{fig:phiF_mesh_highAmpZ}
	\end{subfigure}
	\caption{Spatial contributions to superharmonic resonance excitation. 
			For VPRNM at 0.5 g, plots are (a) tangent $ X $ and (b) normal $ Z $ directions. 
			For VPRNM at approximately 6.0 g (highest amplitude), plots are (c) tangent $ X $ and (d) normal $ Z $ directions. 
									For each amplitude, the color bar is shown above the plots, and values less than the lowest color bar value are shown with the darkest color.}
	\label{fig:mesh_superharmonic_excite}
\end{figure}

The spatial distribution of the superharmonic resonance excitation is shown in \Cref{fig:mesh_superharmonic_excite}. 
Here, the linear mode shape is interpolated to the quadrature points where the friction traction contributions to $ \boldsymbol{F_{nq,broad}} $ are evaluated. The product between the interpolated mode shape and the friction tractions is plotted as the internal superharmonic modal forcing flux $ |\boldsymbol{\eta_{vprnm}}/A| $.
The influence of the mesoscale topology (see \Cref{fig:hbrb_surface}) is evident in the left and right most columns of elements where the initial gap is largest and thus the contributions to the superharmonic excitation are small. 
As the nonlinearity is activated over more of the interface at the higher amplitude level, more of the interface contributes to the excitation of the superharmonic resonance (see \Cref{fig:phiF_mesh_lowAmpX} versus \Cref{fig:phiF_mesh_highAmpX}). 
Larger portions of the interface also contribute to the tangent $ X $ excitation of the superharmonic resonance compared to the normal $ Z $ excitation of the interface (see \Cref{fig:phiF_mesh_highAmpX} versus \Cref{fig:phiF_mesh_highAmpZ}). 
This is consistent with expectations that clapping and complete separation in the normal direction excite the superharmonic resonance as complete separation is expected towards the edges of the joint.

\FloatBarrier

\subsubsection{Requirement of VPRNM Modal Filter} \label{sec:results_modal_filter}

\begin{figure}[h!]
	\centering
	\includegraphics[width=0.75\linewidth]{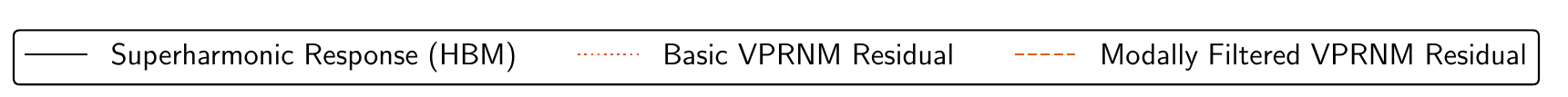}
	\begin{subfigure}[c]{0.45\linewidth } 
		\centering
		\includegraphics[width=\linewidth]{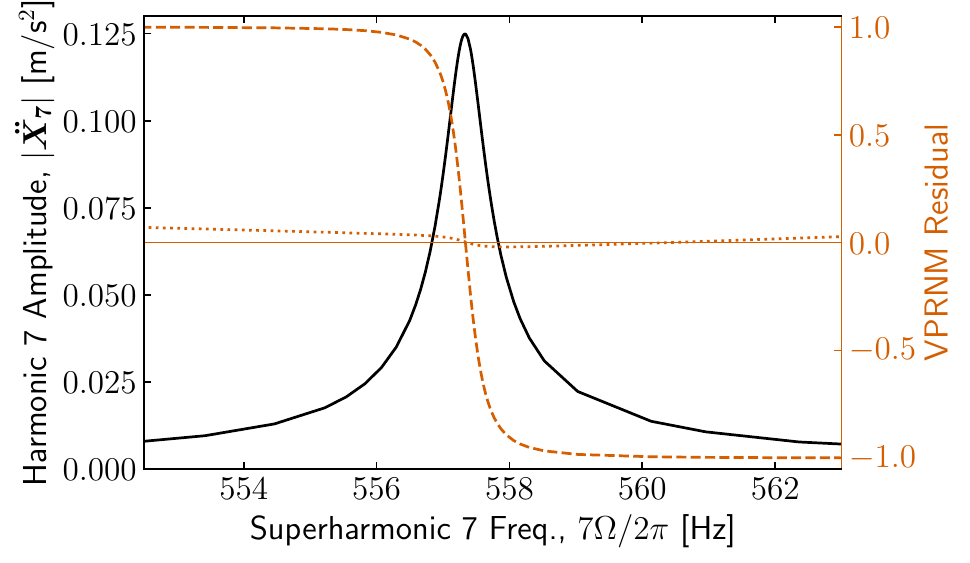}
		\caption{}
		\label{fig:vprnm_res0p5g}
	\end{subfigure}
	\quad
	\begin{subfigure}[c]{0.45\linewidth } 
		\centering
		\includegraphics[width=\linewidth]{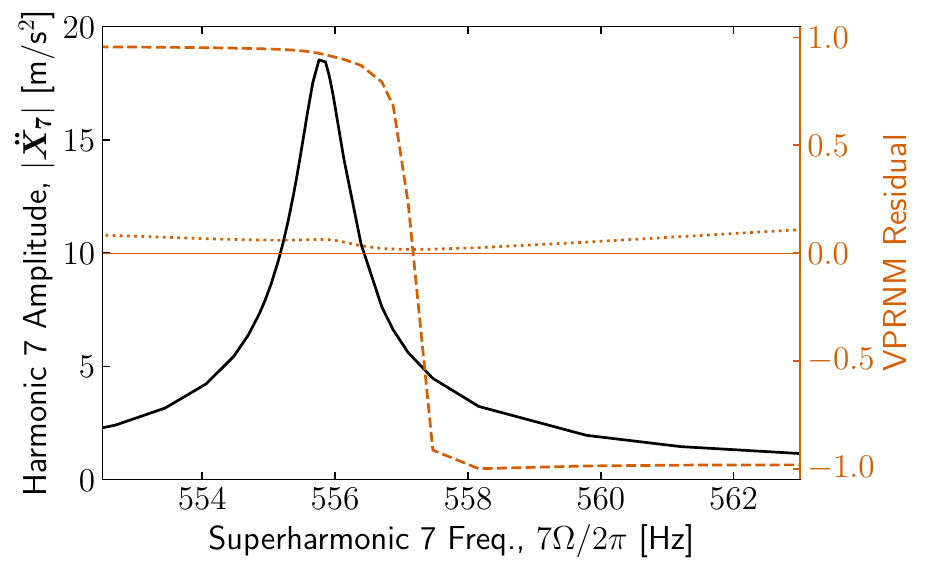}
		\caption{}
		\label{fig:vprnm_res6g}
	\end{subfigure}
	\caption{VPRNM residuals calculated along HBM solution for (a) 0.5 g and (b) 6.0 g cases with high preload. The basic VPRNM equation shows two solutions at 0.5 g near the superharmonic resonance and no solutions at 6.0 g. The modally filtered VPRNM solution shows a single solution in both cases.}
	\label{fig:vprnm_residual_plots}
\end{figure}

As introduced in \Cref{sec:theory_high_dim}, a modal filter for VPRNM is required for the simulations of the HBRB.
This section calculates the residual of the VPRNM equation along HBM solutions from \Cref{sec:hbrb_hbm_res} to provide insights into the effect of the modal filter on VPRNM as shown in \Cref{fig:vprnm_residual_plots}. 
For the 0.5 g case in \Cref{fig:vprnm_res0p5g}, VPRNM with and without the modal filter shows residuals with the same solution point corresponding to the superharmonic resonance indicating that the modal filter is not needed in this case. 
For both forms of VPRNM, a positive residual can be interpreted as the superharmonic responding in-phase with the excitation $\boldsymbol{F_{nq,broad}}$ in \eqref{eq:vprnm_rewrite_hbm}. At the solution point, the residual transitions from positive to negative corresponding to a transition to out-of-phase motion by analogy to linear systems.
For the 0.5 g case, the basic VPRNM residual becomes positive again after 560 Hz, likely due to a transition between resonance phenomena. On the other hand, the modally filtered VPRNM residual remains negative for the plotted regime.
The small magnitude of the basic VPRNM residual in \Cref{fig:vprnm_residual_plots} indicates that much of the motion is orthogonal to the excitation at all frequencies.\footnote{The modally filtered VPRNM residual does not provide such insights given the different normalization.} This is not surprising given that the nonlinear forces are only located at the interface where the vibration displacements are small.

For the 6.0 g case in \Cref{fig:vprnm_res6g}, the basic VPRNM residual has no solutions over the considered frequency regime, so the modally filtered VPRNM equation must be used.
Here, it appears that the secondary resonance phenomena that causes a shift from a negative to positive VPRNM residual in \Cref{fig:vprnm_res0p5g} has become strong enough to cover up the superharmonic resonance for the basic VPRNM equation.
It is important to emphasize that the VPRNM equation only has nonzero forces $\boldsymbol{F_{nq,broad}}$ at the nonlinearity in the interface and thus the basic VPRNM equation only has contributions from the small local displacements at the interface. 
Therefore, other resonance phenomena can easily distort these small displacements while the response near the tip of the beam shows a clear superharmonic resonance.
Likewise, other work has suggested that excitation location must be carefully considered to ensure that a single mode can be isolated for velocity feedback testing \cite{scheelNonlinearModalTesting2022}, but here the excitation is due to the internal nonlinear forces, so it cannot be changed, requiring the use of modal filtering.
These effects could also explain the errors in the VPRNM backbone and VPRNM ROM for the higher amplitude levels in \Cref{sec:hbrb_hbm_res}.
The potential for superharmonic motion to get hidden by other vibration phenomena was also partially observed in the reduced phase shift (less than $\pi$ crossing the superharmonic resonance) of DOFs 2 and 3 in \Cref{fig:3dof_h3phase}.

\FloatBarrier

\subsection{Computation Time}
\FloatBarrier

\Cref{tab:computation_time_hbrb} provides an overview of the computation times for the HBRB.
The VPRNM ROM construction takes 24\% longer than the average time to calculate a single HBM solution.\footnote{Including computation time for the continuation with respect to amplitude to initialize the HBM solutions at higher amplitudes.} However, comparing to the time to calculate all five HBM solutions, the VPRNM ROM construction is 4 times faster than HBM.
Once constructed, the VPRNM ROM is over 780,000 times faster than HBM for calculating responses at new amplitude levels.
Furthermore, VPRNM ROM construction provides higher resolution force magnitudes around the primary resonance of the first mode compared to significant faceting with HBM (see \Cref{fig:rom_highPre_force}). 
Additionally, the EPMC backbone for the third bending mode can be used to understand the primary resonance of the third bending mode, which is not considered with HBM here.
To resolve the primary resonances of the first and third bending modes with HBM, the total computational cost is expected to increase by a factor of 2 to 3 while no additional calculations are needed for VPRNM and EPMC ROM construction. Thus, the VPRNM and EPMC ROM construction would be 8-12 times faster than HBM to characterize the two primary resonances and one superharmonic resonance.

\begin{table}[h!]
	\centering
	\caption{Computation time for HBRB at high preload. In the third column the continuation time is normalized by the static solution time for the given run. 
	ROM evaluation is for all HBM amplitude levels without parallelization. 
	HBM (harmonics 0 and 1-9), EPMC (harmonics 0 and 1-3), and VPRNM (harmonics 0 and 1-9) simulations were run on individual nodes of a heterogeneous cluster (16-40 cores, 32 GB-1.5TB of RAM, 2.10-2.60 GHz processor) and the static solution time was timed during each run as a reference for how the speed varied between nodes.
	The VPRNM ROM utilizes the same desktop computer as \Cref{tab:3dof_time} (6 core Intel i7-10710U CPU, 32 GB of RAM, 1.10 GHz processor).
	The HBM continuation varying the amplitude provides initial guesses for the simulations at higher amplitude constant accelerations (see \Cref{sec:solver_consider}).
	}
	\label{tab:computation_time_hbrb}
	\resizebox{\textwidth}{!}{ 	\begin{tabular}{cccc}
		\hline\hline 
		& Static (sec) & Continuation (hr) & Norm. Continuation (hr/sec)  \\ \hline \hline
		\multicolumn{4}{l}{Backbone and Initialization Calculations}\\ \hline\hline
		HBM Varying Amplitude & 20.4 & 1.70 & 301 \\
		EPMC Mode 1 & 19.7 & 0.192 & 35.2 \\
		EPMC Mode 3 & 14.2 & 0.220 & 55.6 \\
		VPRNM & 19.9 & 9.13 & 1650 \\
		VPRNM ROM Construction & NA & 9.55 & 1740 \\
		\hline\hline
		\multicolumn{4}{l}{FRC Calculations}\\ \hline
		HBM 0.5 g & 19.4 & 8.38 & 1560 \\
		HBM 1.0 g & 20.9 & 10.1 & 1730 \\
		HBM 2.0 g & 19.0 & 10.6 & 2010 \\
		HBM 4.0 g & 20.1 & 4.30 & 771 \\
		HBM 6.0 g & 19.2 & 3.53 & 661 \\
		HBM Total (Including Initialization) & NA & 38.5 & 7030 \\
		VPRNM ROM Evaluation (All Levels) & NA & 4.86e-5 (0.175 sec)  & 8.97e-3 \\
				\hline\hline
	\end{tabular}}
\end{table}

\FloatBarrier

\section{Conclusions} \label{sec:mdof_conclusions}

The major contributions of this work are the extension of variable phase resonance nonlinear modes (VPRNM) to superharmonic resonances in high dimensional systems and the creation of a new reduced order model based on VPRNM (VPRNM ROM) for efficiently calculating superharmonic resonance behavior.
This modeling approach provides significant computational speedups compared to the harmonic balance method (HBM) and provides insights into the phenomena of superharmonic resonances for both the Half Brake-Reu{\ss} Beam (HBRB) and a 3 degree of freedom (DOF) system.
Specific insights into superharmonic resonances are
\begin{itemize}
	\item Superharmonic resonances evolve with amplitude and must be characterized at multiple amplitude levels.
	\item Amplitude control of the first harmonic allows for the superharmonic resonance to be treated with single nonlinear mode theory through the VPRNM ROM.
	\item Compatibility between the excitation of the superharmonic and the mode shape for a potential superharmonic resonance is required for a superharmonic resonance to occur. 
	This compatibility can be expressed as nonzero modal forcing from $\boldsymbol{F_{nq,broad}}$ (see \eqref{eq:fbroad}).
		\item The modeling results of the HBRB suggest that the nonlinear forces in both the normal (from clapping) and tangential (from friction) directions play an important role in exciting the superharmonic resonance. Furthermore, VPRNM shows higher errors when the normal direction forces are more important, which is consistent with the single DOF analyses in \cite{porterTrackingSuperharmonic2024}.
\end{itemize}
In addition to these insights into superharmonic resonances, this work contributes:
\begin{itemize}
	\item The VPRNM ROM that efficiently captures superharmonic resonance behavior with speedups of up to 4 times compared to HBM for VPRNM ROM construction and speedups of up to 780,000 times compared to HBM for VPRNM ROM evaluation.
	\item Further validation of the physics-based friction model proposed in \cite{porterPredictivePhysicsbasedFriction2023} considering a second jointed structure, and excellent agreement is found for the prediction of linear modal frequencies.
	\item New experimental measurements of the HBRB utilizing instrumented bolts to control preload, shaker ring down to isolate individual modes, and stepped-sine testing to identify superharmonic resonances.
\end{itemize}

Future work could address limitations of VPRNM at some amplitude levels and further investigate applications to other systems.
Additional investigation and validation of the physics-based friction model from \cite{porterPredictivePhysicsbasedFriction2023} is necessary to understand errors in the regime with less nonlinearity investigated here.
Lastly, the VPRNM ROM could be adapted to extrema tracking \cite{razeTrackingAmplitudeExtrema2024} to provide insights into more systems. Such an adaptation would also require addressing second derivatives for the present piecewise linear friction models and concerns about computation time.

Through the development of the VPRNM ROM and the investigation of both a 3 DOF system and a real jointed connection, this work provides significant insights into superharmonic resonances. 
Python codes for simulations and applications of the VPRNM ROM are made available \cite{porterTMDSimPy}.

\section*{Acknowledgments}
\label{sec:acknowledgments}

Funding: 
This material is based upon work
supported by the U.S. Department of Energy, Office of
Science, Office of Advanced Scientific Computing
Research, Department of Energy Computational Science
Graduate Fellowship under Award Number(s) DE-SC0021110.
This work was supported in part by the Big-Data Private-Cloud Research Cyberinfrastructure MRI-award funded by NSF under grant CNS-1338099 and by Rice University. The authors are thankful for the support of the National
Science Foundation under Grant Number 1847130.

This report was prepared as an account of
work sponsored by an agency of the United States
Government. Neither the United States Government nor
any agency thereof, nor any of their employees, makes
any warranty, express or implied, or assumes any legal
liability or responsibility for the accuracy, completeness,
or usefulness of any information, apparatus, product, or
process disclosed, or represents that its use would not
infringe privately owned rights. Reference herein to any
specific commercial product, process, or service by trade
name, trademark, manufacturer, or otherwise does not
necessarily constitute or imply its endorsement,
recommendation, or favoring by the United States
Government or any agency thereof. The views and
opinions of authors expressed herein do not necessarily
state or reflect those of the United States Government or
any agency thereof.

\printbibliography

\appendix

\section{Additional Nonlinear Solution Considerations}

\subsection{Prestress Solution} \label{sec:prestress_theory}

For the HBRB, static forces in the form of bolt tensions are applied (see \Cref{sec:model_procedure} for modeling details). As an initial step and as a basis of the initial guesses (see \Cref{sec:solver_consider}), a static displacement solution $ \boldsymbol{x_s} $ is found prior to solving the full equations of motion \eqref{eq:eom} by solving
\begin{equation}
	\boldsymbol{K} \boldsymbol{x_s} + \boldsymbol{T} \boldsymbol{f_{nl}}\left(\boldsymbol{Qx_s}\right) - \boldsymbol{F_{ext, 0}} = \boldsymbol{0}	.
\end{equation}
The linearized stiffness around $ \boldsymbol{x_s} $ is used for low amplitude linearized solutions and to set mass and stiffness proportional damping (e.g., see \cite{porterPredictivePhysicsbasedFriction2023}).
For the prestress solution, all tangential friction forces are set to zero as is done in previous work \cite{porterQuantitativeAssessmentModel2022, porterPredictivePhysicsbasedFriction2023}.
For subsequent analyses utilizing AFT, the static solution displacements are used to initialize the friction sliders at that displacement and zero traction. This is an approach to select a unique residual traction field, but an appropriate choice for these residual tractions remains an open area of research \cite{ferhatogluNonuniquenessFrictionForces2021, ferhatogluFrequencyResponseVariability2023}.
For jointed structures such as the HBRB, the effect of residual tractions is expected to be small since the structure lacks the extent of tangential-normal coupling seen in \cite{ferhatogluFrequencyResponseVariability2023}.

\subsection{Extended Periodic Motion Concept (EPMC)} \label{sec:epmc_theory}

The extended periodic motion concept (EPMC) \cite{krackNonlinearModalAnalysis2015} is utilized to calculate nonlinear modal frequency and damping properties for an isolated mode in Sections \ref{sec:3dof_description}, \ref{sec:3dof_noSR}, and \ref{sec:hbrb_epmc_results}. 
EPMC considers harmonic equations similar to HBM (see \Cref{sec:theory_hbm}) with negative mass proportional damping to act as self excitation instead of harmonic forcing. Thus the unknowns at a given modal amplitude $ q_{i} $ are the harmonic displacements, the modal frequency $ \omega_{i} $, and the self excitation factor $ \xi_{i} $
\begin{subequations}
\begin{equation}
	\boldsymbol{K} \boldsymbol{X_0} + \boldsymbol{F_{nl, 0}} - \boldsymbol{F_{ext, 0}} = \boldsymbol{0}
\end{equation}
\begin{equation}
	(-n \omega_{i}^2 \boldsymbol{M} + \boldsymbol{K}) \boldsymbol{X_{nc}} + n \omega_{i} (\boldsymbol{C} - \xi_{i} \boldsymbol{M}) \boldsymbol{X_{ns}} + + \boldsymbol{F_{nl, nc}} = \boldsymbol{0} \ \ 
	\forall n \in \{1, \dots, H\}
\end{equation}
\begin{equation}
	(-n \omega_{i}^2 \boldsymbol{M} + \boldsymbol{K}) \boldsymbol{X_{ns}} - n \omega_{i} (\boldsymbol{C}- \xi_{i} \boldsymbol{M}) \boldsymbol{X_{nc}} + \boldsymbol{F_{nl, ns}} = \boldsymbol{0} \ \ 
	\forall n \in \{1, \dots, H\}	
\end{equation}
\begin{equation} \label{eq:EPMC_phase}
	\boldsymbol{R} \boldsymbol{X_{1c}} = 0
\end{equation}
\begin{equation}
	\boldsymbol{X_{1c}}^T \boldsymbol{M}\boldsymbol{X_{1c}} + \boldsymbol{X_{1s}}^T \boldsymbol{M}\boldsymbol{X_{1s}} 
	-q_{i}^2
	= 0
\end{equation}
\end{subequations}
Here, \eqref{eq:EPMC_phase} is an arbitrary phase constraint defined with an $ 1 \times N $ vector $ \boldsymbol{R} $.
The self excitation factor can be converted to a fraction of critical damping factor as
\begin{equation}
	\zeta_{i}(q_{i})  = \dfrac{ \xi_{i}(q_{i}) }{2 \omega_{i}(q_{i})}.
\end{equation}
For EPMC, continuation (see \Cref{sec:theory_hbm}) is done with respect to the modal amplitude $ q_{i} $.
The exact implementation of EPMC is available at \cite{porterTMDSimPy}.

\subsection{Solver Considerations for Large Models} \label{sec:solver_consider}

The HBRB system (see \Cref{sec:hbrb_system}) has 860 DOFs, making it more challenging to solve than the 3 DOF system. The present section provides additional insights into the solver configurations and initial guesses that were used throughout \Cref{sec:hbrb_results}.
For all of the HBM and VPRNM solutions for the HBRB, amplitude control is applied to the acceleration (i.e., $ k=2 $ in \eqref{eq:amp_phase_constraint}) at the shaker location to match the experimental tests.

For HBM at the lowest amplitude levels, initial guesses are constructed based on the linearized frequency response around the static solution (see \Cref{sec:prestress_theory}). The phase of this solution is rotated to match the phase constraint of \eqref{eq:amp_phase_constraint}, and the amplitude is appropriately scaled. 
For higher acceleration levels, this linear initial guess did not always converge. Therefore, continuation is conducted at a constant frequency ($ \Omega= 490$ rad/s, 78 Hz) with respect to $ A_1 $ in \eqref{eq:amp_phase_constraint} to obtain initial solutions at higher amplitude levels. The point from this continuation with the closest amplitude level to the desired acceleration is used to initialize the continuation with respect to frequency. 
In practice, this continuation was relatively coarse and did not have too much trouble converging even when the nearest solution point from the amplitude continuation was approximately 10\% different than the desired amplitude level.

For VPRNM, the $ n=7 $ superharmonic resonance of the third bending mode with frequency $ \omega_{3} $ is of interest, so the initial guess is constructed at forcing frequency 
\begin{equation}
	\Omega = \omega_{3} / n.
\end{equation}
As with the HBM solution, the linearized model is utilized to calculate a frequency response of the first harmonic, which is scaled and rotated to satisfy \eqref{eq:amp_phase_constraint}. 
Then, the nonlinear forces exciting harmonic $ n $ are calculated using AFT and the predicted displacements for the zeroth and first harmonics.
These nonlinear forces in the form of $ \boldsymbol{F_{nq,broad}} $ are then treated as external forcing of the higher harmonic and are used with the linearized system to calculate an initial guess of the response of the harmonic $ n $.
Thus, the initial guess for VPRNM includes the prestress static displacements, the linear frequency response displacements for the first harmonic based on the external excitation, and linear frequency response displacements for the superharmonic number $ n $ based on the internal nonlinear forces $ \boldsymbol{F_{nq,broad}} $.

A gradient based solver is utilized for the prestress solution, initial solution to the HBM or VPRNM equations, and for each continuation step. 
The exact solver implementations are available at \cite{porterTMDSimPy}.
For the prestress solution and the initial point, a full Newton-Raphson scheme is utilized. 
For VPRNM, line search similar to \cite{matthiesSolutionNonlinearFinite1979} is necessary for the initial point since individual steps tended to overshoot the solution and diverge when line search was not utilized.
During continuation, line search is turned off for VPRNM since it requires additional computation time and can require more iterations to find the solution.

For continuation of the HBM and EPMC solutions, the BFGS algorithm (named for Broyden, Fletcher, Goldfarb, and Shanno) based on algorithm 7.4 of \cite{nocedalNumericalOptimization2006} is implemented. 
BFGS is a quasi-Newton approach with a rank two update to the Jacobian on BFGS steps. For the present problems, the full Jacobian matrix is calculated and factored every other step starting from the first step. Using BFGS on alternating steps, the low rank update to the Jacobian is utilized to determine a new step for the nonlinear solver. 
The benefit of BFGS is that the Jacobian need not be evaluated or factored at every step. In practice, calculation of the Jacobian for the HBRB model is on the order of 10-100 times more expensive than calculating just the residual vector, so the BFGS iterations are much faster than the full Newton-Raphson iterations. Asymptotically for problems with $ M $ unknowns, BFGS iterations do not require refactoring the Jacobian matrix, replacing an $ \mathcal{O}(M^3) $ matrix factoring operation with $ \mathcal{O}( M^2) $ vector operations.
For VPRNM, initial tests indicated that while BFGS could converge when the modal filter was applied, the simulation time increased because more iterations and smaller continuation steps were required.

\section{Summary of VPRNM ROM Steps} \label{sec:rom_summary_steps}

This section provides an outline of the steps for the VPRNM ROM without the derivation details provided in \Cref{sec:full_rom_theory}. 
The steps are
\begin{enumerate}
	\item Interpolate the VPRNM solution to the desired amplitude based on \Cref{eq:interp_vprnm}.
	
	\item Interpolate the EPMC solution for the mode corresponding to the superharmonic resonance to match the VPRNM amplitude with \eqref{eq:interp_super_EPMC}.
	
	\item Approximate the modal force exciting the superharmonic resonance with \eqref{eq:superharmonic_force}.
	
	\item Utilize the single nonlinear mode EPMC ROM with phase information described in \Cref{sec:phase_amp_rom} to calculate the superharmonic resonance response for harmonic $ n $. Note that the forcing frequencies in \eqref{eq:epmc_amp_only_rom} are the frequency of the superharmonic resonance $ \Omega_S $ for this step corresponding to $ n $ times the external forcing frequencies $ \Omega $. Also, the nonlinear modal frequency of the superharmonic mode is scaled as given by \eqref{eq:scale_superharmonic_mode_freq} to match the VPRNM solution. Modal amplitudes of the superharmonic response are also saved at this point as $ q_S(\Omega_S) $.
	
	\item Rotate the superharmonic EPMC ROM from the previous step to match the VPRNM phase using \eqref{eq:super_phase_diff} and \eqref{eq:super_phase_rot}.

	\item Interpolate the EPMC solution for the fundamental mode to the desired amplitude with \eqref{eq:interp_epmc_fund}. Note that this EPMC backbone solution should be calculated excluding the harmonic corresponding to the superharmonic resonance.
	
	\item Rotate the EPMC mode for the single interpolated point of the fundamental mode to match the VPRNM phase with \eqref{eq:fund_phase_diff} and \eqref{eq:fund_phase_rot}.
	
	\item Combine the displacements from the zeroth harmonic of the interpolated VPRNM solution, the fundamental interpolated and rotated EPMC point from \eqref{eq:fund_phase_rot}, and the rotated superharmonic response from \eqref{eq:super_phase_rot}.
	
	\item Calculate the force to achieve constant amplitude for the fundamental EPMC ROM utilizing \eqref{eq:epmc_amp_control} at external forcing frequencies of $ \Omega = \Omega_S / n $.
	
	\item Optionally, calculate a force correction to the fundamental external forcing magnitude with 
	\eqref{eq:rom_deltaF} and \eqref{eq:vprnm_rom_force}.
\end{enumerate}

\section{3 Degree of Freedom System Formulation} \label{sec:3dof_deriv}

This section describes how physical mass and stiffness matrices are constructed to achieve the desired mode shapes and modal frequencies of the 3 DOF system used in \Cref{sec:3dof_sys}.
The desired mass-normalized mode shapes for the system are
\begin{equation}
	\phi_1 = \begin{bmatrix}
		1 \\ 2 \\ 3 
	\end{bmatrix}, 
	\ \ 
	\phi_a = \begin{bmatrix}
		2 \\ 1 \\ -1 
	\end{bmatrix}, 
	\ \ 
	\phi_b = \begin{bmatrix}
		-2 \\ 1 \\ 1 
	\end{bmatrix}
\end{equation}
where $ \phi_1 $ is the first mode, $ \phi_a $ is the second mode shape for the first 3 DOF example, and $ \phi_b $ is the third mode shape for the first 3 DOF example (see \Cref{sec:3dof_description}). For the second 3 DOF example (see \Cref{sec:3dof_noSR}), $ \phi_b $ is the second mode shape and $ \phi_a $ is the third mode shape. 
For each system, the mode shape vectors are arranged in order as column vectors in the mode shape matrix $ \boldsymbol{\Phi} $, and the modal frequencies are chosen to be $ \omega = $ 1.0 rad/s, 3.0 rad/s, and 7.5 rad/s.
These frequencies are put in the spectral matrix as
\begin{equation}
	\boldsymbol{\Lambda} = 
	\begin{bmatrix}
		1.0^2 & 0 & 0 \\
		0 & 3.0^2 & 0 \\
		0 & 0 & 7.5^2
	\end{bmatrix}.
\end{equation}
The mass and stiffness matrices to give the desired properties are then
\begin{equation}
	\boldsymbol{M} = 
	\boldsymbol{\Phi}^{-T} \boldsymbol{\Phi}^{-1} 
	, \ \ \ \
	\boldsymbol{K_{tuned}} = 
	\boldsymbol{\Phi}^{-T} \boldsymbol{\Lambda} \boldsymbol{\Phi}^{-1} .
\end{equation}
Mass proportional damping with factor of 0.01 is used
\begin{equation}
	\boldsymbol{C} = 0.01 \boldsymbol{M}.
\end{equation}
The Iwan element nonlinearity (see next) contributes a linear stiffness $ k_t $ at low amplitude and zero stiffness in the limit of infinite amplitude. Therefore, half of the stiffness of the Iwan element is subtracted from $ \boldsymbol{K_{tuned}} $ such that the desired frequency and mode shapes will be approximately achieved somewhere in the nonlinear regime. 
Therefore the purely linear stiffness matrix is calculated as
\begin{equation}
	\boldsymbol{K}
	=
	\boldsymbol{K_{tuned}}
	- 
	(0.5)
	\boldsymbol{T} k_t \boldsymbol{Q}
	=
	\boldsymbol{K_{tuned}}
	- 
	0.5
	\begin{bmatrix}
		0 & 0 & 0 \\
		0 &  k_t & -k_t \\
		0 & -k_t &  k_t
	\end{bmatrix}
\end{equation}
using the $ \boldsymbol{Q} $ and $ \boldsymbol{T} $ from \Cref{sec:3dof_sys}.

\subsection{Four-Parameter Iwan Model} \label{sec:iwan_details}

The four-parameter Iwan element has a distribution of sliders of
\cite{segalmanFourParameterIwanModel2005} (the implementation of \cite{porterTrackingSuperharmonic2024, porterTMDSimPy} is used)
\begin{equation} \label{eq:mdof_iwan_distrib}
	\rho(\phi) = \dfrac{F_s (\chi + 1)}{\phi_{max}^{\chi + 2}  \left( \beta + \frac{\chi + 1}{\chi + 2} \right) } \phi^\chi 
		+ 
		\dfrac{F_s \beta}{\phi_{max}  \left( \beta + \frac{\chi + 1}{\chi + 2} \right) } 
	\delta(\phi-\phi_{max})
\end{equation}
where $ \delta(\cdot) $ is the Dirac delta function, $ \phi $ is the displacement when a given slider starts to slip, $ \phi_{max} $ is the displacement that causes all sliders to start slipping and is calculated as 
\begin{equation} \label{eq:mdof_phi_max}
	\phi_{max} = \dfrac{F_s (1 + \beta)}{k_t \left( \beta + \frac{\chi + 1}{\chi + 2} \right)}.
\end{equation}
Additionally, $ F_s $ describes the total force required to cause all sliders to slip and $ \chi $ and $ \beta $ control the shape of the hysteresis loop.
The Iwan element in the 3 DOF model uses parameters of $ k_t = 0.6 $ N/m, $ F_s = 10.0 $~N, $ \chi = -0.5 $, $ \beta = 0 $, and 100 discrete sliders.
Each discrete slider has an associated force at each instance of time as given by
\begin{equation} \label{eq:mdof_slider_f_iwan}
	f_\phi = \begin{cases}
		\underbrace{x - x_0 + f_{\phi,0}  }_{f_{\phi,stuck}} & |f_{\phi,stuck}| < \phi
		\\
		\phi \text{sgn}(f_{\phi,stuck}) & Otherwise
	\end{cases}
\end{equation}
where $ x $ is the current displacement, $ x_0 $ is the displacement at the previous instant, and $ f_{\phi,0} $ is the force of the given slider at the previous instant in time. 
The total force from all of the sliders is calculated as
\begin{equation}\label{eq:mdof_iwan_integral}
	f_{nl} = 
	\int_{0}^{\phi_{max}}   f_\phi	\rho(\phi) d\phi  .
\end{equation}

\section{Additional Half Brake-Reu{\ss} Beam Results}

Supplementing \Cref{sec:hbrb_results}, addition results for the HBRB are provided here. First, \Cref{sec:extra_epmc_appendix} provides additional EPMC backbones for medium and low preload levels. 
Then, additional stepped sine results are provided (\Cref{sec:extra_stepped_sine_appendix}). HBM convergence is shown in \Cref{sec:hbm_convg} and additional HBM and VPRNM ROM results are shown in \Cref{sec:extra_hbm_vprnm}.

\subsection{Additional EPMC Results} \label{sec:extra_epmc_appendix}

\begin{figure}[h!]
	\centering
	\includegraphics[width=0.55\linewidth]{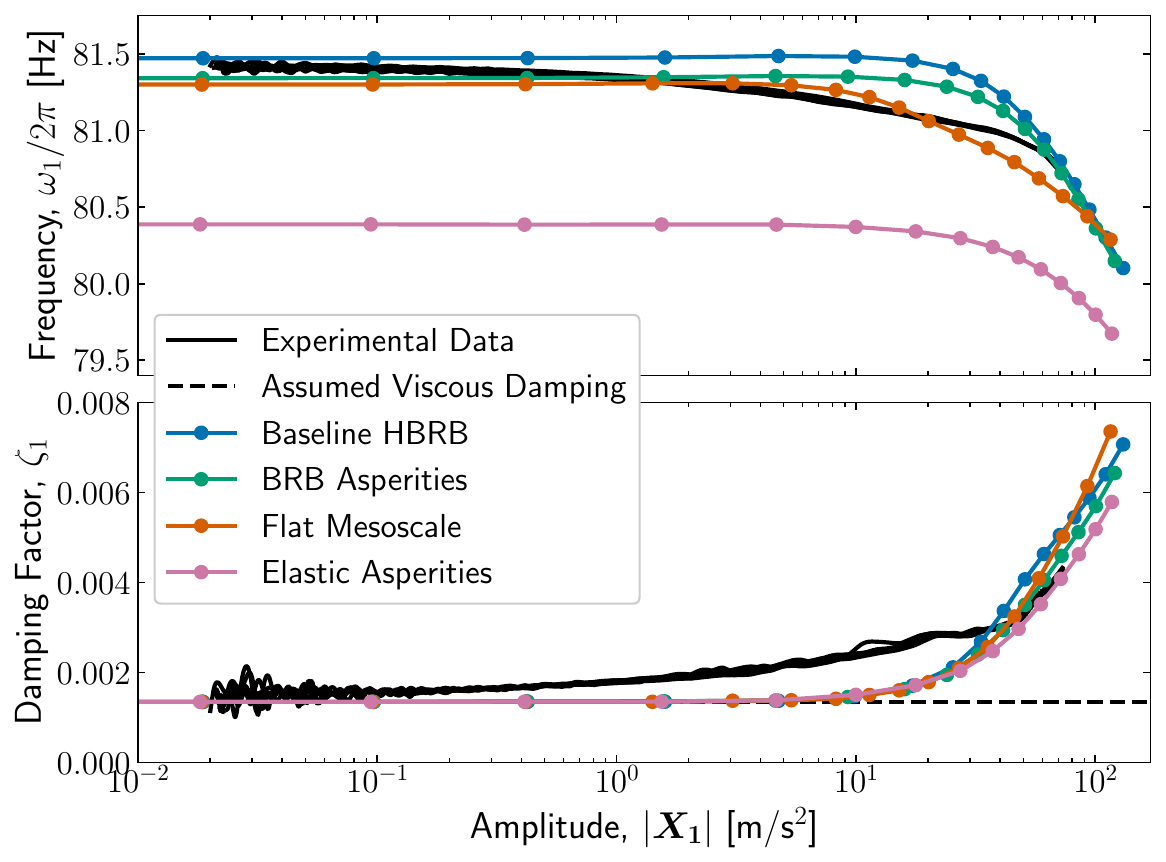}
	\caption{\EpmcPhysicsCaption{low}
	}
	\label{fig:hbrb_epmc_low_physics}
\end{figure}

\begin{figure}[h!]
	\centering
	\includegraphics[width=0.55\linewidth]{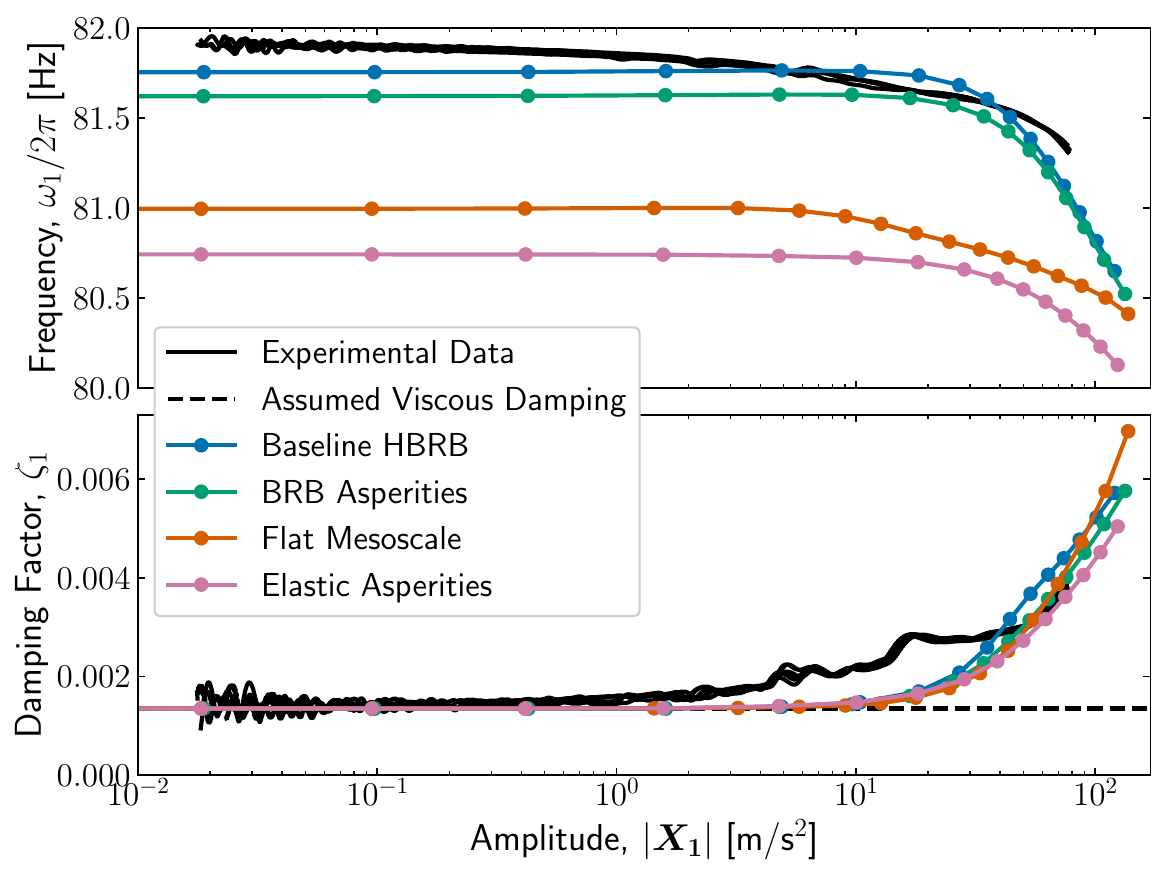}
	\caption{\EpmcPhysicsCaption{medium}
	}
	\label{fig:hbrb_epmc_medium_physics}
\end{figure}

Figures \ref{fig:hbrb_epmc_low_physics} and \ref{fig:hbrb_epmc_medium_physics} show the different effects of neglecting parts of the friction model physics for low and medium preload levels respectively. Here, the baseline HBRB model corresponds to all of the physics included (e.g., the model used in \Cref{sec:hbrb_results}). The `BRB Asperities' case uses the asperity properties from \cite{porterPredictivePhysicsbasedFriction2023} rather than those identified in this work. Given the good prediction quality, utilizing statistics from other structures manufactured in the same way is reasonable. 
For both preload values, using purely elastic asperities resulted in a decrease in frequency as previously observed in \cite{porterPredictivePhysicsbasedFriction2023}.
Lastly, considering only a flat mesoscale topology (`Flat') results in a decrease in frequency as was observed in \cite{porterPredictivePhysicsbasedFriction2023}. 
The effect of excluding the mesoscale topology is less than was observed in \cite{porterPredictivePhysicsbasedFriction2023} due to the more flat mesoscale topology (see \Cref{fig:hbrb_surface}). 
Furthermore, excluding the mesoscale topology changes how the system transitions from sticking to slipping with a more graduate decrease in frequency starting at a lower amplitude level.

\FloatBarrier

\FloatBarrier

\begin{figure}[h!]
	\centering
	\includegraphics[width=0.55\linewidth]{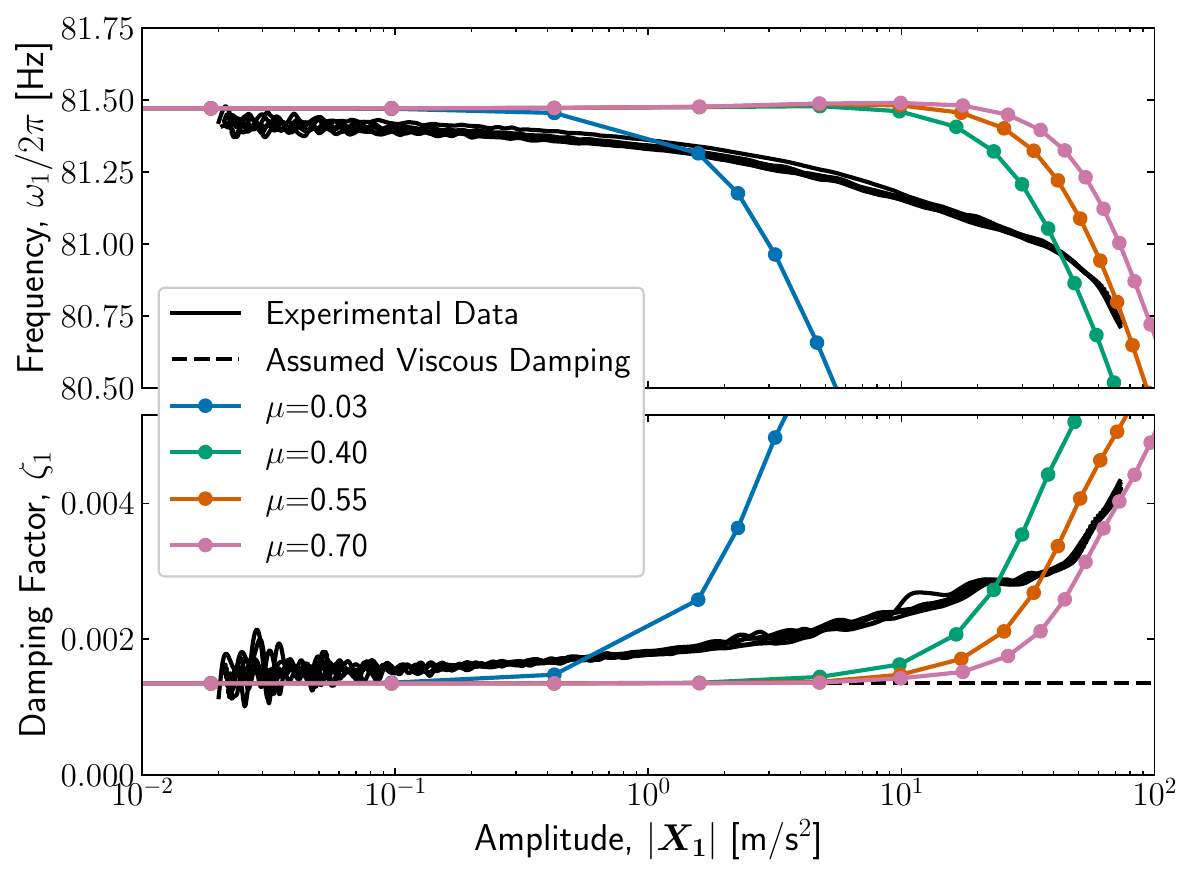}
	\caption{
		First bending mode EPMC results for different friction coefficients $ \mu $ compared to experiments using the lowest preload value and the tip accelerometer. A friction coefficient of $ \mu=0.55 $ is chosen for subsequent models.
	}
	\label{fig:hbrb_epmc_low}
\end{figure}

\begin{figure}[h!]
	\centering
	\includegraphics[width=0.55\linewidth]{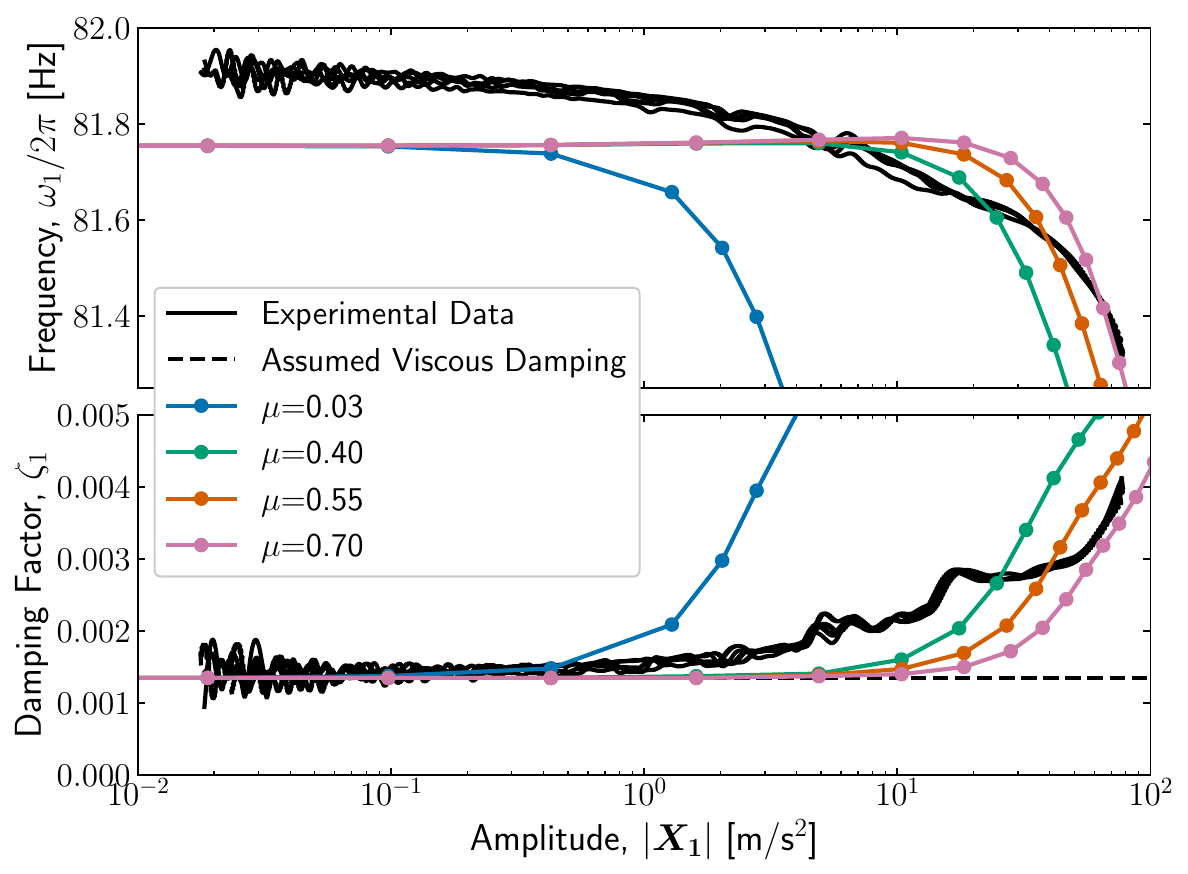}
	\caption{First bending mode EPMC results for different friction coefficients $ \mu $ compared to experiments using the medium preload value and the tip accelerometer. A friction coefficient of $ \mu=0.55 $ is chosen for subsequent models.}
	\label{fig:hbrb_epmc_medium}
\end{figure}

\begin{figure}[h!]
	\centering
	\includegraphics[width=0.55\linewidth]{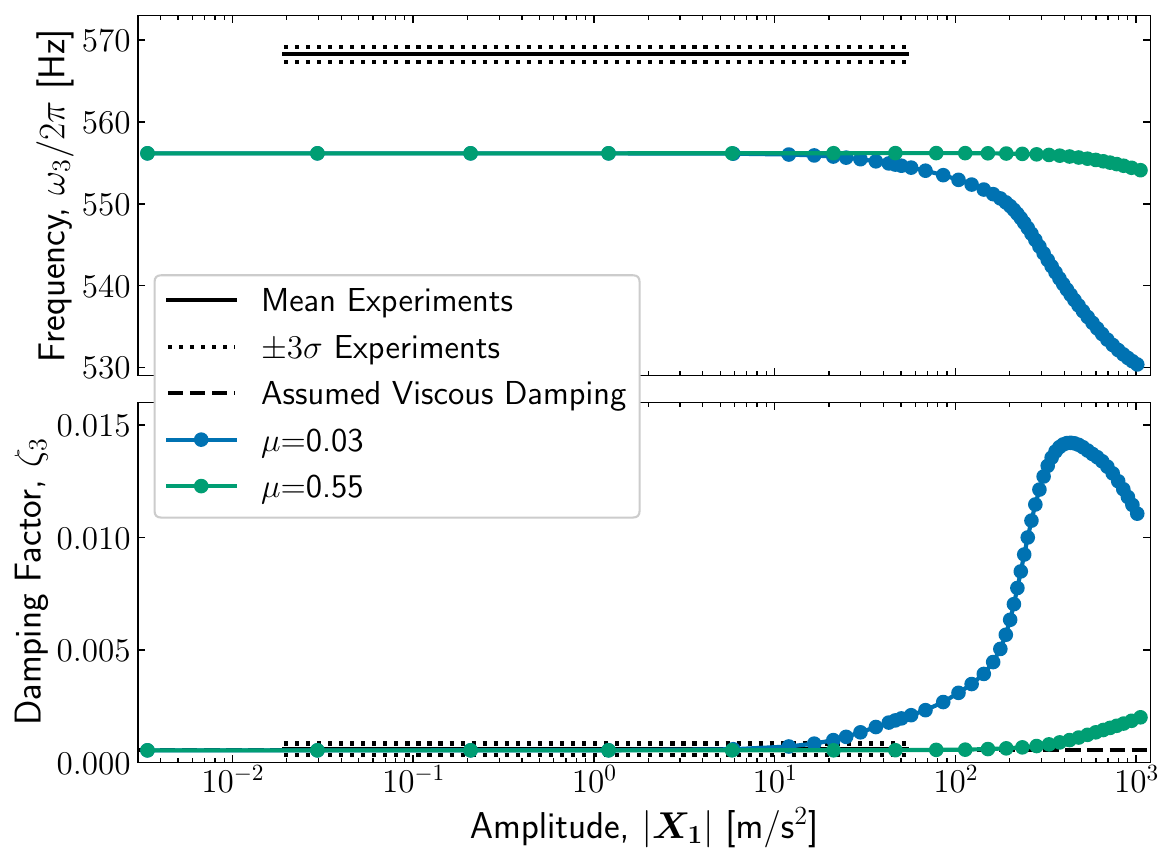}
	\caption{
		Third bending mode EPMC results for different friction coefficients $ \mu $ compared to experiments (with standard deviation $\sigma$) using the lowest preload value and the tip accelerometer. A friction coefficient of $ \mu=0.55 $ is chosen for subsequent models.	
		Since the linear damping for the model is chosen as the mean experimental damping, the mean experimental damping falls directly below the plotted models. }
	\label{fig:hbrb_epmc_third_low}
\end{figure}

\begin{figure}[h!]
	\centering
	\includegraphics[width=0.55\linewidth]{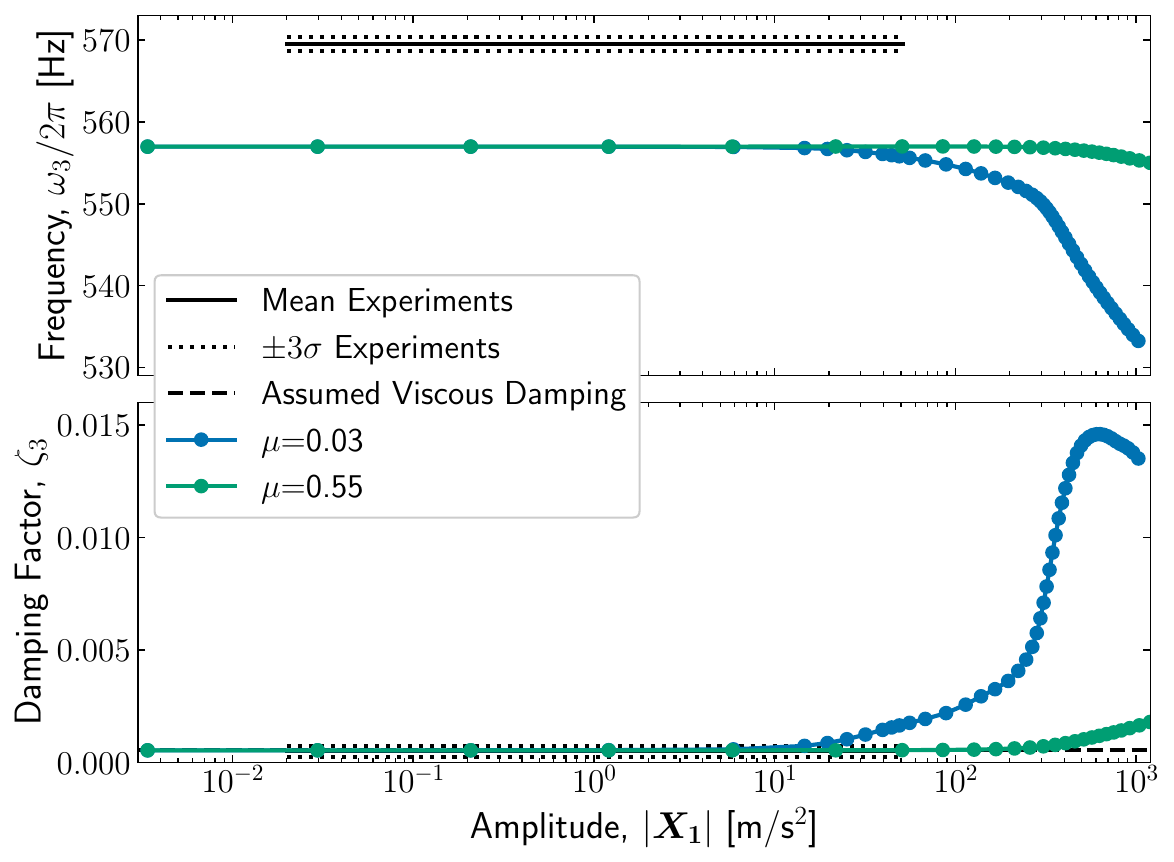}
	\caption{
		Third bending mode EPMC results for different friction coefficients $ \mu $ compared to experiments (with standard deviation $\sigma$) using the medium preload value and the tip accelerometer. A friction coefficient of $ \mu=0.55 $ is chosen for subsequent models.
		Since the linear damping for the model is chosen as the mean experimental damping, the mean experimental damping falls directly below the plotted models. }
	\label{fig:hbrb_epmc_third_med}
\end{figure}

EPMC results for the first bending mode at low (\Cref{fig:hbrb_epmc_low}) and medium (\Cref{fig:hbrb_epmc_medium}) preload are consistent with those at high preload in \Cref{sec:hbrb_epmc_results}.
Similar trends also are observed for the third bending mode at low (\Cref{fig:hbrb_epmc_third_low}) and medium (\Cref{fig:hbrb_epmc_third_med}) preload values as the high preload case (\Cref{sec:hbrb_epmc_results}).

\FloatBarrier

\FloatBarrier

\FloatBarrier

\subsection{Additional Stepped-Sine Test Results} \label{sec:extra_stepped_sine_appendix}

\begin{figure}[h!]
	\centering
	\begin{subfigure}[c]{\linewidth } 
		\centering
		\includegraphics[width=\linewidth]{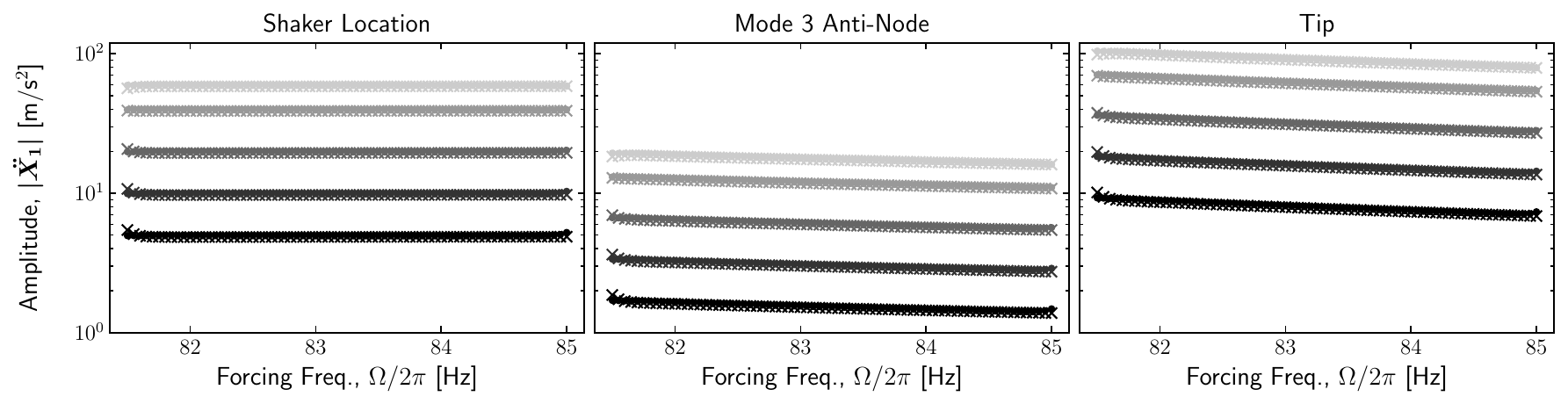}
		\caption{}
			\end{subfigure}
	\\
	\begin{subfigure}[c]{\linewidth } 
		\centering
		\includegraphics[width=\linewidth]{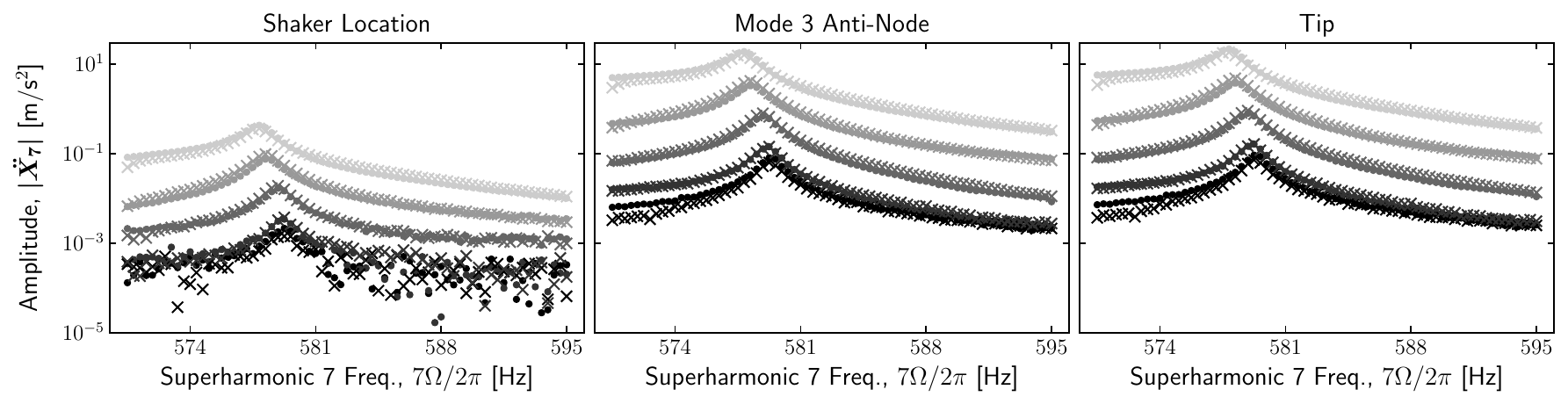}
		\caption{}
			\end{subfigure}
	\caption{
		HBRB experimental results for the low preload case stepped-sine response of (a) first harmonic and (b) seventh harmonic. 
		From dark to light, the tests are at amplitudes of 0.5 g, 1.0 g, 2.0 g, 4.0 g and 6.0 g for the first harmonic at the shaker.
		Tests are conducted sweeping from low to high frequency ($ \boldsymbol{\cdot} $) and high to low frequency ($ \boldsymbol{\times} $).
		The forcing frequency bounds are the same for (a) and (b), but (b) plots the superharmonic response frequency on the x-axis corresponding to seven times the forcing frequency.
	}
	\label{fig:exp_lowPre}
\end{figure}

\begin{figure}[h!]
	\centering
	\begin{subfigure}[b]{0.49\linewidth } 
		\centering
		\includegraphics[width=\linewidth]{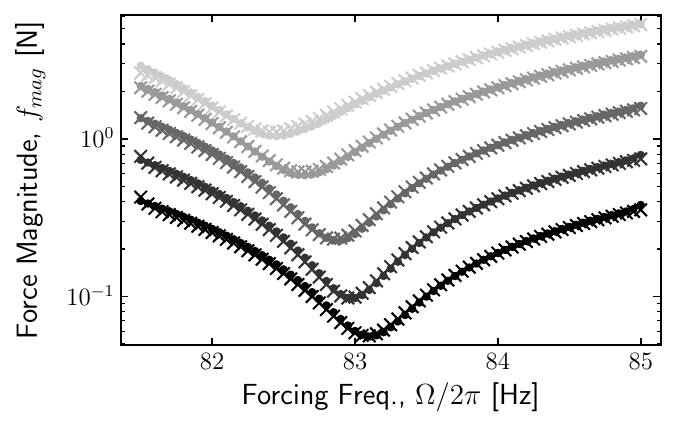}
		\caption{}
			\end{subfigure}
	\hfill
	\begin{subfigure}[b]{0.49\linewidth } 
		\centering
		\includegraphics[width=\linewidth]{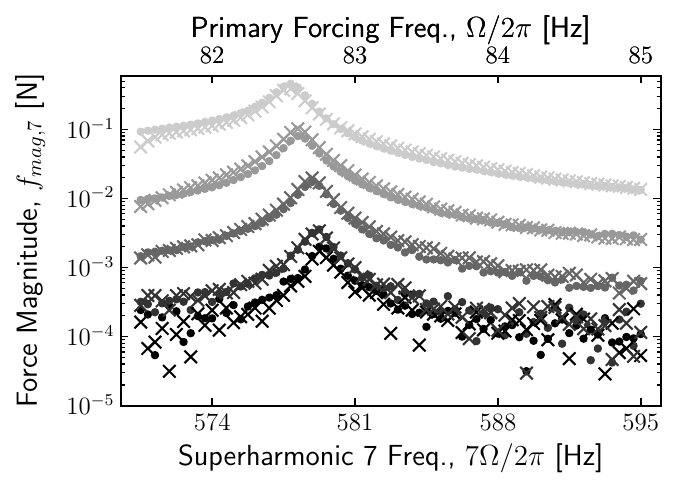}
		\caption{}
			\end{subfigure}
	\caption{
		HBRB experimental force measurements of (a) first harmonic and (b) seventh harmonic for the low preload case stepped-sine tests.
		From dark to light, the tests are at amplitudes of 0.5 g, 1.0 g, 2.0 g, 4.0 g and 6.0 g for the first harmonic at the shaker.
		Tests are conducted sweeping from low to high frequency ($ \boldsymbol{\cdot} $) and high to low frequency ($ \boldsymbol{\times} $).
		The forcing frequency bounds are the same for (a) and (b), but (b) plots the superharmonic response frequency on the x-axis corresponding to seven times the forcing frequency.
		The first harmonic force in (a) is the desired input to the structure while the seventh harmonic force in (b) is undesired structure-shaker interaction.
	}
	\label{fig:exp_lowPre_force}
\end{figure}

\begin{figure}[h!]
	\centering
	\begin{subfigure}[c]{\linewidth } 
		\centering
		\includegraphics[width=\linewidth]{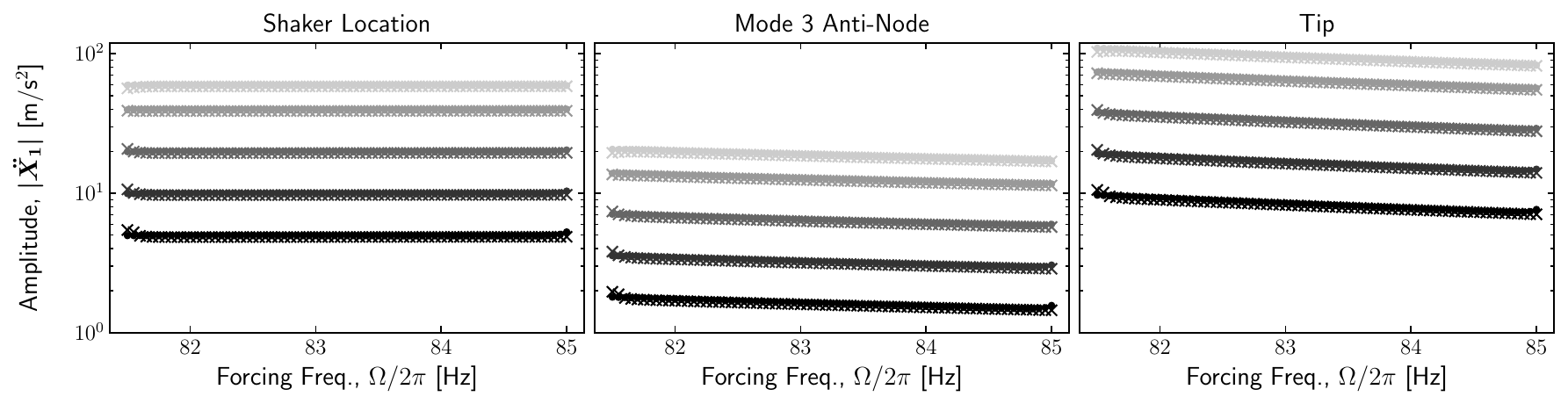}
		\caption{}
			\end{subfigure}
	\\
	\begin{subfigure}[c]{\linewidth } 
		\centering
		\includegraphics[width=\linewidth]{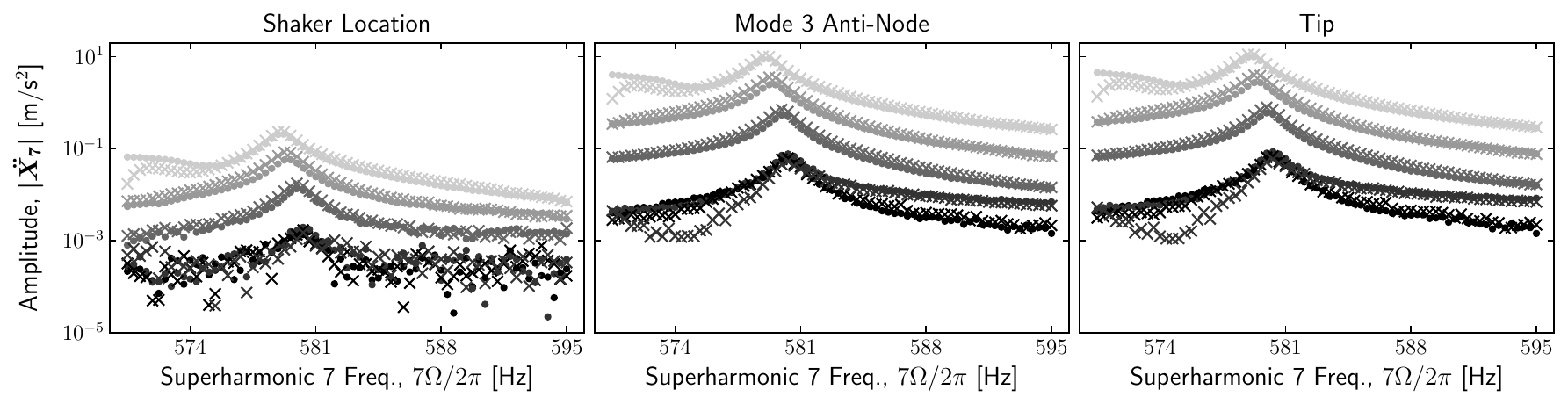}
		\caption{}
			\end{subfigure}
	\caption{HBRB experimental results for the medium preload case stepped-sine response of (a) first harmonic and (b) seventh harmonic. 
		From dark to light, the tests are at amplitudes of 0.5 g, 1.0 g, 2.0 g, 4.0 g and 6.0 g for the first harmonic at the shaker.
		Tests are conducted sweeping from low to high frequency ($ \boldsymbol{\cdot} $) and high to low frequency ($ \boldsymbol{\times} $).
		The x-axes are consistent between the top and bottom plot and are just scaled by the superharmonic number to reflect the response frequency.}
	\label{fig:exp_medPre}
\end{figure}

\begin{figure}[h!]
	\centering
	\begin{subfigure}[b]{0.49\linewidth } 
		\centering
		\includegraphics[width=\linewidth]{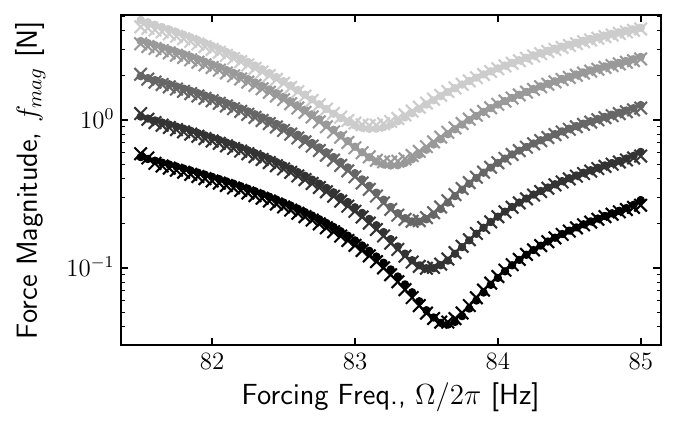}
		\caption{}
			\end{subfigure}
	\hfill
	\begin{subfigure}[b]{0.49\linewidth } 
		\centering
		\includegraphics[width=\linewidth]{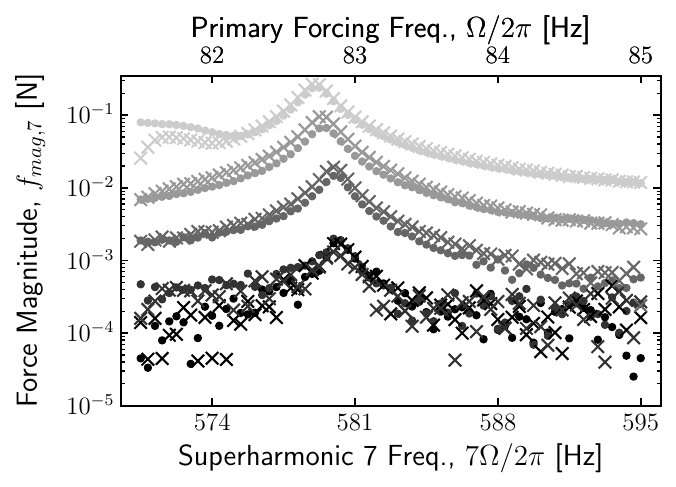}
		\caption{}
			\end{subfigure}
	\caption{
		HBRB experimental force measurements of (a) first harmonic and (b) seventh harmonic for the medium preload case stepped-sine tests.
		From dark to light, the tests are at amplitudes of 0.5 g, 1.0 g, 2.0 g, 4.0 g and 6.0 g for the first harmonic at the shaker.
		Tests are conducted sweeping from low to high frequency ($ \boldsymbol{\cdot} $) and high to low frequency ($ \boldsymbol{\times} $).
		The forcing frequency bounds are the same for (a) and (b), but (b) plots the superharmonic response frequency on the x-axis corresponding to seven times the forcing frequency.
		The first harmonic force in (a) is the desired input to the structure while the seventh harmonic force in (b) is undesired structure-shaker interaction.
	}
	\label{fig:exp_medPre_force}
\end{figure}

In addition to the experimental stepped-sine test results presented in \Cref{sec:exp_stepped_sine}, this section provides results for the low and medium preload levels. 
The responses for the first and seventh harmonics at low preload are shown in \Cref{fig:exp_lowPre}. For the low preload, the shaker applied first and seventh harmonic forces as shown in \Cref{fig:exp_lowPre_force} with the seventh harmonic forces being undesirable and occurring due to structure shaker interaction. 
Phase information for the superharmonic resonance at low preload level is shown in \Cref{fig:exp_lowPre_h7Phase}.
Likewise, the responses of the structure for the medium preload case are shown in \Cref{fig:exp_medPre} with shaker forces shown in \Cref{fig:exp_medPre_force} and superharmonic resonance phase information shown in \Cref{fig:exp_medPre_h7Phase}.

\FloatBarrier

\newcommand{\SteppedSinePhaseCaption}[1]{HBRB experimental phase of the seventh superharmonic resonance for the {#1} preload case stepped-sine tests. 
From dark to light, the tests are at amplitudes of 0.5 g, 1.0 g, 2.0 g, 4.0 g and 6.0 g for the first harmonic at the shaker.
Tests are conducted sweeping from low to high frequency ($ \boldsymbol{\cdot} $) and high to low frequency ($ \boldsymbol{\times} $).
Bottom is zoomed in version of top plot. }

\begin{figure}[h!]
	\centering
	\includegraphics[width=\linewidth]{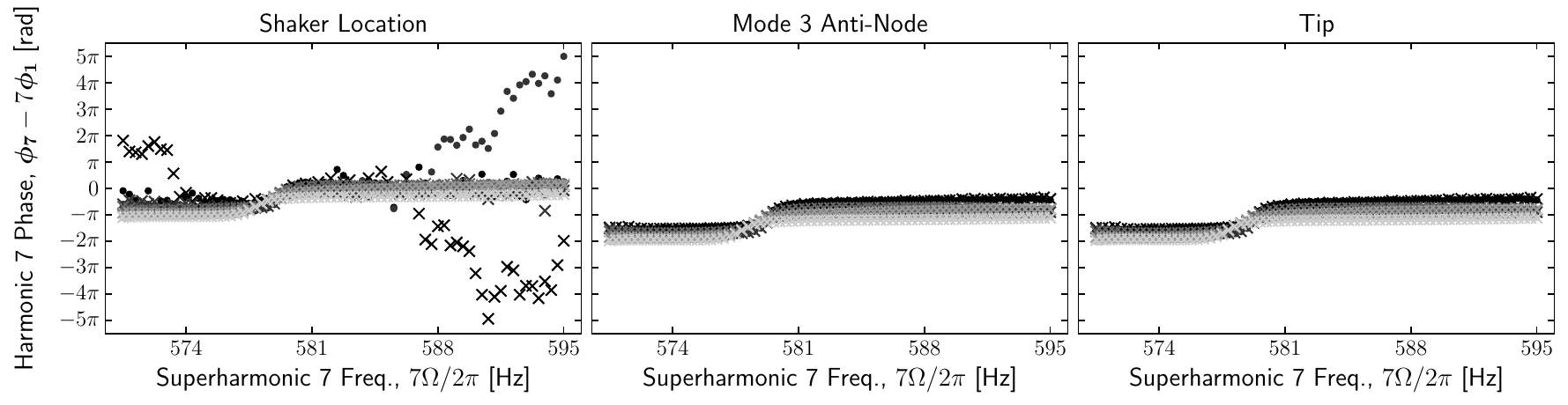}
	\\
	\includegraphics[width=\linewidth]{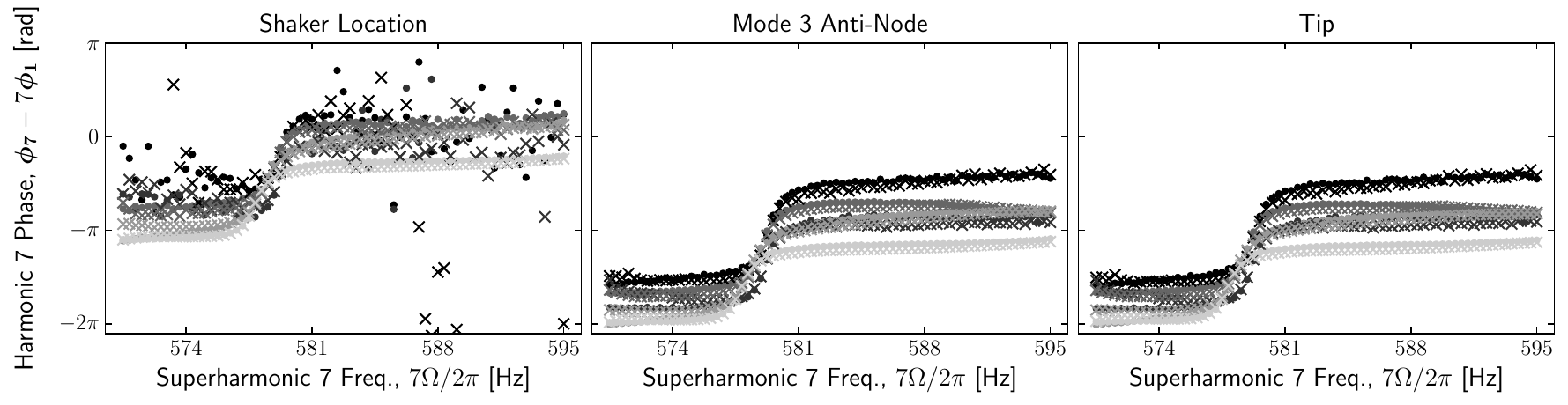}
	\caption{\SteppedSinePhaseCaption{low}
	}
	\label{fig:exp_lowPre_h7Phase}
\end{figure}

\FloatBarrier

\begin{figure}[h!]
	\centering
	\includegraphics[width=\linewidth]{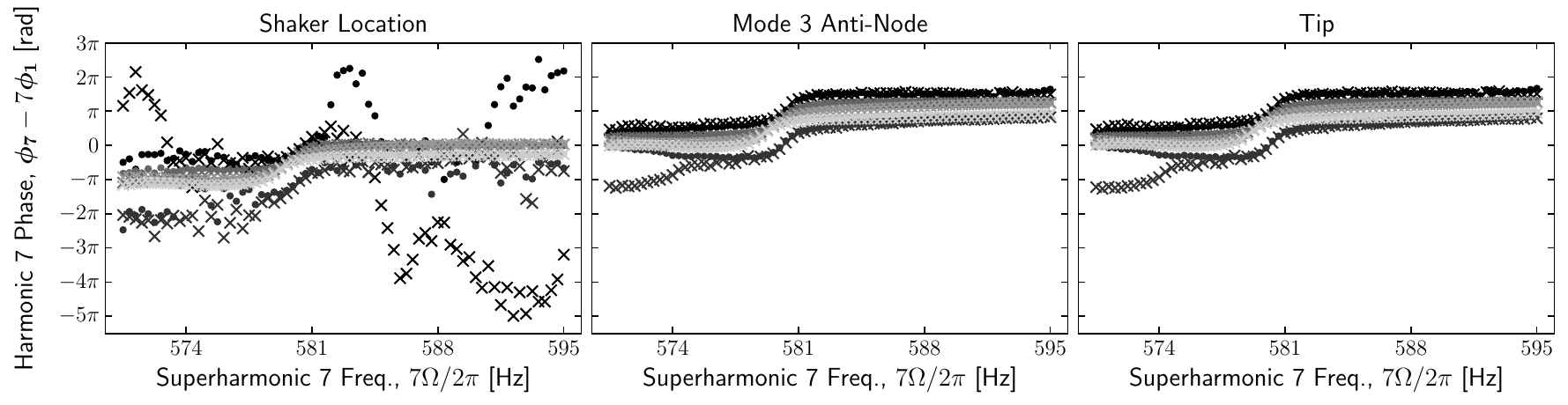}
	\\
	\includegraphics[width=\linewidth]{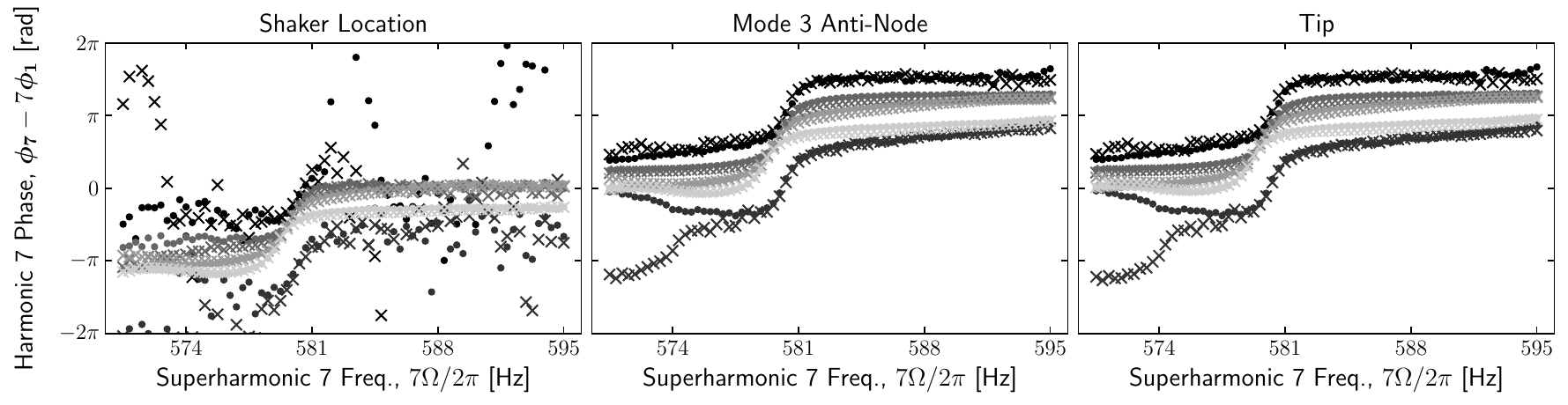}
	\caption{\SteppedSinePhaseCaption{medium}
	}
	\label{fig:exp_medPre_h7Phase}
\end{figure}

\FloatBarrier

\subsection{Harmonic Balance Method Convergence Results}\label{sec:hbm_convg}

At the lowest three acceleration levels (0.5 g, 1.0 g, 2.0 g), HBM results with harmonics 0 and 1-7 or 0 and 1-9 show converged behavior for low (\Cref{fig:hbm_lowPre_lowA_convg}), medium (\Cref{fig:hbm_medPre_lowA_convg}), and high (\Cref{fig:hbm_highPre_lowA_convg}) preload levels.
The solutions for the external force magnitudes have significant faceting because the amplitude and phase constraint HBM formulation described in \Cref{sec:theory_high_dim} result in the solution changing minimally over the primary resonance and large adaptive continuation steps are taken (reducing computation time). 
Since the primary resonance is not the main interest here, this reduction in computation time is considered an advantage. 
Also, the phase of the seventh harmonic is converged between using harmonics 0 and 1-7 or 0 and 1-9 at the lowest three acceleration levels and all three preload levels.
For the 4.0 g and 6.0 g acceleration cases, using harmonics 0 and 1-9 is consistent with harmonics 0 and 1-15 for low (\Cref{fig:hbm_lowPre_highA_convg}), medium (\Cref{fig:hbm_medPre_highA_convg}), and high (\Cref{fig:hbm_highPre_highA_convg}) preload values. For these cases, using harmonics 0 and 1-7 is not sufficient.

Initial simulations considering harmonics 0, 1, and 7 showed no distinguishable differences for using 64, 128, or 256 time steps in the AFT evaluations. Therefore, 128 time steps are used for AFT evaluations to be consistent with \cite{porterPredictivePhysicsbasedFriction2023}.

\FloatBarrier

\newcommand{\HBMlowAConvgCaption}[1]{{#1} preload HBM solutions for 0.5 g, 1.0 g, and 2.0 g cases show converged results using harmoncs 0 and 1-7 or 0 and 1-9. Results are for (a) external force magnitude to achieve amplitude control and (b) superharmonic resonance response amplitude.}

\begin{figure}[h!]
	\centering
	\begin{subfigure}[c]{0.25\linewidth } 
		\includegraphics[width=1\linewidth, trim={0cm 0.2cm 0cm 0}, clip]{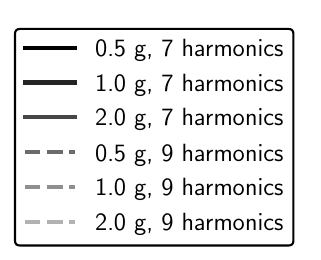}
	\end{subfigure}
	\quad\quad
	\begin{subfigure}[c]{0.4\linewidth } 
		\centering
		\includegraphics[width=\linewidth]{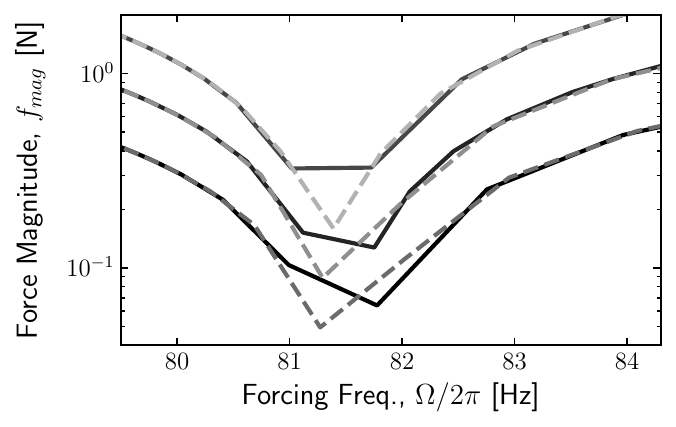}
		\caption{}
	\end{subfigure}
	\\
	\begin{subfigure}[c]{\linewidth } 
		\centering
		\includegraphics[width=\linewidth]{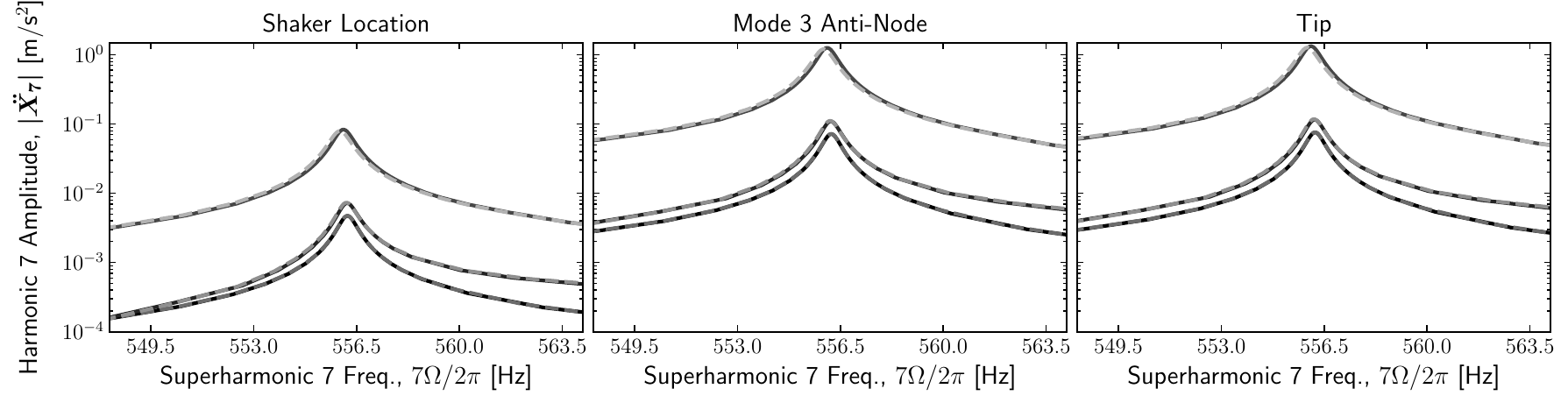}
		\caption{}
			\end{subfigure}
	\caption{\HBMlowAConvgCaption{Low}	}
	\label{fig:hbm_lowPre_lowA_convg}
\end{figure}

\begin{figure}[h!]
	\centering
	\begin{subfigure}[c]{0.25\linewidth } 
		\includegraphics[width=1\linewidth, trim={0cm 0.2cm 0cm 0}, clip]{hbm_highPre_lowA_convg_h7_legend.pdf}
	\end{subfigure}
	\quad\quad
	\begin{subfigure}[c]{0.4\linewidth } 
		\centering
		\includegraphics[width=\linewidth]{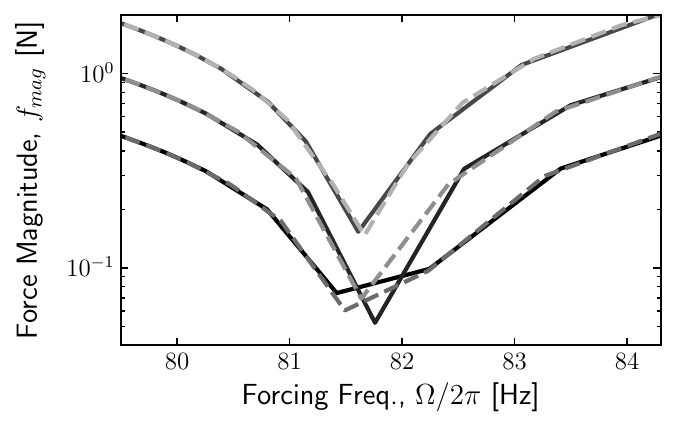}
		\caption{}
	\end{subfigure}
	\\
	\begin{subfigure}[c]{\linewidth } 
		\centering
		\includegraphics[width=\linewidth]{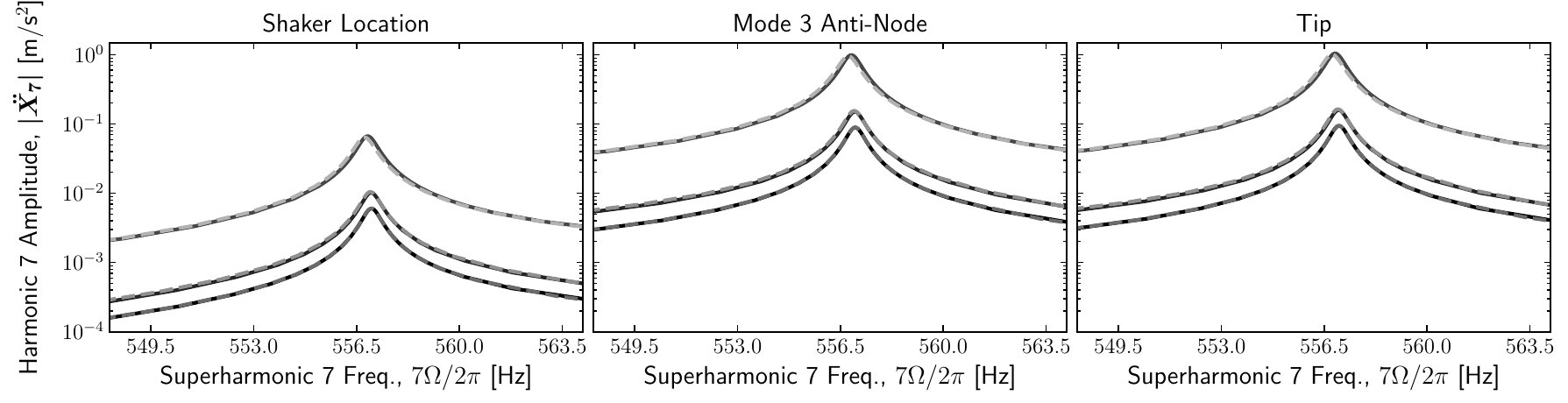}
		\caption{}
			\end{subfigure}
	\caption{\HBMlowAConvgCaption{Medium}
	}
	\label{fig:hbm_medPre_lowA_convg}
\end{figure}

\begin{figure}[h!]
	\centering
	\begin{subfigure}[c]{0.25\linewidth } 
		\includegraphics[width=1\linewidth, trim={0cm 0.2cm 0cm 0}, clip]{hbm_highPre_lowA_convg_h7_legend.pdf}
	\end{subfigure}
	\quad\quad
	\begin{subfigure}[c]{0.4\linewidth } 
		\centering
		\includegraphics[width=\linewidth]{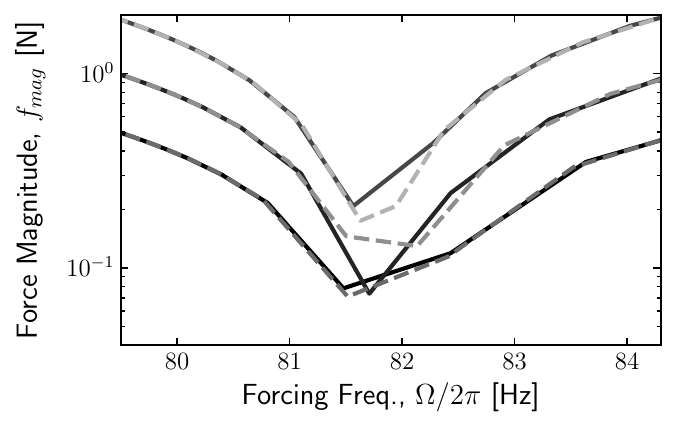}
		\caption{}
	\end{subfigure}
	\\
	\begin{subfigure}[c]{\linewidth } 
		\centering
		\includegraphics[width=\linewidth]{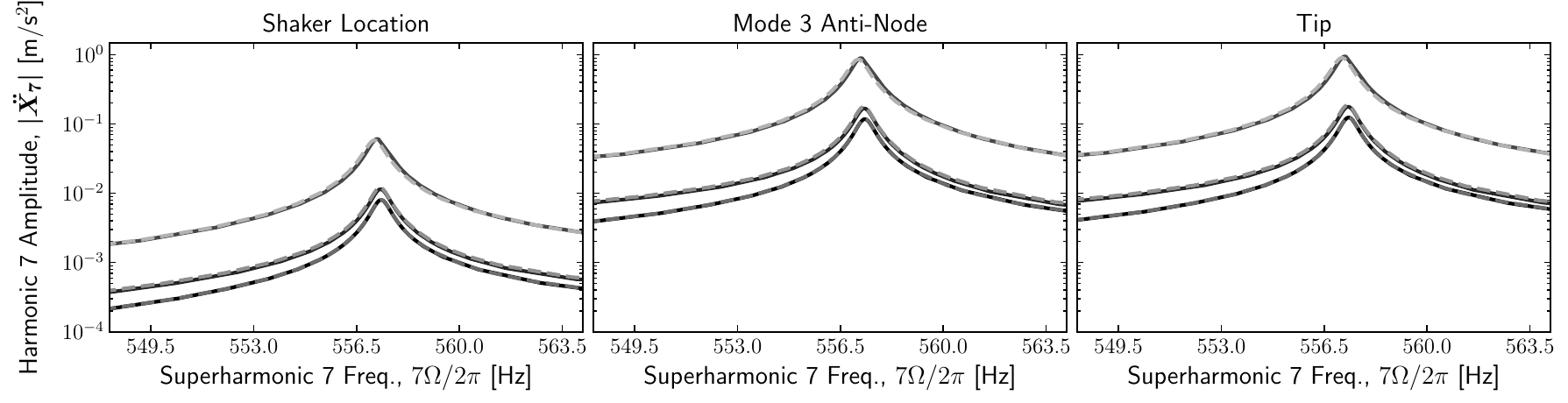}
		\caption{}
			\end{subfigure}
	\caption{\HBMlowAConvgCaption{High}
	}
	\label{fig:hbm_highPre_lowA_convg}
\end{figure}

\FloatBarrier

\newcommand{\HBMhighAConvgCaption}[1]{{#1} preload HBM solutions for 4.0 g and 6.0 g cases show converged results using harmoncs 0 and 1-9 or 0 and 1-15. Results are for (a) external force magnitude to achieve amplitude control and (b) superharmonic resonance response amplitude.}

\begin{figure}[h!]
	\centering
	\begin{subfigure}[c]{0.25\linewidth } 
		\includegraphics[width=1\linewidth, trim={0cm 0.2cm 0cm 0}, clip]{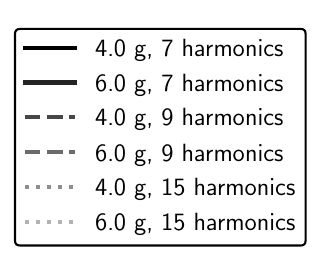}
	\end{subfigure}
	\quad\quad
		\begin{subfigure}[c]{0.4\linewidth } 
		\centering
		\includegraphics[width=\linewidth]{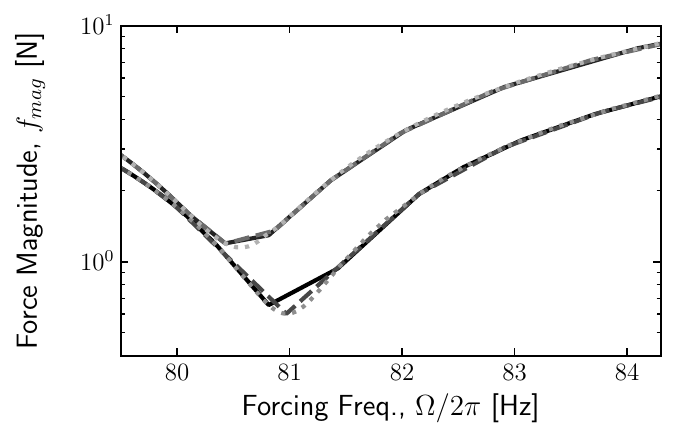}
		\caption{}
			\end{subfigure}
	\\
	\begin{subfigure}[c]{\linewidth } 
		\centering
		\includegraphics[width=\linewidth]{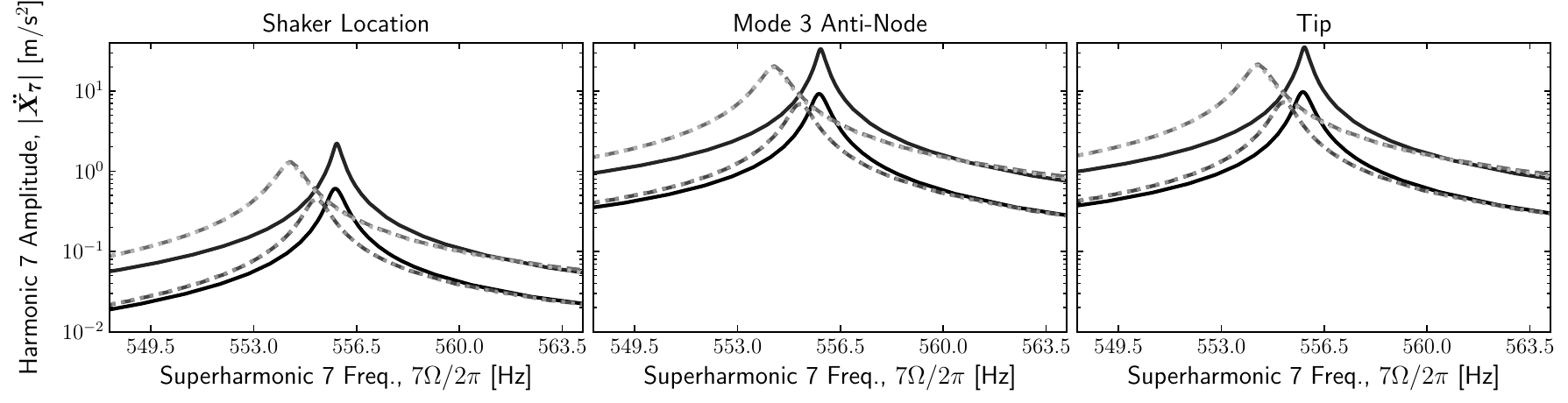}
		\caption{}
			\end{subfigure}
	\caption{\HBMhighAConvgCaption{Low}
	}
	\label{fig:hbm_lowPre_highA_convg}
\end{figure}

\begin{figure}[h!]
	\centering
	\begin{subfigure}[c]{0.25\linewidth } 
		\includegraphics[width=1\linewidth, trim={0cm 0.2cm 0cm 0}, clip]{hbm_highPre_highA_convg_h7_legend.pdf}
	\end{subfigure}
	\quad\quad
		\begin{subfigure}[c]{0.4\linewidth } 
		\centering
		\includegraphics[width=\linewidth]{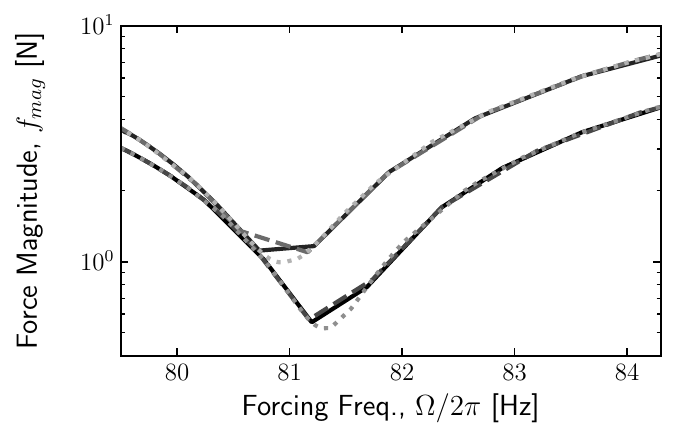}
		\caption{}
			\end{subfigure}
	\\
	\begin{subfigure}[c]{\linewidth } 
		\centering
		\includegraphics[width=\linewidth]{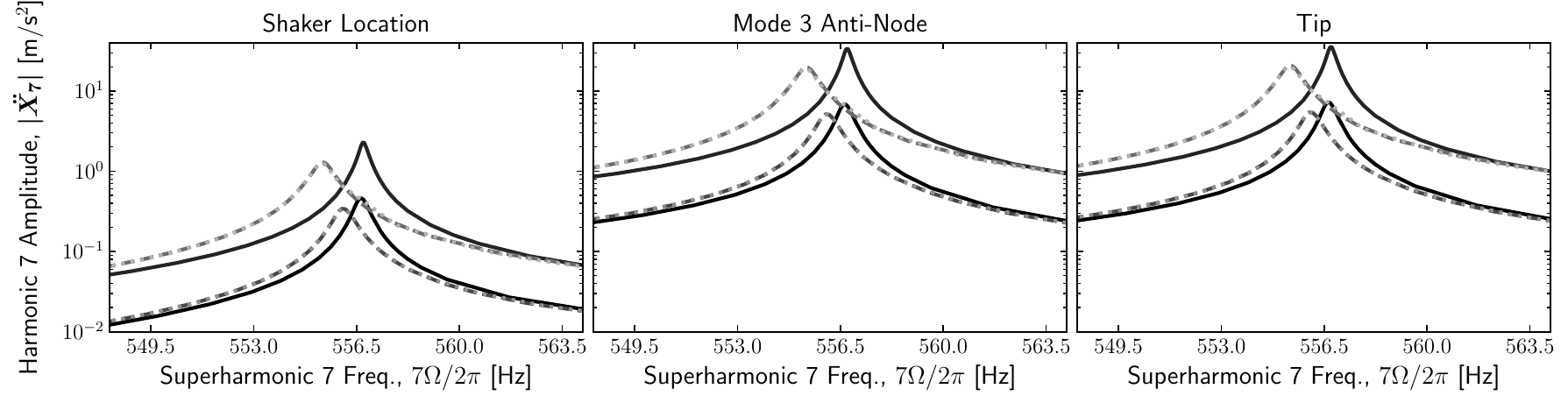}
		\caption{}
			\end{subfigure}
	\caption{\HBMhighAConvgCaption{Medium}
	}
	\label{fig:hbm_medPre_highA_convg}
\end{figure}

\begin{figure}[h!]
	\centering
	\begin{subfigure}[c]{0.25\linewidth } 
		\includegraphics[width=1\linewidth, trim={0cm 0.2cm 0cm 0}, clip]{hbm_highPre_highA_convg_h7_legend.pdf}
	\end{subfigure}
	\quad\quad
		\begin{subfigure}[c]{0.4\linewidth } 
		\centering
		\includegraphics[width=\linewidth]{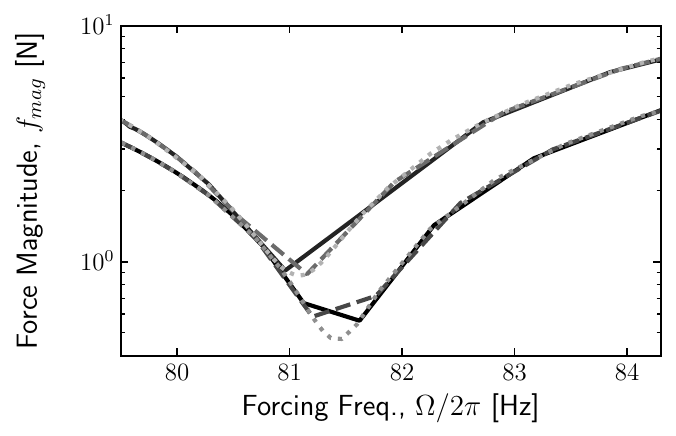}
		\caption{}
		\label{fig:hbm_highPre_highA_convg_force}
	\end{subfigure}
	\\
	\begin{subfigure}[c]{\linewidth } 
		\centering
		\includegraphics[width=\linewidth]{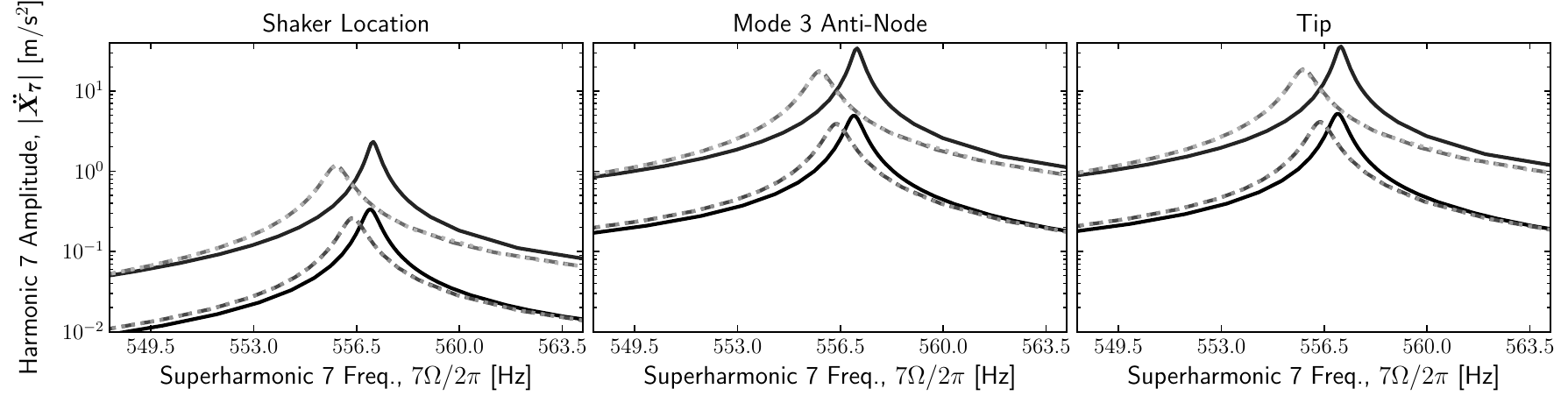}
		\caption{}
			\end{subfigure}
	\caption{\HBMhighAConvgCaption{High}
	}
	\label{fig:hbm_highPre_highA_convg}
\end{figure}

\FloatBarrier

\subsection{Additional Harmonic Balance and VPRNM ROM Results} \label{sec:extra_hbm_vprnm}

\begin{figure}[h!]
	\centering
	\includegraphics[width=0.85\linewidth]{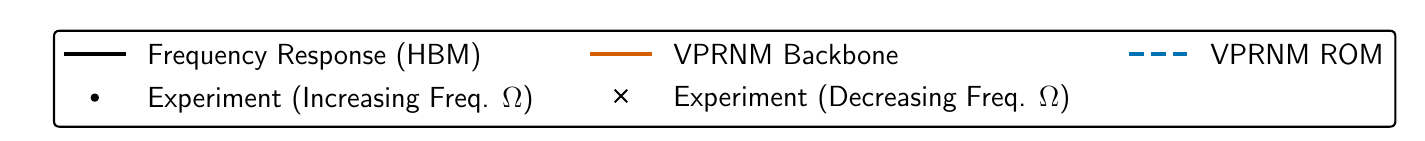}
	\begin{subfigure}[c]{\linewidth } 
		\centering
		\includegraphics[width=\linewidth]{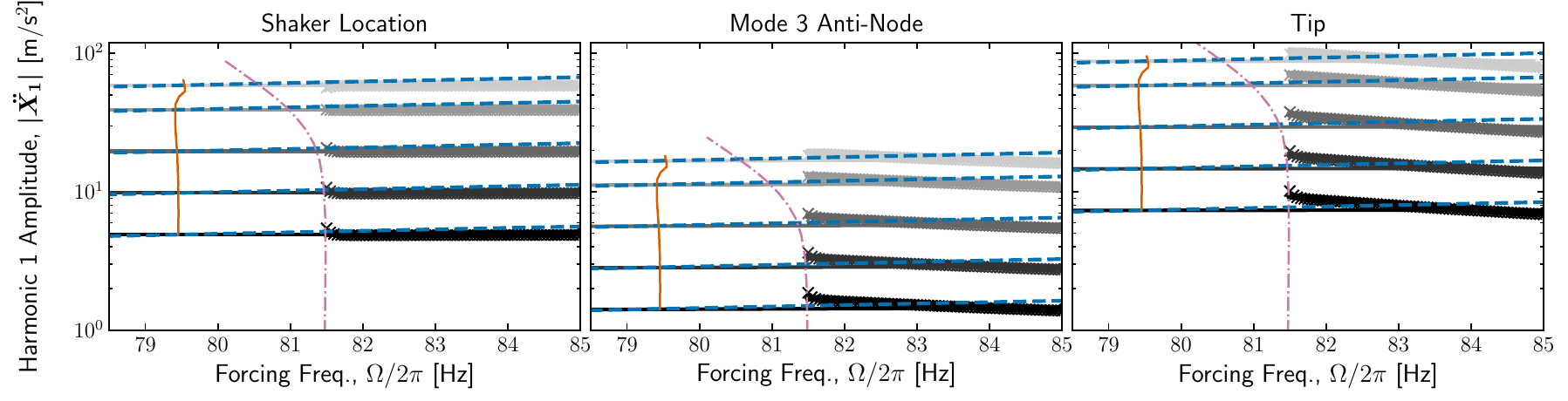}
		\caption{}
			\end{subfigure}
	\\
	\begin{subfigure}[c]{\linewidth } 
		\centering
		\includegraphics[width=\linewidth]{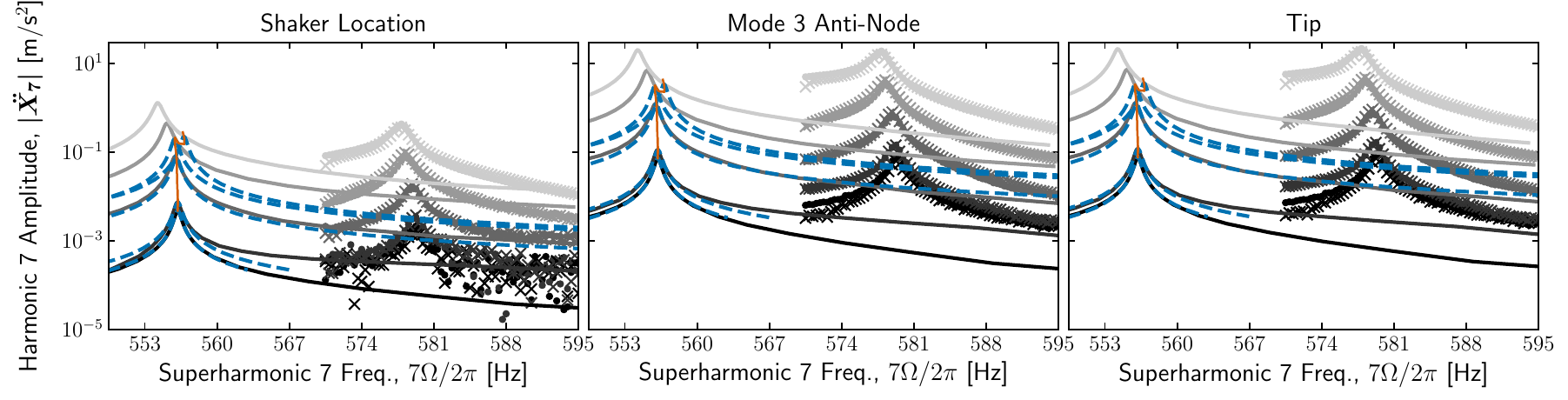}
		\caption{}
		\label{fig:rom_lowPre_h7}
	\end{subfigure}
	\caption{\HbrbAmpCaption{low}
											}
	\label{fig:rom_lowPre}
\end{figure}

\begin{figure}[h!]
	\centering
	\includegraphics[width=0.85\linewidth]{hbrb_all_legend.pdf}
	\includegraphics[width=0.49\linewidth]{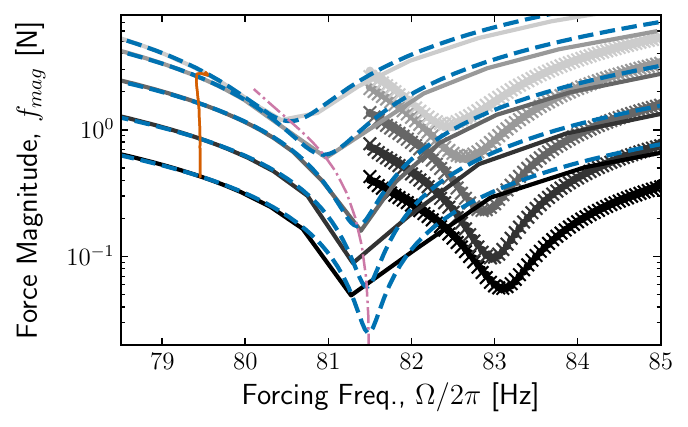}
	\caption{\HbrbForceCaption{low}{\Cref{fig:rom_lowPre}}
							}
	\label{fig:rom_lowPre_force}
\end{figure}

\begin{figure}[h!]
	\centering
	\includegraphics[width=0.65\linewidth]{hbrb_rom_no_bb_legend.pdf}
	\includegraphics[width=\linewidth]{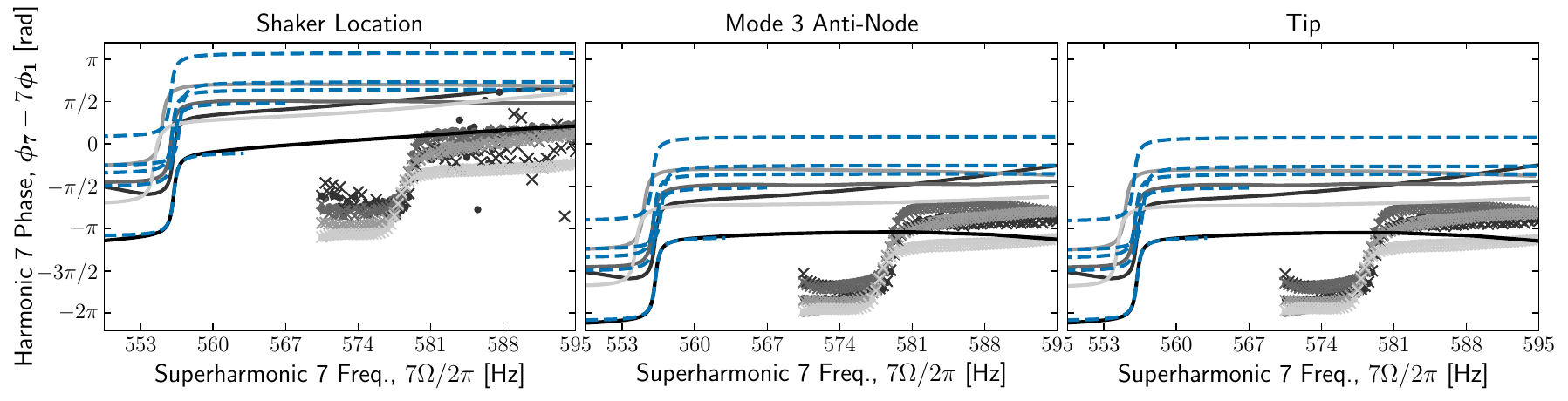}
	\caption{
		\HbrbPhaseSevenCaption{low}{}
									}
	\label{fig:rom_lowPre_h7phase}
\end{figure}

\begin{figure}[h!]
	\centering
	\includegraphics[width=0.85\linewidth]{hbrb_rom_legend.pdf}
	\begin{subfigure}[c]{\linewidth } 
		\centering
		\includegraphics[width=\linewidth]{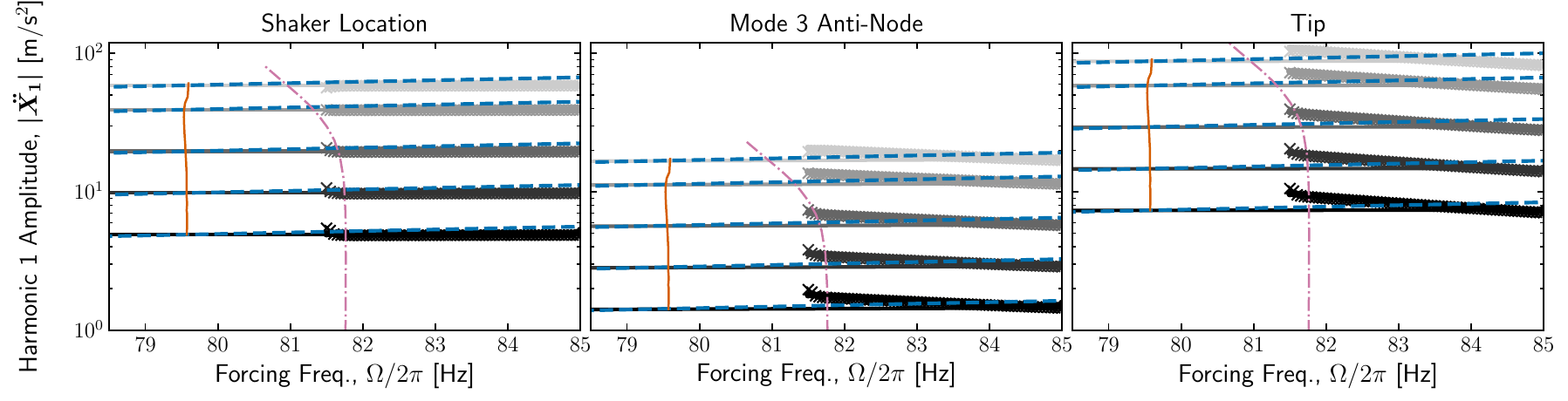}
		\caption{}
			\end{subfigure}
	\\
	\begin{subfigure}[c]{\linewidth } 
		\centering
		\includegraphics[width=\linewidth]{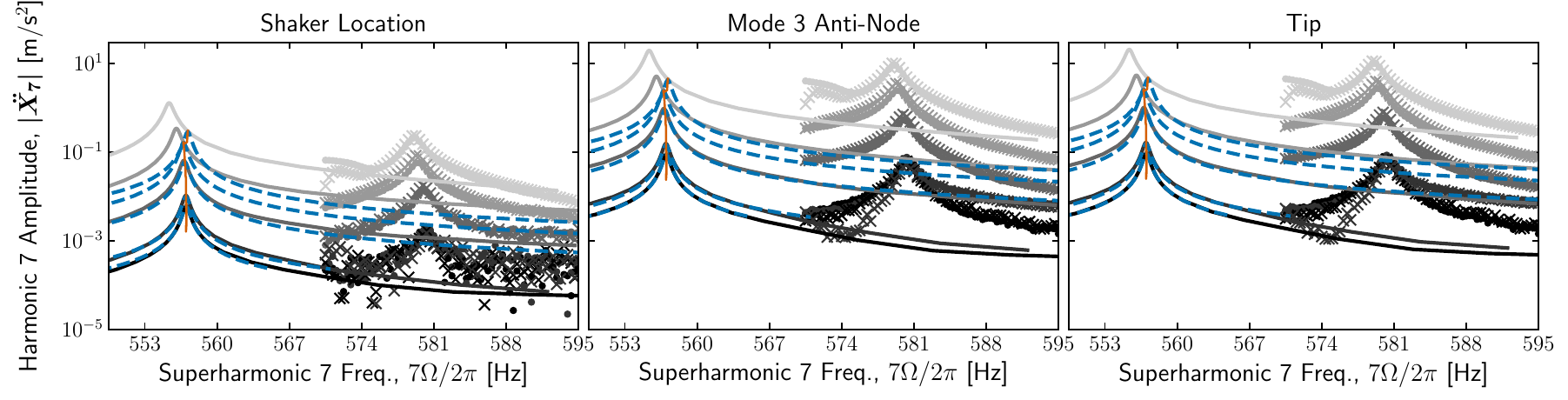}
		\caption{}
		\label{fig:rom_mediumPre_h7}
	\end{subfigure}
	\caption{\HbrbAmpCaption{medium}
											}
	\label{fig:rom_mediumPre}
\end{figure}

\begin{figure}[h!]
	\centering
	\includegraphics[width=0.85\linewidth]{hbrb_all_legend.pdf}
	\includegraphics[width=0.49\linewidth]{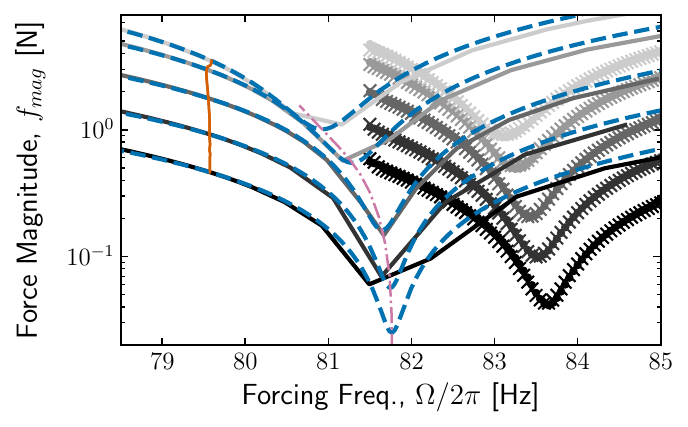}
	\caption{\HbrbForceCaption{medium}{\Cref{fig:rom_mediumPre}}
							}
	\label{fig:rom_mediumPre_force}
\end{figure}

\begin{figure}[h!]
	\centering
	\includegraphics[width=0.65\linewidth]{hbrb_rom_no_bb_legend.pdf}
	\includegraphics[width=\linewidth]{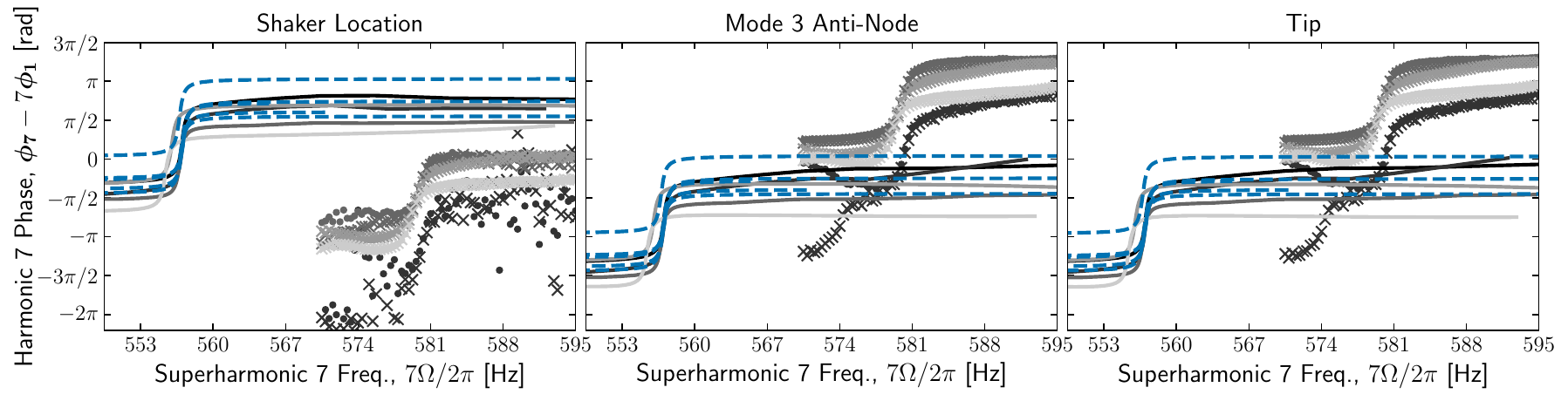}
	\caption{\HbrbPhaseSevenCaption{medium}{\PhaseClipping}
									}
	\label{fig:rom_mediumPre_h7phase}
\end{figure}

This section presents HBM and VPRNM ROM results for low and medium preload to supplement the high preload results presented in \Cref{sec:hbrb_hbm_res}. As with the high preload case, the low and medium preload cases show good qualitative agreement between experiments and models. The VPRNM ROM captures the behavior from the HBM results for the 0.5 g, 1.0 g, and 2.0 g cases, but shows some errors for the 4.0 g and 6.0 g cases. 
For low preload, \Cref{fig:rom_lowPre} shows the fundamental and superharmonic resonance responses; \Cref{fig:rom_lowPre_force} shows the force magnitude to achieve control; \Cref{fig:rom_lowPre_h7phase} shows the phase of the seventh harmonic.
For medium preload, the response amplitudes (\Cref{fig:rom_mediumPre}), forcing magnitudes (\Cref{fig:rom_mediumPre_force}) and seventh harmonic phases (\Cref{fig:rom_mediumPre_h7phase}) are likewise similar to the high preload cases in \Cref{sec:hbrb_hbm_res}.

\FloatBarrier

\FloatBarrier

\FloatBarrier

\FloatBarrier

\FloatBarrier

\FloatBarrier

\end{document}
\typeout{get arXiv to do 4 passes: Label(s) may have changed. Rerun}

%% file: preamble.tex
\Large
\noindent \textbf{Efficient Model Reduction and Prediction of Superharmonic Resonances in Frictional and Hysteretic Systems}

\vspace{24 pt}

\large

\noindent Justin H. Porter and Matthew R. W. Brake

\noindent Department of Mechanical Engineering, Rice University, Houston,
TX 77005

\vspace{1 in}

\normalsize

\noindent This is a preprint of the journal paper 
``Efficient Model Reduction and Prediction of Superharmonic Resonances in Frictional and Hysteretic Systems.'' This paper is currently under review at the journal Mechanical Systems and Signal Processing.
\\
\\\\

\thispagestyle{empty}